\documentclass[a4paper,14pt]{extarticle} 
  
\newcommand{\mbf}[1]{\protect\makebox{\boldmath$#1$}}
\newcommand{\eus}{\EuScript}
\newcommand{\unitint}{[ \, -1, \, 1 \, ]}

\usepackage[left=3cm,right=1.5cm,top=2cm,bottom=2cm,bindingoffset=0cm]{geometry}
\usepackage{setspace} 

\usepackage{cmap}					
\usepackage[T2A]{fontenc}			
\usepackage[utf8]{inputenc}			
\usepackage[english,russian]{babel}	

\usepackage[margin=10pt,font=small,labelfont=bf,center]{caption}
\DeclareCaptionLabelSeparator{dot}{. }
\newenvironment{definition}%
{\par\addvspace{\medskipamount}\addtocounter{DefNum}{1} 
\noindent\textbf{Определение {\thesection}.\arabic{DefNum}.}}%
{\par\medskip} 

\usepackage{geometry} 
\usepackage{mathtext} 
\usepackage{euscript} 
\usepackage{mathrsfs} 
\usepackage{graphicx}
\usepackage{amsmath,amsthm,amsfonts,amssymb} 
\usepackage{icomma} 

\usepackage{array,tabularx,tabulary,booktabs} 
\usepackage{longtable}  
\usepackage{multirow} 
\usepackage{hhline}

\usepackage{graphicx}  
\graphicspath{{images/}}  
\usepackage{wrapfig} 
\usepackage{epstopdf} 

\usepackage{lscape} 

\usepackage[usenames]{color}
\usepackage{color}
\usepackage{colortbl}

\makeatletter
\renewcommand{\@biblabel}[1]{#1.}
\makeatother

\usepackage{tocloft}

\usepackage{titlesec} 
\titlelabel{\thetitle.\quad} 

\usepackage{indentfirst}

\numberwithin{equation}{section}
\numberwithin{figure}{section}
\numberwithin{table}{section}

\usepackage{soul}
\usepackage{csquotes}

\usepackage{pdfpages}

\begin{document}

\newcommand{\tocsecindent}{\hspace{7mm}}

\begin{titlepage}

\begin{center}
{\small МИНИСТЕРСТВО НАУКИ И ВЫСШЕГО ОБРАЗОВАНИЯ \\РОССИЙСКОЙ ФЕДЕРАЦИИ} \\
\vspace{1em}
{\small ФЕДЕРАЛЬНОЕ ГОСУДАРСТВЕННОЕ АВТОНОМНОЕ ОБРАЗОВАТЕЛЬНОЕ УЧРЕЖДЕНИЕ ВЫСШЕГО ОБРАЗОВАНИЯ}\\
\vspace{1em}
{\small «НОВОСИБИРСКИЙ НАЦИОНАЛЬНЫЙ ИССЛЕДОВАТЕЛЬСКИЙ ГОСУДАРСТВЕННЫЙ УНИВЕРСИТЕТ» (НОВОСИБИРСКИЙ ГОСУДАРСТВЕННЫЙ УНИВЕРСИТЕТ, НГУ)} \\
\end{center}

\vspace{1.5em}
\raggedright{\small Факультет \qquad \qquad Механико-математический факультет}\\*
\vspace{0.8em}
\raggedright{\small Кафедра \; \; \; \quad \qquad Математического моделирования} \\*
\vspace{1.5em}
\raggedright{\small Направление подготовки \qquad Прикладная математика и информатика}\\*
\vspace{1.5em}

\begin{center}
{\textbf{\small ВЫПУСКНАЯ КВАЛИФИКАЦИОННАЯ РАБОТА МАГИСТРА}} \\
\vspace{0.8em}
{\small Скорика Дмитрия Александровича}\\
\end{center}

\vspace{0.6em}
\begin{center}
{\small Тема работы: Развитие линейной функциональной арифметики и её приложение \\к решению задач интервального анализа.}\\*
\end{center}
\vspace{2cm}

\begin{small}
	\begin{tabular}{m{8.5cm} l}
	
		\textbf{<<К защите допущен>>} & \textbf{Научный руководитель} \\ \vspace{1mm}
		
		Заведующий кафедрой, & Д. ф.-м. н., \\ \vspace{1mm}
		
		д. т. н., доц. & проф. каф. ММод ММФ НГУ \\ \vspace{1mm}
		
		\underline{ Барахнин В. Б. } / \underline{\phantom{ПодписьПодпись}} & \underline{ Шарый С. П. } / \underline{\phantom{ПодписьПодпись}} \\
		
		\footnotesize{(Фамилия И. О.)} \hspace{1cm} \footnotesize{(Подпись)} & \footnotesize{(Фамилия И. О.)} \hspace{1cm} \footnotesize{(Подпись)} \\ \vspace{1mm}
		
		<<\makebox[0.7cm][l]{\dotfill}>>\makebox[5cm][l]{\dotfill 2022 г.} & <<\makebox[0.7cm][l]{\dotfill}>>\makebox[5cm][l]{\dotfill 2022 г.}

	\end{tabular}
\end{small}

\vspace{3.2cm}
\raggedright{\small \ \qquad \qquad \qquad \qquad \qquad \qquad \qquad \qquad Дата защиты:  \makebox[0.07\textwidth][l]{«\dotfill»}\makebox[0.3\textwidth][l]{\dotfill 2022 г.}} \\
\vspace{0.5cm}
\begin{center}
{\small Новосибирск, 2022}
\end{center}

\thispagestyle{empty}
\end{titlepage}
\clearpage


\setcounter{page}{2}

\clearpage
\tableofcontents

\newpage
\newcounter{DefNum}[subsection]
\newcounter{LemmNum}[subsection]
\setcounter{DefNum}{0} 
\newtheorem{lemma}[theorem]{Лемма}

\section*{Реферат}
\addcontentsline{toc}{section}{Реферат}

\textbf{Тема работы:} Развитие линейной функциональной арифметики и её приложение к решению задачи интервального анализа.

\textbf{Объём работы} составляет 63 страниц, список использованной литературы включает в себя 11 источников, в работе приводятся 4 таблицы и 42 рисунка.

\textbf{Ключевые слова:} интервальный анализ, интервал, функциональный интервал, интервальная арифметика

Работа посвящена построению нового вида интервалов --- функциональных интервалов. Эти интервалы построены на идее расширения границ с чисел на функции. Функциональные интервалы показали себя перспективными для дальнейшего изучения и использования, поскольку имеют более богатые алгебраические свойства по сравнению с классическими интервалами.

В работе была построена линейная функциональная арифметика от одной переменной. Данная арифметика была применена для решения таких задач интервального анализа, как минимизация функции на интервале и нахождение нулей функции на интервале.

Результаты численных экспериментов для линейной функциональной арифметики показали высокий порядок сходимости и более высокую скорость работы алгоритмов при использовании интервалов нового вида, несмотря на то, что при вычислениях не использовалась информация о производной функции.

Также в работе исследовалась модификация алгоритмов минимизации функции от нескольких переменных, основанная на использовании функциональных интервалов от нескольких переменных. В результате было получено ускорение алгоритмов, но только до некоторого числа неизвестных.

\clearpage

\section{Введение}

\subsection{О числовых неопределённостях}

При решении практических задач часто возникают численные неопределённости, имеющие различную 
природу. Одни, например, происходят из неточности измерительных приборов или человеческого 
фактора. Другие~--- из-за использования на ЭВМ арифметики с ограниченной точностью 
(или арифметики с плавающей точкой).

Эти неопределённости могут достаточно сильно отклонить результат вычислений от истинного. 
Также существуют неопределённости, которые мы не можем устранить. Примерами могут служить 
атомные массы химических элементов, для значений которых известны лишь доверительные 
интервалы, гарантированно содержащих их. 
  
Одним из естественных подходов описания числовых неопределённостей, которые принимают 
значения в каком-либо промежутке, является указание интервала, который гарантированно 
содержит в себе все возможные значения. Интервалы просты в описании и использовании. 

Поэтому для решения задач, в которых возникают параметры с неопределённостями, возникает 
необходимость математического аппарата, который позволит оперировать интервальными 
величинами точно так же, как и <<точечными>>. Таким аппаратом является \textit{интервальный 
анализ} \cite{Shary}, \cite{HansenWalster}, \cite{AlefeldHerzberger}.

Одним из основных инструментов интервального анализа является классическая интервальная 
арифметика. Впервые этот инструмент анализа появился в середине 20-го столетия и был рассмотрен 
в работе японского ученого Теруо Сунаги в Японии \cite{Sunaga}. 
  
Зачастую бывает, что использование классической интервальной арифметики в вычислениях 
недостаточно, поскольку она дает достаточно грубые интервальные оценки при больших по объёму 
вычислениях. К тому же она имеет скудные алгебраические свойства \cite{AlefeldHerzberger}, 
\cite{Shary}, \cite{Skorik}. Это является мотивацией к разработке иных видов интервалов 
и интервальных  арифметик, которые позволят улучшить их свойства или приобрести новые. 
  
Цель данной работы --- построение нового вида интервалов, так называемых функциональных 
интервалов. Подробно будут рассмотрены объекты их простейшего подсемейства --- 
однопараметрические  линейные функциональные интервалы. Понятие функционального интервала 
будет введено далее. 
  
  
\subsection[Мотивации использования функциональных интервалов]%
           {Мотивации использования\\  функциональных интервалов}

Рассмотрим какую-либо функцию $f(x)$ на интервале $[ \, a, \, b \, ]$:

\begin{figure}[ht]
	\begin{center}
		
		\includegraphics[width = 0.40 \linewidth]{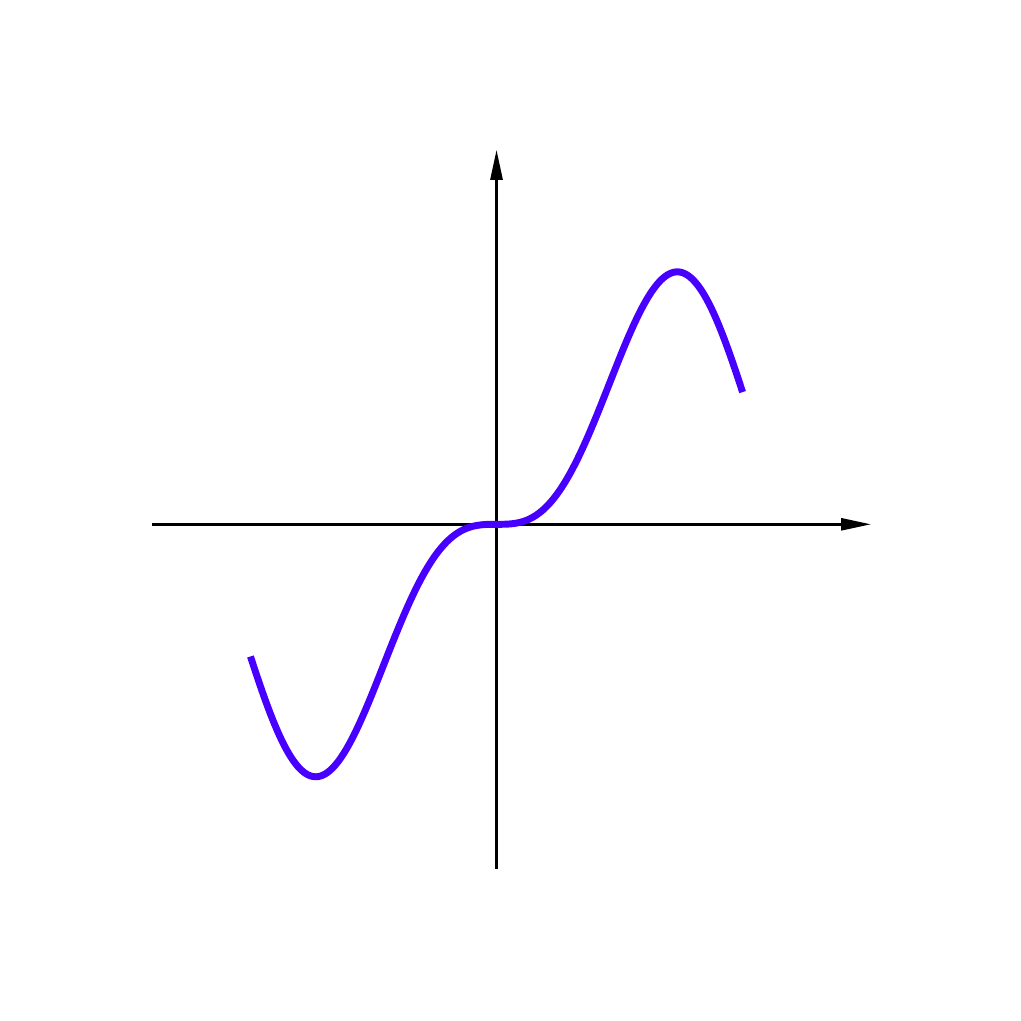}
		\caption{Пример функции $f(x)$.}
		\label{fig:functional_interval_example_1}
	
	\end{center}
\end{figure}

Если необходимо использовать функцию в дальнейших вычислениях, иногда имеет смысл приблизить
график функции и далее использовать это приближение. Простейший способ сделать это --- 
заключить график в интервал значений $f(x)$ на интервале $x \in [ \, a, \, b \, ]$. График функции $f(x)$ при $x \in [ \, a, \, b \, ]$ расположен 
в интервале $\big[ \, \min \, \text{ran}_{[a, \, b]} \, f(x), \, \max \, \text{ran}_{[a, \, b]}\, f(x) \, \big]:$

\begin{figure}[ht]
	\begin{center}
		
		\includegraphics[width = 0.40 \linewidth]{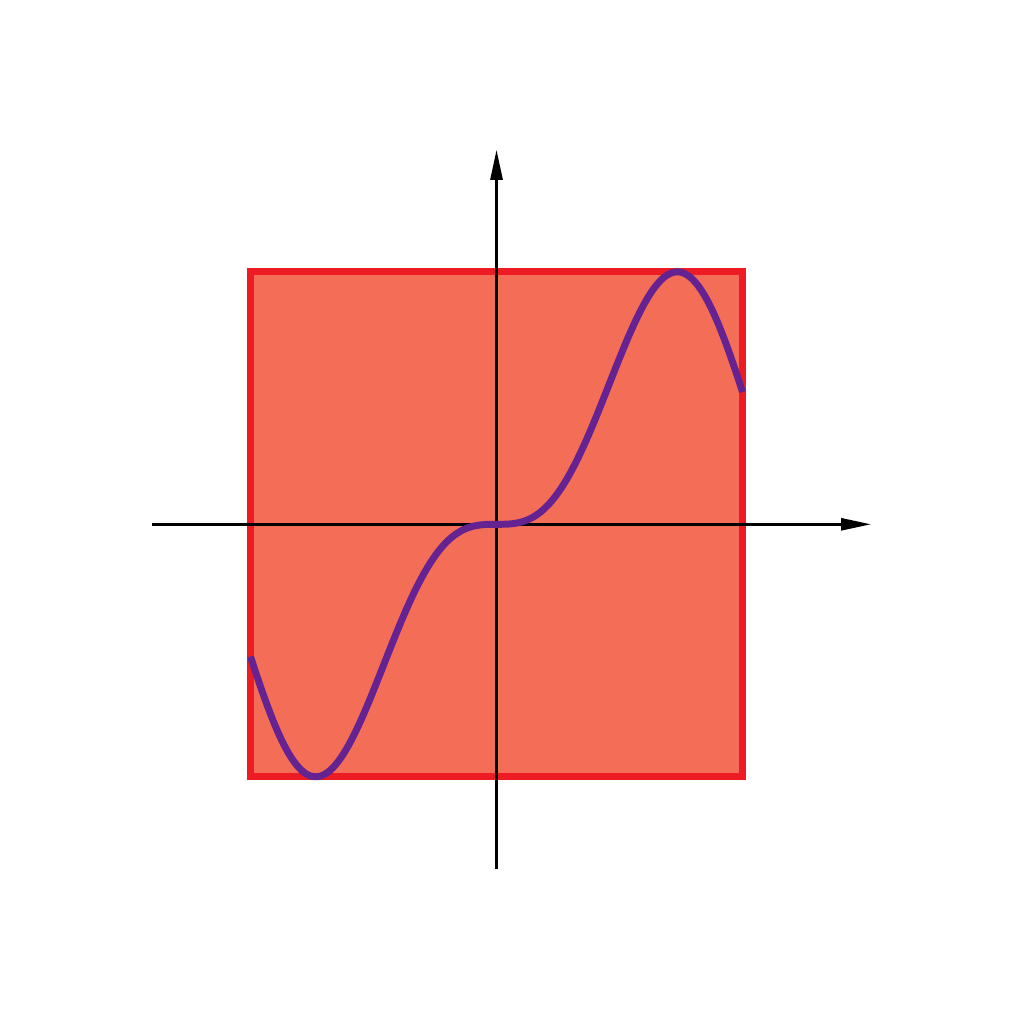}
		\caption{Пример оценки графика функции $f(x)$ интервалом\\ области значений функции.}
		\label{fig:functional_interval_example_1}
	
	\end{center}
\end{figure}

Однако видно, что в зависимости от характера функции $f(x)$, такой способ может достаточно давать достаточно грубую оценку графика $f(x)$. Поэтому иногда имеет смысл придумать более точную оценку, возможно и получаемую более сложным путём. Например, график данной функции $f(x)$ можно приблизить сверху и снизу квадратичными полиномами:

\begin{figure}[ht]
	\begin{center}
		
		\includegraphics[width = 0.40 \linewidth]{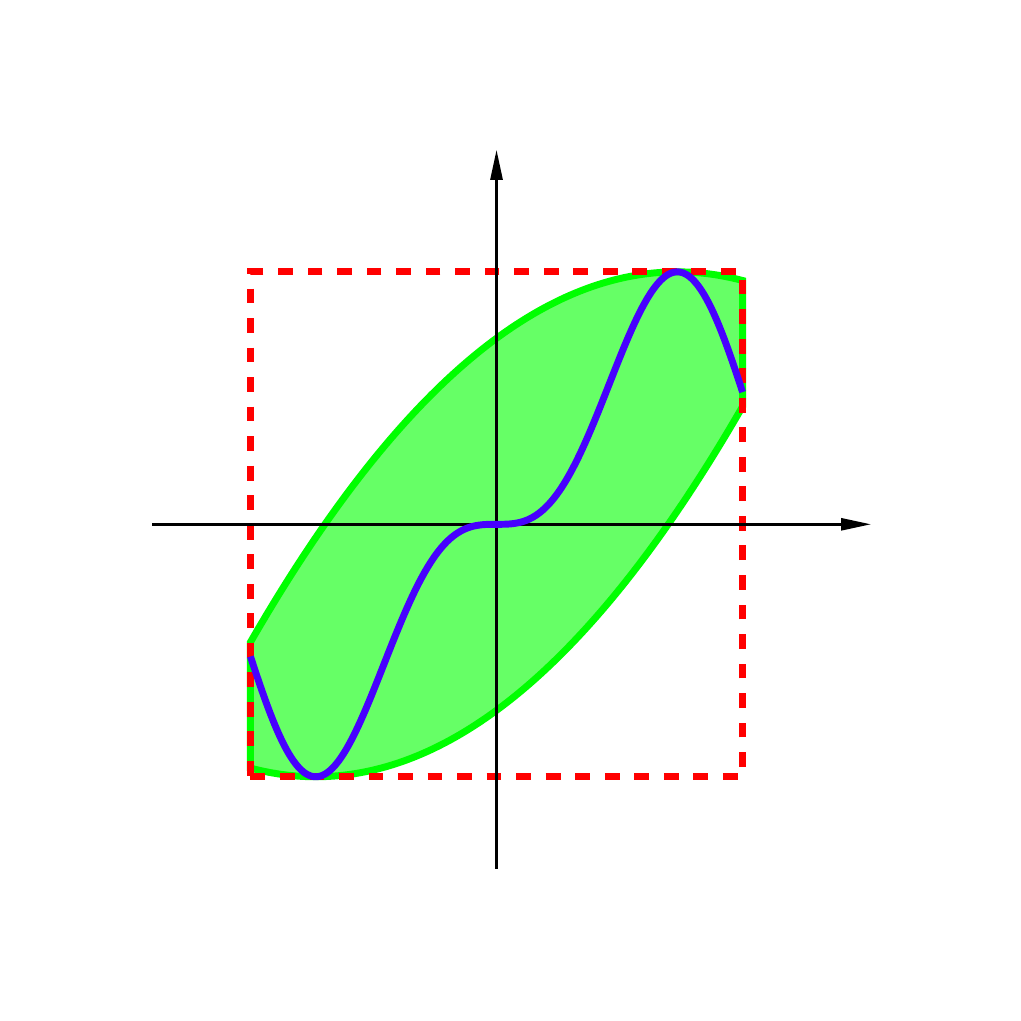}
		\caption{Пример приближения графика функции $f(x)$ \\квадратичными полиномами.}
		\label{fig:functional_interval_example_1}
	
	\end{center}
\end{figure}

Таким образом, использование функциональных интервалов перспективно в плане точности оценивания графиков функций.

\subsection{Понятие функционального интервала}

\begin{definition} \textit{Функциональным интервалом} назовём интервал, нижняя и верхняя границы которого представляются некоторыми функциями 
    \begin{equation*}
        \eus{L} : \unitint ^ {n} \rightarrow \mathbb{R} \qquad \text{и} \qquad \eus{U} : \unitint ^ {n} \rightarrow \mathbb{R}
    \end{equation*} 
соответственно, которые удовлетворяют свойству
    \begin{equation*}
    	\forall x_{i} \in \unitint, \qquad i = 1, \dots, n, \qquad \eus{L}(x_{1}, \dots, x_{n}) \leq \eus{U}(x_{1}, \dots, x_{n}).
    \end{equation*}
\end{definition}

Функции $\eus{L}$ и $\eus{U}$ будем называть нижней и верхней \textit{граничной функцией} соответственно. Область значений аргументов этих функций $x_{i}$ условна. Можно рассматривать любые интервалы изменения значений этих переменных. При этом любой интервал изменения переменной можно линейным преобразованием перевести в интервал $\unitint$.

\noindentСемейство функциональных интервалов будем обозначать $\mathbb{F}\mathbb{R}(x_{1}, \dots, x_{n})$, где $x_{1}, \dots, x_{n}$ --- аргументы функций $\eus{L}$ и $\eus{U}$. Введём обозначения таких интервалов согласно стандартам обозначений \cite{IntervalNotation}. Какой-либо конкретный функциональный интервал будем обозначать $\mbf{x}$. Его граничную функцию $\eus{L}$ будем обозначать как $\underline{\mbf{x}}$, а функцию $\eus{U}$, как $\overline{\mbf{x}}$.

Семейство функциональных интервалов без аргументов будем обозначать, как 
$\mathbb{F}\mathbb{R}_{0}$. Удобно мыслить это семейство функциональных интервалов, 
как параметрическое семейство классических интервалов $\big( \mathbb{I}\mathbb{R} \equiv 
\mathbb{F}\mathbb{R}_{0} \big)$. Также отметим, что семейство интервалов $\mathbb{I}\mathbb{R}$  будем называть классическими интервалами. А классические интервалы с арифметическими операциями между ними (описаны в книгах \cite{HansenWalster,Shary}), будем называть \textit{классической интервальной арифметикой}.

Для записи функциональных интервалов предлагается следующая нотация: вместо числовых границ 
интервалов будем писать на месте левого конца интервала функцию $\eus{L}$, а вместо правого 
конца --- функцию $\eus{U}$. Например, функциональный интервал $\mbf{x}$, у которого 
$\eus{L}(x) = x$, а $\eus{U}(x) = \sin(x) + 0.2$, будет записываться как $\mbf{x} 
= \big[ \, x, \, \sin(x) + 0.2 \, \big]$. 
  
\begin{figure}[ht]
	\begin{center}
	
	    \unitlength=1mm  
		
		\begin{picture}(60, 60)		
		
		    \put(0, 0){\includegraphics[height = 60mm]{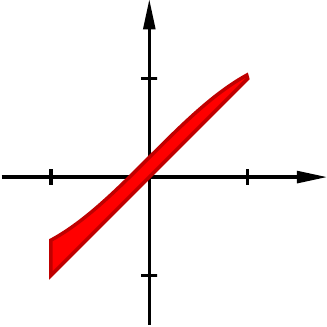}} 
		    
		    \put(57, 22){$x$}
		    \put(30, 57){$y$}
		    
		    \put(30, 44){$1$}
		    \put(30, 21){$0$}
		    \put(30, 8){$-1$}
		    
		    \put(5, 21){$-1$}
		    \put(45, 21){$1$}
		
		\end{picture}
		
		\caption{Пример функционального интервала $\mbf{x} = \big[ \, x, \, \text{sin}(x) + 0.2 \, \big]$.}
	    \label{fig:functional_interval_example} 
	
	\end{center}
\end{figure}

\subsection{Понятие линейного функционального интервала.}

\begin{definition} \textit{Линейным функциональным интервалом} назовём функциональный интервал, 
у которого функции $\eus{L}$ и $\eus{U}$ представляют собой линейные функции вида 
\begin{equation*}
\begin{array}{l} 
\eus{L}(x_{1}, \dots, x_{n}) = c^{\eus{L}} + \sum\limits_{i = 1}^{n} a_{i}^{\eus{L}} x_{i},
\\[4mm] 
\eus{U}(x_{1}, \dots, x_{n}) = c^{\eus{U}} + \sum\limits_{i = 1}^{n} a_{i}^{\eus{U}} x_{i}, 
\end{array} 
\qquad \forall i = 1, \dots, n, \qquad x_{i} \in \unitint. 
\end{equation*} 
  
Семейство таких интервалов будем обозначать $\mathbb{L}\mathbb{F}\mathbb{R}(x_{1}, \dots, 
x_{n})$, где $x_{1}, \dots, x_{n}$ --- аргументы функций $\eus{L}(x_{1}, \dots, x_{n})$ и 
$\eus{U}(x_{1}, \dots, x_{n})$. Семейство функциональных интервалов без аргументов будем  
обозначать, как $\mathbb{L}\mathbb{F}\mathbb{R}_{0}$ 
\end{definition}

Отметим, что альтернативное определение линейного функционального интервала было приведено в работе С.А. Югая \cite{Ugay}.

В работе будут рассмотрены однопараметрические интервалы семейства $\mathbb{L}\mathbb{F}\mathbb{R}(x)$. Интервал $\mbf{x}$ этого семейства будет записываться 
в виде 
\begin{equation*}
	\mbf{x} = \big[ \, \underline{a} x + \underline{b}, \, \overline{a} x + \overline{b} \, \big], \qquad x \in \unitint.
\end{equation*}
  
  
\subsection{Вспомогательные определения}

\begin{definition} 
Для функциональных интервалов $\mbf{x}$ и $\mbf{y}$ будем говорить, что \textit{интервал 
$\mbf{x}$ меньше интервала $\mbf{y}$} и писать $\mbf{x} < \mbf{y}$, если верхняя граница 
$\mbf{x}$ меньше  нижней границы $\mbf{y}$ для всех допустимых значений аргументов функций 
$\eus{L}$ и $\eus{U}$. Иными словами, для $\mbf{x}$, $\mbf{y}\in\mathbb{F}\mathbb{R}$ 
полагаем 
\begin{equation*}
\mbf{x} < \mbf{y} \ 
\Leftrightarrow   \ 
( \forall x_{i}\in\unitint, \, i = 1, \dots, n ) \;\ 
\overline{\mbf{x}(x_{1}, \dots, x_{n})} < \underline{\mbf{y}(x_{1}, \dots, x_{n})}. 
\end{equation*}

\end{definition}

\begin{definition} Аналогично для функциональных интервалов $\mbf{x}$ и $\mbf{y}$ будем говорить, что \textit{интервал $\mbf{x}$ больше интервала $\mbf{y}$}, если нижняя граница $\mbf{x}$ больше верхней границы $\mbf{y}$ для всех допустимых значений параметров функций $\eus{L}$ и $\eus{U}$, то есть
\begin{equation*}
\mbf{x} > \mbf{y} \ 
\Leftrightarrow   \ 
( \forall x_{i}\in\unitint, \, i = 1, \dots, n ) \;\ 
\underline{\mbf{x}(x_{1}, \dots, x_{n})} > \overline{\mbf{y}(x_{1}, \dots, x_{n})}. 
\end{equation*} 

\end{definition}

Обозначим через $F_{\downarrow}$ функцию,
\begin{equation*}
    F_{\downarrow} \, : \, \mathbb{F}\mathbb{R} (x_{1}, \dots, x_{n}) \times \mathbb{F}\mathbb{R} (x_{1}, \dots, x_{n}) \rightarrow \mathbb{F}\mathbb{R} (x_{1}, \dots, x_{n}),
\end{equation*} которая строит оценку снизу для двух интервалов семейства $\mathbb{F}\mathbb{R} (x_{1}, \dots, x_{n})$, то есть
\begin{equation*}
\left.\begin{array}{l}
    \forall \mbf{x}, \mbf{y} \in \mathbb{F}\mathbb{R}(x_{1}, \dots, x_{n}), \\
    \forall x_{i} \in \unitint, \qquad i = 1, \dots, n, 
\end{array}\right. \qquad \left\{\begin{array}{l}
    F_{\downarrow} ( \mbf{x}, \mbf{y} ) \leq \mbf{x}(x_{1}, \dots, x_{n}), \\
    F_{\downarrow} ( \mbf{x}, \mbf{y} ) \leq \mbf{y}(x_{1}, \dots, x_{n}).
\end{array}\right.
\end{equation*}

Аналогично введём функцию $F_{\uparrow}$,
\begin{equation*}
    F_{\uparrow} \, : \, \mathbb{F}\mathbb{R} (x_{1}, \dots, x_{n}) \times \mathbb{F}\mathbb{R} (x_{1}, \dots, x_{n}) \rightarrow \mathbb{F}\mathbb{R} (x_{1}, \dots, x_{n}),
\end{equation*} которая строит оценку сверху для двух интервалов семейства $\mathbb{F}\mathbb{R}$. Для неё верно
\begin{equation*}
\left.\begin{array}{l}
    \forall \mbf{x}, \mbf{y} \in \mathbb{F}\mathbb{R}(x_{1}, \dots, x_{n}), \\
    \forall x_{i} \in \unitint, \qquad i = 1, \dots, n, 
\end{array}\right. \qquad \left\{\begin{array}{l}
    F_{\uparrow} ( \mbf{x}, \mbf{y} ) \geq \mbf{x}(x_{1}, \dots, x_{n}), \\
    F_{\uparrow} ( \mbf{x}, \mbf{y} ) \geq \mbf{y}(x_{1}, \dots, x_{n}).
\end{array}\right.
\end{equation*}

\begin{definition} Областью значений $\text{ran}$ функционального интервала будем называть точный интервал его значений, то есть
\begin{equation*}
    \text{ran} \, \mbf{x}(x_{1}, \dots, x_{n}) = \Big[ \, \min\limits_{x_{1}, \dots, x_{n} \in \unitint} \underline{\mbf{x}(x_{1}, \dots, x_{n})}, \, \max\limits_{x_{1} , \dots, x_{n} \in \unitint} \overline{\mbf{x}(x_{1}, \dots, x_{n})}\Big].
\end{equation*}
\end{definition}

\noindent\textit{Примечание.} 

\noindentЗаметим, что для любого интервала $\mbf{x} \in \mathbb{L}\mathbb{F}\mathbb{R}_{0}$ область его значений, как функционального интервала, совпадает с интервалом $\mbf{x}$ в смысле классического интервала.
  
  
\subsection{Центрированные формы интервалов $\mathbb{I}\mathbb{R}$}

Пусть дана вещественная функция $f(x)$. Согласно \cite{Shary}, интервальная оценивающая функция $f(\mbf{x})$ имеет центрированную форму с центром $c$, если она представима в виде
\begin{equation*}
    f_{c}(\mbf{x}) = f(c) + \mbf{g} \, (\mbf{x} - c), \qquad g \in \mathbb{I}\mathbb{R}. 
\end{equation*}

Если $f(x)$ --- дифференцируемая функция, то в силу теоремы о среднем
\begin{equation*}
    f(x) - f(c) = f'(\xi) \, (x - c). 
\end{equation*}
Отсюда получаем, что
\begin{equation*}
    f_{c}(\mbf{x}) = f(c) + f'(\mbf{x}) \, (\mbf{x} - c).
\end{equation*}

Отличительная особенность центрированных форм заключается в квадратичной сходимости интервальной оценки функции с уменьшением интервала области определения оцениваемой функции $f(x)$.

\subsection{Вспомогательные утверждения} 
Докажем вспомогательные утверждения. 

\begin{lemma}[о верхней оценке линейных функций на $\unitint$]
\end{lemma} 
\noindentПусть
\begin{equation*}
	\begin{array}{c}
		f_{1}(x) = a_{1} x + b_{1}, \qquad f_{2}(x) = a_{2} x + b_{2}, \qquad x \in \unitint, \\ [1mm]
		M_{+1} = \text{max} \big\{ \, f_{1}(1), \, f_{2}(1) \, \big\}, \qquad M_{-1} = \text{max} \big\{ \, f_{1}(-1), \, f_{2}(-1) \, \big\}.
	\end{array}
\end{equation*}
Обозначим
\begin{equation}
\label{eq:M}
    M(x) = \frac{x}{2}\,\big( M_{+1} - M_{-1} \big) + \frac{1}{2} \big( M_{+1} + M_{-1} \big)
\end{equation}
Тогда
\begin{equation*}
\forall x \in \unitint \qquad \left\{\begin{array}{l}
    f_{1}(x) \leq M(x),\\
    f_{2}(x) \leq M(x).
\end{array}\right.
\end{equation*}

\begin{figure}[ht]
    \begin{center}
    
    \unitlength=1mm  
    
    \begin{picture}(80,80)		
    
        \put(0, 0){\includegraphics[height = 80mm]{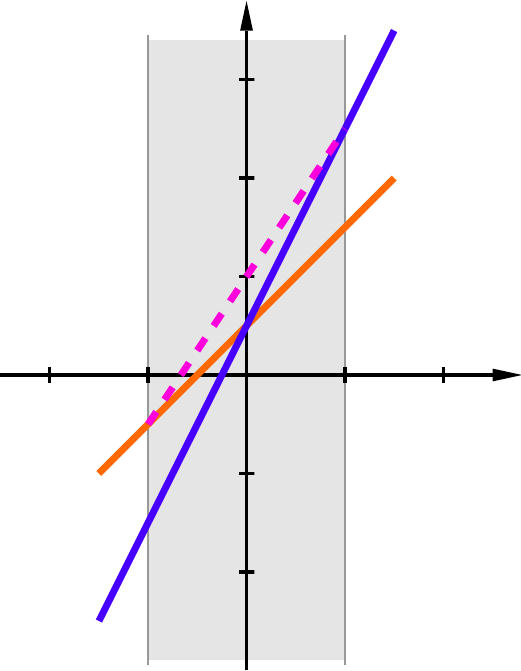}} 
        \put(50, 76){$y = x + 0.5$}
        \put(50, 57){$y = 2x + 0.5$}
        \put(-2, 43){$y = 1.5x + 1$}
        \put(60, 30){$x$}
        \put(32, 77){$y$}
        \put(32, 69){$3$}
        \put(32, 57){$2$}
        \put(32, 29){$0$}
        \put(32, 22){$-1$}
        \put(32, 10){$-2$}
        
        \put(1, 29){$-2$}
        \put(11, 29){$-1$}
        \put(42, 29){$1$}
        \put(52, 29){$2$}
        
    \end{picture} 	
	
	\caption{В результате применения Леммы 1 к отрезкам \\
		    $\big\{ (x, \, y) \mid x \in \unitint, \, y = x + 0.5 \, \big\}$ и $\big\{ (x, \, y) \mid x \in \unitint, \, y = 2x + 0.5 \, \big\}$ \\
		будет получен отрезок $\big\{ (x, \, y) \mid x \in \unitint, \, y = 1.5x + 1 \, \big\}$.}
		\label{fig:unity_up_example}
    	
    \end{center}
\end{figure}

\noindent\textbf{Доказательство}

\noindentРассмотрим четыре случая, исчерпывающие все возможные ситуации отношения функций $f_{1}(x)$ и $f_{2}(x)$.\\
\underline{\textit{Случай 1.}}\\
Пусть
\begin{equation*}
	\forall x \in \unitint \qquad f_{1}(x) \geq f_{2}(x).
\end{equation*}
Тогда
\begin{equation*}
	\begin{array}{l}
		\text{max} \{ \, f_{1}(1), \, f_{2}(1) \, \} = f_{1}(1) = a_{1} + b_{1}, \\ 
		\text{max} \{ \, f_{1}(-1), f_{2}(-1) \, \} = f_{1}(-1) = -a_{1} + b_{1}.
	\end{array}
\end{equation*}

Рассмотрим выражение (\ref{eq:M})
\begin{equation*}
	\begin{array}{l}
		M(x) = \frac{x}{2} \, \Big( \text{max} \big\{ \, f_{1}(1), \, f_{2}(1) \, \big\} - \text{max} \big\{ \, f_{1}(-1), \, f_{2}(-1) \, \big\} \Big) \\	 
		\hspace{2cm} + \, \frac{1}{2} \, \Big( \text{max} \big\{ \, f_{1}(1), \, f_{2}(1) \, \big\} + \text{max} \big\{ \, f_{1}(-1), \, f_{2}(-1) \, \big\} \Big) = \\
		= \frac{x}{2} \, (a_{1} + b_{1} + a_{1} - b_{1}) + \frac{1}{2} \, (a_{1} + b_{1} - a_{1} + b_{1}) = \\
		 = a_{1} x + b_{1} = f_{1}(x).
	\end{array}
\end{equation*}
Очевидно, что
\begin{equation*}
	\forall x \in \unitint \qquad M(x) = f_{1}(x) \geq f_{1}(x), \\
\end{equation*}
а также по предположению
\begin{equation*}
	\forall x \in \unitint \qquad M(x) = f_{1}(x) \geq f_{2}(x).
\end{equation*}
Заключение леммы верно, перейдём к рассмотрению второго случая.
\underline{\textit{Случай 2.}}\\
Пусть
\begin{equation*}
	\forall x \in \unitint \qquad f_{1}(x) \leq f_{2}(x).
\end{equation*}
Переобозначим функции как
\begin{equation*}
	f ^ {*}_{1}(x) = f_{2}(x), \qquad f ^ {*}_{2}(x) = f_{1}(x).
\end{equation*}
Тогда можно перейти к рассмотрению \underline{\textit{случая 1}} для функций $f_{1}^{*}(x)$ и $f_{2}^{*}(x)$.
\underline{\textit{Случай 3.}}\\
Пусть
\begin{equation*}
	\exists x^{*} \in \unitint\text{ такая, что }\\
	\left\{\begin{array}{lcl}
		\forall x \in [ \, -1, \, x^{*} \, ] & \qquad & f_{1}(x) \geq f_{2}(x),\\
		\forall x \in [ \, x^{*}, \, 1 \, ] & \qquad & f_{1}(x) \leq f_{2}(x).
	\end{array}\right.
\end{equation*}

\noindentТогда
\begin{equation*}
	\begin{array}{l}
		\text{max} \, \{ \, f_{1}(1), \, f_{2}(1) \, \} = f_{2}(1) = a_{2} + b_{2},\\
		\text{max} \, \{ \, f_{1}(-1), \, f_{2}(-1) \, \} = f_{1}(-1) = -a_{1} + b_{1}.
	\end{array}
\end{equation*}

Рассмотрим выражение
\begin{equation*}
	\begin{array}{l}
	    M(x) - f_{1}(x) =\\
		= \frac{x}{2}\,\big( M_{+1} - M_{-1} \big) + \frac{1}{2}\,\big( M_{+1} + M_{-1} \big) - f_{1}(x) =\\
		 = \frac{x}{2}\,(a_{2} + b_{2} + a_{1} - b_{1}) + \frac{1}{2}\,(a_{2} + b_{2} - a_{1} + b_{1}) - f_{1}(x) =\\
		= \frac{x}{2}\,(a_{2} + b_{2} + a_{1} - b_{1}) + \frac{1}{2}\,(a_{2} + b_{2} - a_{1} + b_{1}) - \frac{2}{2}\,a_{1} x - \frac{2}{2} \, b_{1} = \\
		= \frac{x}{2} \, \big( (a_{2} + b_{2} - a_{1} - b_{1}) + (a_{2} + b_{2} - a_{1} - b_{1}) \big) =\\
		= \frac{1}{2} \, (a_{2} + b_{2} - a_{1} - b_{1})\,(x + 1) =\\
		= \frac{1}{2}\,\big( f_{2}(1) - f_{1}(1) \big)\,(x + 1) \geq 0.
	\end{array}
\end{equation*}

\noindentАналогично рассмотрим выражение
\begin{equation*}
	\begin{array}{l}
	    M(x) - f_{2}(x) =\\
		= \frac{x}{2}\,\big( M_{+1} - M_{-1} \big) + \frac{1}{2}\,\big( M_{+1} + M_{-1} \big) - f_{2}(x) = \\
		= \frac{x}{2}\,(a_{2} + b_{2} + a_{1} - b_{1}) + \frac{1}{2}\,(a_{2} + b_{2} - a_{1} + b_{1}) - f_{2}(x) =\\
		= \frac{x}{2}\,(a_{2} + b_{2} + a_{1} - b_{1}) + \frac{1}{2}\,(a_{2} + b_{2} - a_{1} + b_{1}) - \frac{2}{2}\,a_{2} x - \frac{2}{2}\,b_{2} = \\
		= \frac{x}{2} \, \big( (-a_{2} + b_{2} + a_{1} - b_{1}) + (a_{2} - b_{2} - a_{1} + b_{1}) \big) =\\
		= \frac{1}{2}\,(-a_{2} + b_{2} + a_{1} - b_{1})\,(x - 1) = \\
		= \frac{1}{2}\,\big(f_{2}(-1) - f_{1}(-1) \big) \, (x - 1) = \\
		= \frac{1}{2}\,\big( f_{1}(-1) - f_{2}(-1) \big) \, (1 - x) \geq 0.
	\end{array}
\end{equation*}
Таким образом получаем, что
\begin{equation*}
    \forall x \in [ \, -1, \, 1 \, ] \qquad \left\{\begin{array}{l}
        f_{1}(x) \leq M(x), \\
        f_{2}(x) \leq M(x).
    \end{array}\right.
\end{equation*}
\underline{\textit{Случай 4.}}\\
Пусть
\begin{equation*}
	\exists x^{*} \in \unitint\text{ такая, что }
	\left\{\begin{array}{lcl}
		\forall x \in [ \, -1, \, x^{*} \, ] & \qquad & f_{1}(x) \leq f_{2}(x), \\
		\forall x \in [ \, x^{*}, \, 1 \, ] & \qquad & f_{1}(x) \geq f_{2}(x).
	\end{array}\right.
\end{equation*}
Переобозначим функции как
\begin{equation*}
	f_{1}^{*}(x) = f_{2}(x), \qquad f_{2}^{*}(x) = f_{1}(x).
\end{equation*}

Тогда заметим, что можем перейти к рассмотрению \underline{\textit{случая 3}} для функций $f_{1}^{*}(x)$ и $f_{2}^{*}(x)$. \hfill $\blacksquare$

\begin{lemma} [о нижней оценке линейных функций на $\unitint$]
\end{lemma}
\noindentПусть 
\begin{equation*}
	\begin{array}{l}
		f_{1}(x) = a_{1} x + b_{1}, \qquad x \in \unitint.\\
		f_{2}(x) = a_{2} x + b_{2}, \qquad x \in \unitint, \\
		M_{+1} = \min \{ \, f_{1}(1), \, f_{2}(1) \, \}, \\
		M_{-1} = \min \{ \, f_{1}(-1), \, f_{2}(-1) \, \}.
	\end{array}
\end{equation*}
Обозначим
\begin{equation*}
    M(x) = \frac{x}{2}\,\big( M_{+1} - M_{-1} \big) + \frac{1}{2} \big( M_{+1} + M_{-1} \big)
\end{equation*}
Тогда
\begin{equation*}
    \forall x \in \unitint \qquad \left\{\begin{array}{l}
        f_{1}(x) \geq M(x), \\
        f_{2}(x) \geq M(x).
    \end{array}\right.
\end{equation*}

\begin{figure}[ht]
    \begin{center}
    
    \unitlength=1mm  
    
    \begin{picture}(80,80)		
    
        \put(0, 0){\includegraphics[height = 80mm]{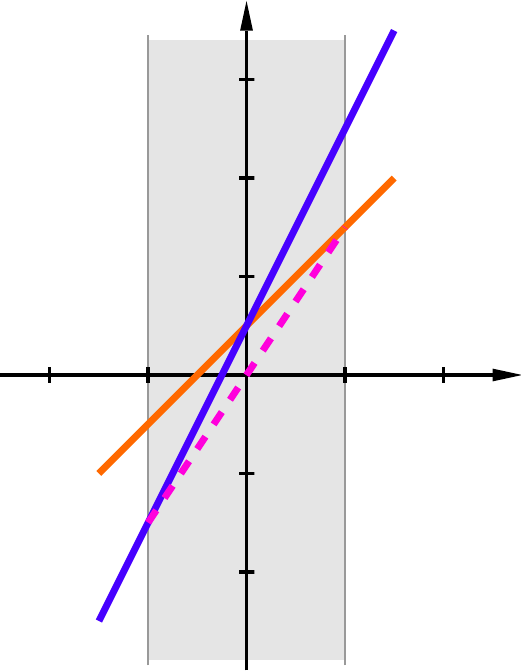}} 
        
        \put(50, 76){$y = x + 0.5$}
        \put(50, 57){$y = 2x + 0.5$}
        \put(37, 40){$y = 1.5x$}
        
        \put(60, 30){$x$}
        \put(32, 77){$y$}
        \put(32, 69){$3$}
        \put(32, 57){$2$}
        \put(32, 29){$0$}
        \put(32, 22){$-1$}
        \put(32, 10){$-2$}
        
        \put(1, 29){$-2$}
        \put(11, 29){$-1$}
        \put(42, 29){$1$}
        \put(52, 29){$2$}
        
    \end{picture} 	
	
	\caption{В результате применения Леммы 2 к отрезкам \\
		    $\big\{ (x, \, y) \mid x \in \unitint, \, y = x + 0.5 \, \big\}$ и $\big\{ (x, \, y) \mid x \in \unitint, \, y = 2x + 0.5 \, \big\}$ \\
		будет получен отрезок $\big\{ (x, \, y) \mid x \in \unitint, \, y = 1.5x\, \big\}$.}
		\label{fig:unity_down_example}
    	
    \end{center}
\end{figure}

\noindent\textbf{Доказательство}

\noindentИз Леммы 1 следует, что
\begin{equation*}
	\begin{array}{l}
		-f_{1}(x), \, -f_{2}(x) \leq \\
		\hspace{0.5cm} \leq \frac{x}{2}\,\Big( \text{max} \big\{ -f_{1}(1), \, -f_{2}(1) \, \big\} - \text{max} \big\{ -f_{1}(-1), \, -f_{2}(-1) \, \big\} \Big) \\
		\hspace{1cm} + \, \frac{1}{2}\,\Big( \text{max} \big\{ -f_{1}(1), \, -f_{2}(1) \, \big\} + \text{max} \big\{ -f_{1}(-1), \, -f_{2}(-1) \, \big\} \Big).
	\end{array}
\end{equation*}
Тогда
\begin{equation*}
	\begin{array}{l}
		f_{1}(x), \, f_{2}(x) \geq \\
		\hspace{0.5cm} \geq -\frac{x}{2}\,\Big( \text{max} \big\{ \, -f_{1}(1), \, -f_{2}(1) \, \big\} - \text{max} \big\{ \, -f_{1}(-1), \, -f_{2}(-1) \, \big\} \Big) \\
		\hspace{1cm} - \, \frac{1}{2}\,\Big( \text{max} \big\{ \, -f_{1}(1), \, -f_{2}(1) \, \big\} + \text{max} \big\{ \, -f_{1}(-1), \, -f_{2}(-1) \, \big\} \Big) = \\
		\hspace{0.5cm} = \frac{x}{2}\,\Big(- \text{max} \big\{ -f_{1}(1), \, -f_{2}(1) \, \big\} + \text{max} \big\{ -f_{1}(-1), \, -f_{2}(-1) \, \big\} \Big) \\
		\hspace{1cm} + \, \frac{1}{2}\,\Big(- \text{max} \big\{ \, -f_{1}(1), \, -f_{2}(1) \, \big\} - \text{max} \big\{ \, -f_{1}(-1), \, -f_{2}(-1) \, \big\} \Big).
	\end{array}
\end{equation*}
Заметим, что
\begin{equation*}
	-\text{max} \, \big\{ -f(x) \big\} = \text{min} \, \big\{ f(x) \big\}.
\end{equation*}

\noindentТогда получим, что
\begin{equation*}
	\begin{array}{l}
		f_{1}(x), \, f_{2}(x) \geq \\
		\hspace{1cm} \geq \frac{x}{2}\,\Big( \text{min} \big\{ \, f_{1}(1), \, f_{2}(1) \, \big\} - \text{min} \big\{ \, f_{1}(-1), \, f_{2}(-1) \, \big\} \Big) \\
		\hspace{2cm} + \, \frac{1}{2}\,\Big(\text{min} \big\{ \, f_{1}(1), \, f_{2}(1) \, \big\} + \text{min} \big\{ \, f_{1}(-1), \, f_{2}(-1) \, \big\} \Big).
	\end{array}
\end{equation*}
\hfill $\blacksquare$

\begin{lemma} [о нижней оценке выпуклой функции на $\unitint$]
\end{lemma}
\noindentПусть $f(x)$ --- дважды непрерывно дифференцируемая функция и
\begin{equation*}
	\begin{array}{c}
	     f''(x) \geq 0, \qquad x \in \unitint, \\ [2mm]
	     \min\limits_{x \in \unitint} \, \big| f'(x) \big| = \big|f'(x ^ {*})\big|.
	\end{array}
\end{equation*}
\noindentТогда 
\begin{equation*}
	\forall x \in \unitint \qquad f(x) \geq f(x ^ {*}) + (x - x ^ {*})\,f'(x ^ {*}).
\end{equation*}

\begin{figure}[ht]
    \begin{center}
    
    \unitlength=1mm  
    
    \begin{picture}(80,80)		
    
        \put(0, 0){\includegraphics[height= 80mm]{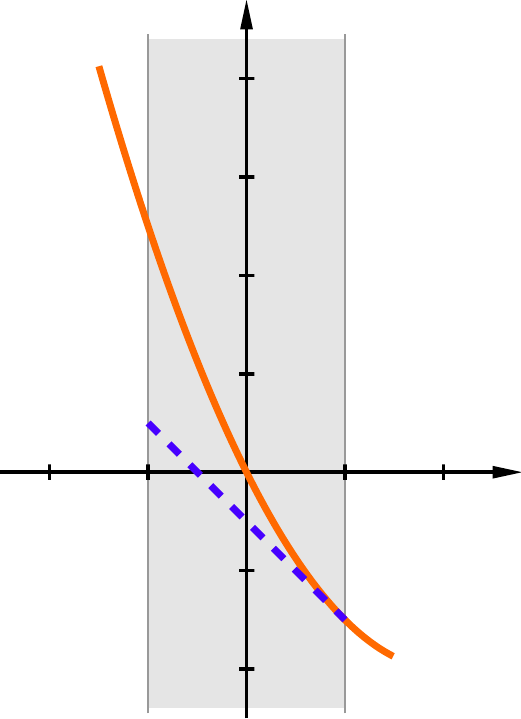}} 
        \put(-10, 34){$y = x + 0.5$}
        \put(46, 6){$y = 0.5 \, x ^ {2} - 2 x$}
        \put(55, 23){$x$}
        \put(30, 78){$y$}
        
        \put(30, 69){$4$}
        \put(30, 59){$3$}
        \put(30, 48){$2$}
        \put(30, 37){$1$}
        \put(30, 29){$0$}
        \put(30, 4){$-2$}
        
        \put(0, 21){$-2$}
        \put(9, 21){$-1$}
        \put(39, 21){$1$}
        \put(48, 21){$2$}
        
    \end{picture} 	
	
	\caption{В результате применения Леммы 3 \\ к функции
		квадратичного полинома $f(x) = 0.5 \, x ^ {2} - 2x, \, x \in \unitint$ \\
		будет получен отрезок $\big\{ (x, \, y) \mid x \in \unitint, \, y = -x - 0.5 \, \big\}$.}
		\label{fig:lemma3_example}
    	
    \end{center}
\end{figure}

\noindent\textbf{Доказательство}\\
Проведём через точку $x = x ^ {*}$ касательный отрезок к графику функции $f(x)$, $x \in \unitint$.

Тогда этот отрезок будет определяться выражением
\begin{equation*}
	g(x) = f(x^{*}) + (x - x^ {*})\,f'(x^{*}).
\end{equation*}
В силу того, что отрезок касательный, имеем
\begin{equation*}
	f(x^{*}) = g(x^{*}), \qquad f'(x^{*}) = g'(x^{*}).
\end{equation*}
Также заметим, что
\begin{equation*}
	f''(x^{*}) \geq 0 \geq g''(x^{*}) = 0.
\end{equation*}
Таким образом, $f(x) \geq g(x) = f(x^{*}) + (x - x ^ {*})\,f'(x ^ {*})$. \hfill $\blacksquare$\\
\textbf{Следствие из Леммы 3}\\
В условиях Леммы 3 верно, что
\begin{equation*}
    \min\limits_{x \in \unitint} f(x) = \min\limits_{x \in \unitint} \big( f(x^{*}) + (x - x^{*}) \, f'(x^{*}) \big).
\end{equation*}
\textbf{Доказательство}\\
Так как $f''(x) > 0$, то уравнение $f'(x) = 0$ имеет единственный корень на вещественной оси $\mathbb{R}$ и не более одного корня на интервале $\unitint$.
Рассмотрим три взаимоисключающих случая:\\
\underline{\textit{Случай 1.}}\\
Пусть корень уравнения $f'(x) = 0$ лежит в интервале $\unitint$. Обозначим его $x_{0}$. 
Заметим, что тогда $x^{*} = x_{0}$, так как 
\begin{equation*}
    \min\limits_{x \in \unitint}\big|f'(x)\big| = 0 = f'(x_{0}) = \big| f'(x_{0}) \big|.
\end{equation*}
Тогда верно
\begin{equation*}
    f(x^{*}) + (x - x^{*}) \, f'(x^{*}) = f(x^{*}).
\end{equation*}

Заметим, что так как $f'(x^{*}) = 0$ и $f''(x) > 0$, то $x^{*}$ --- точка глобального минимума функции, а также, соответственно, на интервале $\unitint$. То есть, 
\begin{equation*}
    f(x^{*}) = \min\limits_{x \in \unitint} f(x).
\end{equation*}

Тогда очевидно, что 
\begin{equation*}
    \min\limits_{x \in \unitint} f(x) = f(x^{*}) = \min\limits_{x \in \unitint} \big( f(x^{*}) + (x - x^{*}) \, f'(x^{*}) \big).
\end{equation*}
\underline{\textit{Случай 2.}}\\
Пусть корень уравнения $f'(x) = 0$ лежит левее интервала $\unitint$. Обозначим его $x_{-}$. То есть $x_{-} < -1$.

Так как $f'(x_{-}) = 0$ и $f''(x) > 0$, тогда $x^{*} = -1$ и
\begin{equation*}
    \forall x \in \unitint \qquad f'(x) > 0,
\end{equation*}
а значит
\begin{equation*}
    \min\limits_{x \in \unitint} f(x) = f(-1).
\end{equation*}
Тогда верно
\begin{equation*}
    \begin{array}{l}
        \min\limits_{x \in \unitint} \big( f(x^{*}) + (x - x^{*}) \, f'(x^{*}) \big) = \\
         \hspace{0.5cm} = \min\limits_{x \in \unitint} \big( f(-1) + (x + 1) \, f'(-1) \big) = \\
         \hspace{1cm} = f(-1) + f'(-1) \min\limits_{x \in \unitint}(x + 1) = \\
         \hspace{1.5cm} = f(-1) + f'(-1) \, (\min\limits_{x \in \unitint} x + 1) = f(-1) = \min\limits_{x \in \unitint} f(x).
    \end{array}
\end{equation*}

\noindent\underline{\textit{Случай 3.}}\\
Пусть корень уравнения $f'(x) = 0$ лежит правее интервала $\unitint$. Обозначим его $x_{+}$. То есть $x_{+} > 1$.

Так как $f'(x_{+}) = 0$ и $f''(x) > 0$, то $x^{*} = 1$ и
\begin{equation*}
    \forall x \in \unitint \qquad f'(x) < 0,
\end{equation*}
а значит
\begin{equation*}
    \min\limits_{x \in \unitint} f(x) = f(1).
\end{equation*}

Тогда верно
\begin{equation*}
    \begin{array}{l}
        \min\limits_{x \in \unitint} \big( f(x^{*}) + (x - x^{*}) \, f'(x^{*}) \big) = \\
         \hspace{0.5cm} = \min\limits_{x \in \unitint} \big( f(1) + (x - 1) \, f'(1) \big) = \\
         \hspace{1cm} = f(1) + f'(1) \max\limits_{x \in \unitint}(x - 1) = \\
         \hspace{1.5cm} = f(1) + f'(1) \, (\max\limits_{x \in \unitint} x - 1) = f(1) = \min\limits_{x \in \unitint} f(x).
    \end{array}
\end{equation*}

\hfill $\blacksquare$

\begin{lemma} [о верхней оценке вогнутой функции на $\unitint$]
\end{lemma} 
\noindentПусть $f(x)$ --- дважды непрерывно дифференцируемая функция и
\begin{equation*}
    \begin{array}{cc}
        f''(x) \leq 0, \qquad x \in \unitint, \\ [2mm]
        \min\limits_{x \in \unitint} \, \big| f'(x) \big| = \big| f'(x ^ {*}) \big|.
    \end{array}    
\end{equation*}
Тогда
\begin{equation*}
	\forall x \in \unitint \qquad f(x) \leq f(x ^ {*}) + (x - x ^ {*})\,f'(x ^ {*}).
\end{equation*}

\begin{figure}[ht]
    \begin{center}
    
    \unitlength=1mm  
    
    \begin{picture}(80,80)		
    
        \put(0, 0){\includegraphics[height = 80mm]{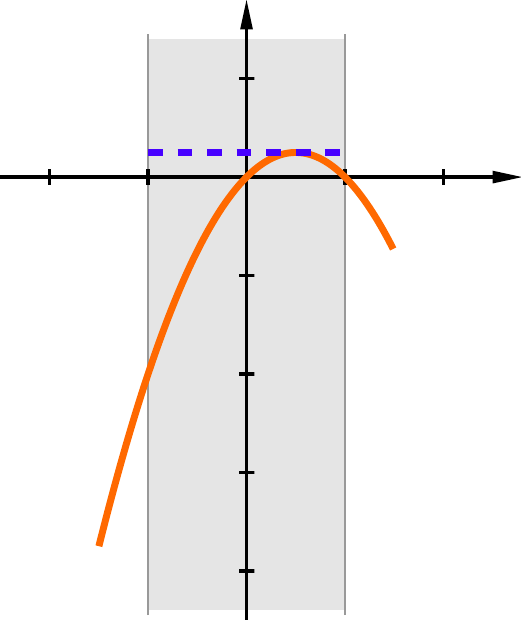}} 
        \put(-8, 59){$y = x + 0.5$}
        \put(-16, 9){$y = -x ^ {2} + x$}
        \put(63, 52){$x$}
        \put(35, 77){$y$}
        
        \put(34, 69){$1$}
        \put(34, 52){$0$}
        \put(34, 43){$-1$}
        \put(34, 31){$-2$}
        \put(34, 18){$-3$}
        \put(34, 5){$-4$}
        
        \put(1, 51){$-2$}
        \put(12, 51){$-1$}
        \put(45, 51){$1$}
        \put(56, 51){$2$}
        
    \end{picture} 	
	
	\caption{В результате применения Леммы 4\\ к функции
		квадратичного полинома $f(x) = -x ^ {2} + x, \, x \in \unitint$ \\
		будет получен отрезок $\big\{ (x, \, y) \mid x \in \unitint, \, y = 0.25 \, \big\}$.}
		\label{fig:lemma4_example}
    	
    \end{center}
\end{figure}

\noindent\textbf{Доказательство}\\
Проведём через точку $x =\,x ^ {*}$ касательный отрезок к графику функции $f(x)$, $x \in \unitint$.

Тогда этот отрезок будет определяться выражением 
\begin{equation*}
	g(x) = f(x^{*}) + (x - x ^ {*})\,f'(x^{*}).
\end{equation*}
В силу того, что отрезок касательный, имеем
\begin{equation*}
	f(x^{*}) = g(x^{*}), \qquad f'(x^{*}) = g'(x^{*}).
\end{equation*}

Также заметим, что
\begin{equation*}
	f''(x^{*}) \leq 0 \leq g''(x^{*}) = 0.
\end{equation*}
Таким образом, $f(x) \leq g(x) = f(x ^ {*}) + (x - x^{*})\,f'(x ^ {*})$. \hfill $\blacksquare$\\
\textbf{Следствие из Леммы 4}\\
В условиях Леммы 4 верно, что
\begin{equation*}
    \max\limits_{x \in \unitint} f(x) = \max\limits_{x \in \unitint} \big( f(x^{*}) + (x - x^{*}) \, f'(x^{*}) \big).
\end{equation*}
\textbf{Доказательство}\\
Так как $f''(x) < 0$, то уравнение $f'(x) = 0$ имеет единственный корень на вещественной оси $\mathbb{R}$ и не более одного корня на интервале $\unitint$.
Рассмотрим три взаимоисключающих случая:\\
\underline{\textit{Случай 1.}}\\
Пусть корень уравнения $f'(x) = 0$ лежит в интервале $\unitint$. Обозначим его $x_{0}$. 
Заметим, что тогда $x^{*} = x_{0}$, так как 
\begin{equation*}
    \min\limits_{x \in \unitint}\big|f'(x)\big| = 0 = f'(x_{0}) = \big| f'(x_{0}) \big|.
\end{equation*}
Тогда верно
\begin{equation*}
    f(x^{*}) + (x - x^{*}) \, f'(x^{*}) = f(x^{*}).
\end{equation*}

Заметим, что так как $f'(x^{*}) = 0$ и $f''(x) < 0$, то $x^{*}$ --- точка глобального максимума функции, а также, соответственно, на интервале $\unitint$. То есть, 
\begin{equation*}
    f(x^{*}) = \max\limits_{x \in \unitint} f(x).
\end{equation*}

Тогда очевидно, что 
\begin{equation*}
    \max\limits_{x \in \unitint} f(x) = f(x^{*}) = \max\limits_{x \in \unitint} \big( f(x^{*}) + (x - x^{*}) \, f'(x^{*}) \big).
\end{equation*}
\underline{\textit{Случай 2.}}\\
Пусть корень уравнения $f'(x) = 0$ лежит левее интервала $\unitint$. Обозначим его $x_{-}$. То есть $x_{-} < -1$.

Так как $f'(x_{-}) = 0$ и $f''(x) < 0$, то $x^{*} = -1$ и
\begin{equation*}
    \forall x \in \unitint \qquad f'(x) < 0,
\end{equation*}
а значит
\begin{equation*}
    \max\limits_{x \in \unitint} f(x) = f(-1).
\end{equation*}

Тогда верно
\begin{equation*}
    \begin{array}{l}
        \max\limits_{x \in \unitint} \big( f(x^{*}) + (x - x^{*}) \, f'(x^{*}) \big) = \\
         \hspace{0.5cm} = \max\limits_{x \in \unitint} \big( f(-1) + (x + 1) \, f'(-1) \big) = \\
         \hspace{1cm} = f(-1) + f'(-1) \min\limits_{x \in \unitint}(x + 1) = \\
         \hspace{1.5cm} = f(-1) + f'(-1) \, (\min\limits_{x \in \unitint} x + 1) = f(-1) = \max\limits_{x \in \unitint} f(x).
    \end{array}
\end{equation*}
\underline{\textit{Случай 3.}}\\
Пусть корень уравнения $f'(x) = 0$ лежит правее интервала $\unitint$. Обозначим его $x_{+}$. То есть $x_{+} > 1$.

Так как $f'(x_{+}) = 0$ и $f''(x) < 0$, то $x^{*} = 1$ и
\begin{equation*}
    \forall x \in \unitint \qquad f'(x) > 0,
\end{equation*}
а значит
\begin{equation*}
    \max\limits_{x \in \unitint} f(x) = f(1).
\end{equation*}

Тогда верно
\begin{equation*}
    \begin{array}{l}
        \max\limits_{x \in \unitint} \big( f(x^{*}) + (x - x^{*}) \, f'(x^{*}) \big) = \\
         \hspace{0.5cm} = \max\limits_{x \in \unitint} \big( f(1) + (x - 1) \, f'(1) \big) = \\
         \hspace{1cm} = f(1) + f'(1) \max\limits_{x \in \unitint}(x - 1) = \\
         \hspace{1.5cm} = f(1) + f'(1) \, (\max\limits_{x \in \unitint} x - 1) = f(1) = \max\limits_{x \in \unitint} f(x).
    \end{array}
\end{equation*}

\hfill $\blacksquare$

\begin{lemma} [о верхней оценке выпуклой функции на $\unitint$]
\end{lemma}

\noindentПусть $f(x)$ --- дважды дифференцируемая функция и
\begin{equation*}
	f''(x) \geq 0, \qquad x \in \unitint.
\end{equation*}
Обозначим
\begin{equation*}
    M(x) = \frac{x}{2}\,\big( f(1) - f(-1) \big) + \frac{1}{2}\,\big( f(1) + f(-1) \big).
\end{equation*}
Тогда 
\begin{equation*}
	\forall x \in \unitint \qquad f(x) \leq M(x).
\end{equation*}

\begin{figure}[ht]
    \begin{center}
    
    \unitlength=1mm  
    
    \begin{picture}(80,80)		
    
        \put(0, 0){\includegraphics[height = 80mm]{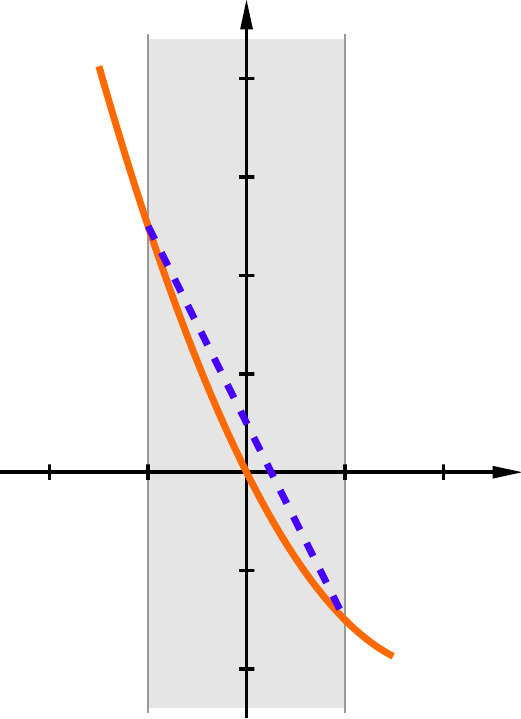}} 
        \put(-23, 71){$y = 0.5x ^ {2} - 2 x$}
        \put(47, 6){$y = -2x + 0.5$}
        \put(55, 23){$x$}
        \put(31, 78){$y$}
        
        \put(30, 70){$4$}
        \put(30, 59){$3$}
        \put(30, 48){$2$}
        \put(30, 37){$1$}
        \put(30, 4){$-2$}
        
        \put(0, 21){$-2$}
        \put(9, 21){$-1$}
        \put(39, 21){$1$}
        \put(48, 21){$2$}
        
    \end{picture} 	
	
	\caption{В результате применения Леммы 5\\ к функции 
		квадратичного полинома $f(x) = 0.5x ^ {2} - 2x, \, x \in \unitint$ \\
		будет получен отрезок $\big\{ \, (x, \, y) \mid x \in \unitint, \, y = -2x + 0.5 \, \big\}$.}
		\label{fig:lemma5_example}
    	
    \end{center}
\end{figure}

\noindent\textbf{Доказательство}

\noindentРассмотрим функцию
\begin{equation*}
	g(x) = f(x) - M(x) = \frac{x}{2}\,\big( f(1) - f(-1) \big) - \frac{1}{2}\,\big( f(1) + f(-1) \big).
\end{equation*}

Заметим, что
\begin{equation*}
	\begin{array}{rcl}
		g(1) & = & f(1) - \frac{1}{2}\,\big( f(1) - f(-1) + f(1) + f(-1) \big) = \\
		& = & f(1) - \frac{2}{2} \, f(1) = 0, \\
		g(-1) & = & f(-1) - \frac{1}{2}\,\big( -f(1) + f(-1) + f(1) + f(-1) \big) = \\
		& = & f(-1) - \frac{2}{2} \, f(-1) = 0.
	\end{array}
\end{equation*}

Получаем, что выражение 
\begin{equation*}
	M(x) = \frac{x}{2}\,\big( f(1) - f(-1) \big) + \frac{1}{2}\,\big( f(1) + f(-1) \big) \qquad x \in [ \, 1, \, 1 \, ]
\end{equation*}
определяет отрезок, соединяющий точки графика функции $\big( -1, \, f(-1) \big)$ и $\big(1, \, f(1) \big)$. Тогда из определения выпуклой функции получаем, что 
\begin{equation*}
    f(x) \leq \frac{x}{2}\,\big( f(1) - f(-1) \big) + \frac{1}{2}\,\big( f(1) + f(-1) \big) = M(x).
\end{equation*} \hfill $\blacksquare$\\
\textbf{Следствие из Леммы 5}\\
В условиях Леммы 5 верно, что
\begin{equation*}
    \min\limits_{x \in \unitint} f(x) = \min\limits_{x \in \unitint} M(x).
\end{equation*}
\textbf{Доказательство}\\
Из того, что $f''(x) \geq 0$ на интервале $\unitint$, следует, что 
\begin{equation*}
    \max\limits_{x \in \unitint} f(x) = \max\big\{ f(-1), \, f(1) \big\}.
\end{equation*}

Так как функция
\begin{equation*}
    M(x) = \frac{x}{2} \, \big( f(1) - f(-1) \big) + \frac{1}{2} \, \big( f(1) + f(-1) \big) \qquad x \in \unitint
\end{equation*}
задает отрезок, который соединяет точки $\big( -1, \, f(-1) \big)$ и $\big( 1, \, f(1) \big)$, то 
\begin{equation*}
    \max\limits_{x \in \unitint} f(x) = \max\limits_{x \in \unitint} \Big(\frac{x}{2} \, \big( f(1) - f(-1) \big) + \frac{1}{2} \, \big( f(1) + f(-1) \big) \Big) = \max\limits_{x \in \unitint} M(x). 
\end{equation*}
\hfill $\blacksquare$

\begin{lemma} [о нижней оценке вогнутой функции на $\unitint$]
\end{lemma}
\noindentПусть $f(x)$ --- дважды дифференцируемая функция и
\begin{equation*}
	f''(x) \leq 0, \qquad x \in \unitint.
\end{equation*}
Обозначим
\begin{equation*}
    M(x) = \frac{x}{2}\,\big( f(1) - f(-1) \big) + \frac{1}{2}\,\big( f(1) + f(-1) \big).
\end{equation*}
Тогда
\begin{equation*}
	\forall x \in \unitint \qquad f(x) \geq M(x).
\end{equation*}

\begin{figure}[ht]
    \begin{center}
    
    \unitlength=1mm  
    
    \begin{picture}(80,80)		
    
        \put(0, 0){\includegraphics[height = 80mm]{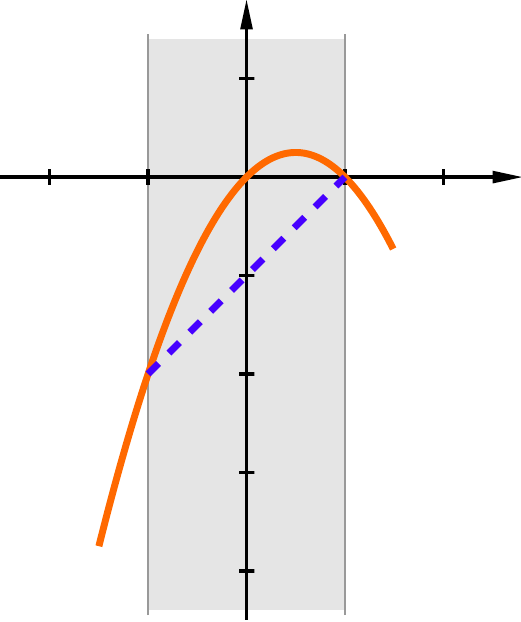}} 
        \put(-4, 30){$y = x - 1$}
        \put(53, 46){$y = -x ^ {2} + x$}
        \put(65, 52){$x$}
        \put(34, 77){$y$}
        
        \put(34, 69){$1$}
        \put(34, 51){$0$}
        \put(34, 43){$-1$}
        \put(34, 31){$-2$}
        \put(34, 18){$-3$}
        \put(34, 5){$-4$}
        
        \put(1, 51){$-2$}
        \put(12, 51){$-1$}
        \put(45, 51){$1$}
        \put(56, 51){$2$}
        
    \end{picture} 	
	
	\caption{В результате применения Леммы 6 \\ к функции
		квадратичного полинома $f(x) = -x ^ {2} + x, \, x \in \unitint$ \\
		будет получен отрезок $\big\{ \, (x, \, y) \mid x \in \unitint, \, y = x - 1 \, \big\}$.}
		\label{fig:lemma6_example}
    	
    \end{center}
\end{figure}

\noindent\textbf{Доказательство}\\
Рассмотрим функцию $-f(x)$. Тогда
\begin{equation*}
    -f''(x) \geq 0.
\end{equation*}

Следовательно к функции $-f(x)$ применима Лемма 5, то есть верно
\begin{equation*}
    -f(x) \leq \frac{x}{2}\,\big( -f(1) + f(-1) \big) + \frac{1}{2}\,\big( -f(1) - f(-1) \big) = - M(x).
\end{equation*}

Заметим, что это неравенство эквивалентно
\begin{equation*}
    f(x) \geq M(x).
\end{equation*}
\hfill $\blacksquare$\\
\textbf{Следствие из Леммы 6}\\
В условиях Леммы 6 верно, что
\begin{equation*}
    \min\limits_{x \in \unitint} f(x) = \min\limits_{x \in \unitint} M(x).
\end{equation*}
\textbf{Доказательство}\\
Из того, что $f''(x) \leq 0$ на интервале $\unitint$, следует, что 
\begin{equation*}
    \min\limits_{x \in \unitint} f(x) = \min\big\{ f(-1), \, f(1) \big\}.
\end{equation*}

Так как функция
\begin{equation*}
    M(x) = \frac{x}{2} \, \big( f(1) - f(-1) \big) + \frac{1}{2} \, \big( f(1) + f(-1) \big) \qquad x \in \unitint
\end{equation*} 
задает отрезок, который соединяет точки $\big( -1, \, f(-1) \big)$ и $\big( 1, \, f(1) \big)$, то 

\begin{equation*}
    \min\limits_{x \in \unitint} f(x) = \min\limits_{x \in \unitint} \Big(\frac{x}{2} \, \big( f(1) - f(-1) \big) + \frac{1}{2} \, \big( f(1) + f(-1) \big) \Big) = \min\limits_{x \in \unitint} M(x).
\end{equation*}

\hfill $\blacksquare$

\begin{lemma}[о выборе опорной точки приближения квадратичного полинома на $\unitint$]
\end{lemma}

\noindentПусть
\begin{equation*}
	\begin{array}{l}
		f(x) = a x^{2} + b x + c, \qquad x \in \unitint, \\ [1mm]
		a \neq 0, \qquad a, \, b, \, c \in \mathbb{R}, \\ [1mm]
		\min\limits_{ \, x \in \unitint } \, \big| f'(x) \big| = \big| f'(x^{*}) \big|, \\ [1mm]
		x_{\textit{extr}} = - b \, / \, (2a).
	\end{array}
\end{equation*}

\noindentТогда
\begin{equation*}
	x^{*} = \left\{ \begin{array}{lr}
		x_{\textit{extr}},\text{ если }x_{\textit{extr}} \in \unitint, & (1)\\
		-1,\text{ если }x_{\textit{extr}} < -1, & (2) \\
		1,\text{ иначе (если }x_{\textit{extr}} > 1). & (3)
	\end{array} \right.
\end{equation*}

\begin{figure}[ht]
	\begin{center}
		
		\unitlength=1mm  
    
        \begin{picture}(160, 50)		
        
            \put(0, 0){\includegraphics[height = 50mm]{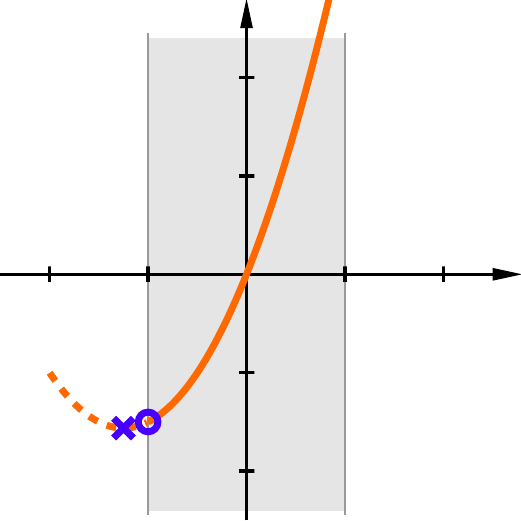}} 
            
            \put(25, 48){\small $y$}
            \put(47, 19){\small $x$}
            
            \put(25, 41){\small $2$}
            \put(25, 32){\small $1$}
            \put(25, 19){\small $0$}
            \put(25, 13){\small $-1$}
            \put(25, 4){\small $-2$}
            
            \put(0, 19){\small $-2$}
            \put(8, 19){\small $-1$}
            \put(34, 19){\small $1$}
            \put(41, 19){\small $2$}
            
            \put(55, 0){\includegraphics[height = 50mm]{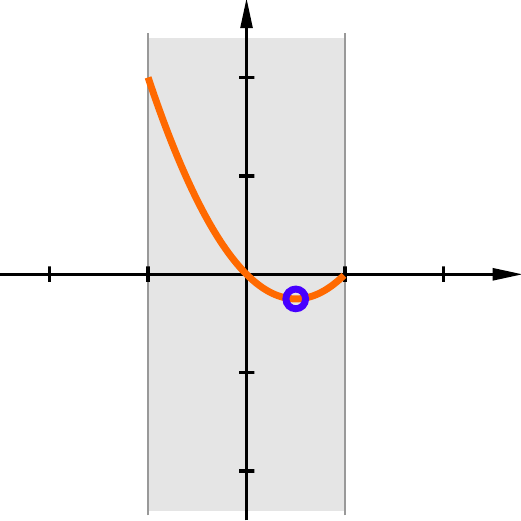}}
            
            \put(80, 48){\small $y$}
            \put(102, 19){\small $x$}
            
            \put(80, 41){\small $2$}
            \put(80, 32){\small $1$}
            \put(80, 19){\small $0$}
            \put(80, 13){\small $-1$}
            \put(80, 4){\small $-2$}
            
            \put(55, 19){\small $-2$}
            \put(63, 19){\small $-1$}
            \put(89, 19){\small $1$}
            \put(96, 19){\small $2$}
            
		    \put(110, 0){\includegraphics[height = 50mm]{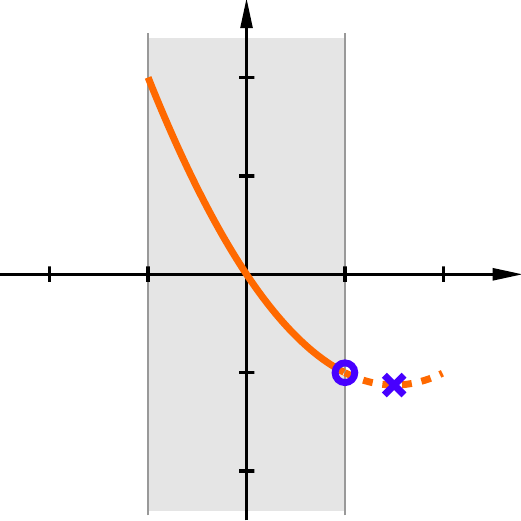}}
		    
		    \put(135, 48){\small $y$}
            \put(157, 19){\small $x$}
            
            \put(135, 41){\small $2$}
            \put(135, 32){\small $1$}
            \put(135, 19){\small $0$}
            \put(135, 13){\small $-1$}
            \put(135, 4){\small $-2$}
            
            \put(110, 19){\small $-2$}
            \put(118, 19){\small $-1$}
            \put(144, 19){\small $1$}
            \put(151, 19){\small $2$}
            
        \end{picture}

		\caption{Случаи положения (1), (2) и (3) для $x^{*}$.\\ Символом '$\times$' показана точка экстремума, а '$\circ$' --- точка, соответствующая $x ^ {*}$. \\На втором рисунке эти точки совпадают.}
		\label{fig:lemma7_example}
	
	\end{center}
\end{figure}

\noindent\textbf{Доказательство}

\noindentЗаметим, что 
\begin{equation*}
	f' \big( x_{\textit{extr}} \big) = 0, \qquad f''(x) = 2 a = \text{const} \neq 0,
\end{equation*}
то есть $x_{\textit{extr}}$ --- единственный корень уравнения $f'(x) = 0$.

Тогда если 
\begin{equation*}
	x_{\textit{extr}} \in \unitint,
\end{equation*}
\noindentто
\begin{equation*}
	\min\limits_{x \in \unitint} \, \big| f'(x) \big| = \big| f'(x_{\textit{extr}}) \big| = 0 \qquad \Rightarrow \qquad x^{*} = x_{\textit{extr}}.
\end{equation*}

Теперь перейдём к случаю $x_{\textit{extr}} \notin \unitint$.\\
Рассмотрим выражение
\begin{equation*}
	\begin{array}{l}
		\big| f'(x) \big|' = f'(x)\,f''(x) \, / \, \big| f'(x) \big| = 2a\,f'(x) \, / \, \big| f'(x) \big| = \\
		= 2a \, \text{sign} \big( f'(x) \big) = 2a \, \text{sign} \big( 2ax + b \big) = \\
		= 2 a \, \text{sign} \big( 2 a\,(x + 0.5 b \, / \, a) \big) = 2 a \,\text{sign} \big(2 a\,(x - x_{\textit{extr}}) \big) = \\
		= 2 a \, \text{sign} (2 a) \, \text{sign} (x - x_{\textit{extr}}) = 2 |a| \, \text{sign}(x - x_{\textit{extr}}).
	\end{array}
\end{equation*}

Пусть 
\begin{equation*}
	x_{\textit{extr}} < -1,
\end{equation*}
тогда так как $x \in \unitint$, то
\begin{equation*}
	\big| f'(x) \big|' = 2 |a| > 0.
\end{equation*}
Значит, 
\begin{equation*}
	\min\limits_{x \in \unitint} \, \big| f'(x) \big| = \big| f'(-1) \big| \qquad \Rightarrow \qquad x^{*} = -1.
\end{equation*}

Пусть теперь
\begin{equation*}
	x_{\textit{extr}} > 1,
\end{equation*}
тогда так как $x \in \unitint$, то
\begin{equation*}
	\big| f'(x) \big|' = - 2 |a| < 0.
\end{equation*}
Значит,
\begin{equation*}
	\min\limits_{x \in \unitint} \, \big| f'(x) \big| = \big| f'(1) \big| \qquad \Rightarrow \qquad x^{*} = 1.
\end{equation*}
\hfill $\blacksquare$

\begin{lemma}[о выборе опорной точки приближения дробно-линейного выражения одной переменной на $\unitint$]
\end{lemma}
Пусть
\begin{equation*}
	\begin{array}{l}
		f(x) = (ax + b) \, / \, (cx + d), \qquad x \in \unitint, \\ [1mm]
		|d| > |c| > 0, \qquad |a| + |b| > 0, \\ [1mm]
		\min\limits_{ \, x \in \unitint \, } \big| f'(x) \big| = \big| f'(x^{*}) \big|.
	\end{array}
\end{equation*}
Тогда
\begin{equation*}
	x^{*} = \text{sign} ( d \, / \, c ).
\end{equation*}
\textbf{Доказательство}\\
Заметим, что из
\begin{equation*}
	|d| > |c| > 0
\end{equation*}
следует, что
\begin{equation*}
	\forall x \in \unitint \qquad 0 \neq c\,x + d.
\end{equation*}

Рассмотрим выражение
\begin{equation*}
	\begin{array}{l}
		\big| f'(x) \big|' = f'(x)\,f''(x) \, / \, \big| f'(x) \big| =\\
		= (ad - bc) \, / \, (cx + d) ^ {2} \, 2c \, (bc - ad) \, / \, (cx + d) ^ {3} \\
		\hspace{2cm} / \, \big| (ad - bc) \, / \, (cx + d) ^ {2} \big| = \\
		= -2c \, (ad - bc) ^ {2} \, / \, (cx + d) ^ {5} \, / \, | ad - bc |\,\big| (cx + d) ^ {2} \big| =  \\
		= - 2c \, |ad - bc| \, / \, (cx + d) ^ {3}.
	\end{array}
\end{equation*}

\noindentЗаметим, что $\big| f'(x) \big|'$ является знакопостоянной при $x \in \unitint$.

\noindentРассмотрим выражение
\begin{equation*}
	\begin{array}{l}
		\text{sign} \big| f'(x) \big|' = \text{sign} ( -2c \, | ad - bc | \, / \, (cx + d) ^ {3} ) = \\
		= - \text{sign} ( c \, / \, (cx + d) ^ {3} \big) = - \text{sign} \big( c \, / \, (cx + d) ) =  \\
		= - \text{sign} ( 1 \, / \, (x + d \, / \, c) ) = - \text{sign} (x + d \, / \, c).
	\end{array}
\end{equation*}

Если
\begin{equation*}
	d \, / \, c > 0,
\end{equation*}
то
\begin{equation*}
	\forall x \in \unitint \qquad \big| f'(x) \big|' < 0.
\end{equation*}
Значит,
\begin{equation*}
	\min\limits_{ \, x \in \unitint \, } \big| f'(x) \big| = \big| f'(1) \big| \qquad \Rightarrow \qquad x^{*} = 1.
\end{equation*}

Перейдём к случаю
\begin{equation*}
	d \, / \, c < 0.
\end{equation*}
Тогда получим, что
\begin{equation*}
	\forall x \in \unitint \qquad \big| f'(x) \big|' > 0.
\end{equation*}

\noindentТаким образом заключаем, что
$$\min_{ \, x \in \unitint \, } \big| f'(x) \big| = \big| f'(-1) \big| \qquad \Rightarrow \qquad x^{*} = -1.$$

Объединяя случаи, получаем, что 
$$x^{*} = \text{sign}(d \, / \, c).$$

$\hfill \blacksquare$

\section{Арифметика линейных функциональных интервалов}

\subsection{Принцип построения арифметических операций}

В интервальном анализе используются несколько принципов построения операций между интервалами \cite{Shary}. Один из них: \textit{результирующий интервал операции должен содержать в себе всевозможные результаты применения этой операции для любых представителей интервалов-операндов}. Будем строить арифметику линейных функциональных интервалов, опираясь на него.

Пусть $\mbf{x}, \mbf{y} \in \mathbb{L}\mathbb{F}\mathbb{R}(x)$, то есть
\begin{equation*}
	\begin{array}{c}
		 \mbf{x} = \big[ \, \underline{a} x + \underline{b}, \, \overline{a} x + \overline{b} \, \big], \qquad x \in \unitint,\\
		 \mbf{y} = \big[ \, \underline{c} x + \underline{d}, \, \overline{c} x + \overline{d} \, \big], \qquad x \in \unitint.
	\end{array}
\end{equation*}

\subsection{Сложение}

\noindent\underline{\textit{Сложение}}

\begin{equation*}
	\begin{array}{r c l}
		 \mbf{x} + \mbf{y} & = &\big[ \, \underline{a} x + \underline{b}, \, \overline{a} x + \overline{b} \, \big] + \big[ \, \underline{c} x + \underline{d}, \, \overline{c} x + \overline{d} \big] = \\ [1mm]
		 & = & \Big[ F_{\downarrow} \big\{ (\underline{a} x + \underline{b}) + (\underline{c} x + \underline{d}), \, (\underline{a} x + \underline{b}) + (\overline{c} x + \overline{d}), \\
		 & & \hspace{1cm} (\overline{a} x + \overline{b}) + (\underline{c} x + \underline{d}), \, (\overline{a} x + \overline{b}) + (\overline{c} x + \overline{d}) \big\}, \\
		 & & \hspace{0.25cm} F_{\uparrow} \big\{ (\underline{a} x + \underline{b}) + (\underline{c} x + \underline{d}), \, (\underline{a} x + \underline{b}) + (\overline{c} x + \overline{d}), \\
		 & & \hspace{1cm} (\overline{a} x + \overline{b}) + (\underline{c} x + \underline{d}), \, (\overline{a} x + \overline{b}) + (\overline{c} x + \overline{d}) \big\} \Big] = \\
		 & = & \Big[ F_{\downarrow} \big\{ (\underline{a} + \underline{c}) \, x + (\underline{b} + \underline{d}), \, (\underline{a} + \overline{c}) \, x + (\underline{b} + \overline{d}), \\
		 & & \hspace{1cm} (\overline{a} + \underline{c}) \, x + (\overline{b} + \underline{d}), \, (\overline{a} + \overline{c}) \, x + (\underline{b} + \underline{d})\big\}, \\
		 & & \hspace{0.25cm} F_{\uparrow} \big\{ (\underline{a} + \underline{c}) \, x + (\underline{b} + \underline{d}), \, (\underline{a} + \overline{c}) \, x + (\underline{b} + \overline{d}), \\
		 & & \hspace{1cm} (\overline{a} + \underline{c}) \, x + (\overline{b} + \underline{d}), \, (\overline{a} + \overline{c}) \, x + (\underline{b} + \underline{d})\big\} \Big] = \\
		 & = & \big[ \, (\underline{a} + \underline{c}) \, x + (\underline{b} + \underline{d}), \, (\overline{a} + \overline{c}) \, x + (\overline{b} + \overline{d}) \, \big].
	\end{array}	
\end{equation*}

Результат процедур $F_{\downarrow}$ и $F_{\uparrow}$ был получен из интерпретации однопараметрического функционального интервала, как семейства классических интервалов. Напомним, что для этих интервалов операция сложения определена через сложение соответствующих концов операндов.

Покажем эффективность проведения операции сложения интервалов семейства $\mathbb{L}\mathbb{F}\mathbb{R}(x)$ с помощью следующей леммы.

\begin{lemma}[об эффективности применения операции сложения в функциональной интервальной арифметике]
\end{lemma}

\noindentПусть
\begin{equation*}
    \begin{array}{cc}
        \mbf{a}, \mbf{b} \in \mathbb{I}\mathbb{R}, \qquad \mbf{c}, \mbf{d} \in \mathbb{L}\mathbb{F}\mathbb{R}(x), \\
        \text{ran} \, \mbf{a} = \text{ran} \, \mbf{c}, \qquad \text{ran} \, \mbf{b} = \text{ran} \, \mbf{d}.
    \end{array}
\end{equation*}
Тогда
\begin{equation*}
    \text{ran} (\mbf{c} + \mbf{d}) \subseteq \text{ran} (\mbf{a} + \mbf{b}),
\end{equation*}
причём 
\begin{equation*}
    \text{ran}(\mbf{c} + \mbf{d}) = \text{ran}(\mbf{a} + \mbf{b}) \qquad \Leftrightarrow \qquad \underline{\mbf{c}}' \, \underline{\mbf{d}}' \geq 0 \text{ и } \overline{\mbf{c}}' \, \overline{\mbf{d}}' \geq 0.
\end{equation*}
\textbf{Доказательство.}
Так как $\mbf{a}, \mbf{b} \in \mathbb{I}\mathbb{R}$, то
\begin{equation*}
    \mbf{a} = [ \, \underline{a}, \, \overline{a} \, ], \qquad \mbf{b} = [ \, \underline{b}, \, \overline{b} \, ].
\end{equation*}
Так как $\mbf{c}, \mbf{d} \in \mathbb{L}\mathbb{F}\mathbb{R}(x)$, то
\begin{equation*}
    \begin{array}{lcl}
        \mbf{c} = [ \, \underline{c_{1}} x + \underline{c_{0}}, \, \overline{c_{1}} x + \overline{c_{0}} \, ], & \qquad & x \in \unitint, \\
        \mbf{d} = [ \, \underline{d_{1}} x + \underline{d_{0}}, \, \overline{d_{1}} x + \overline{d_{0}} \, ], & \qquad & x \in \unitint.
    \end{array}
\end{equation*}

По определению
\begin{equation*}
    \begin{array}{l}
        \mbf{a} + \mbf{b} = [ \, \underline{a} + \underline{b}, \, \overline{a} + \overline{b} \, ], \\
        \mbf{c} + \mbf{d} = \big[ \, (\underline{c_{1}} + \underline{d_{1}}) \, x + (\underline{c_{0}} + \underline{d_{0}}), \, (\overline{c_{1}} + \overline{d_{1}}) \, x + (\overline{c_{0}} + \overline{d_{0}}) \, \big].
    \end{array}
\end{equation*}

Так как
\begin{equation*}
    \begin{array}{cc}
        \text{ran} \, \mbf{a} = \text{ran} \, \mbf{c}, \\
        \text{ran} \, \mbf{b} = \text{ran} \, \mbf{d},
    \end{array}
\end{equation*}
то
\begin{equation*}
    \begin{array}{cc}
        \underline{a} = \underline{c_{0}} - |\underline{c_{1}}|, \qquad \overline{a} = \overline{c_{0}} + |\overline{c_{1}}|, \\
        \underline{b} = \underline{d_{0}} - |\underline{d_{1}}|, \qquad \overline{b} = \overline{d_{0}} + |\overline{d_{1}}|.
    \end{array}
\end{equation*}

Подставим эти соотношения
\begin{equation*}
    \begin{array}{r c l}
        \mbf{a} + \mbf{b}& = & \big[ \, \underline{c_{0}} - |\underline{c_{1}}| + \underline{d_{0}} - |\underline{d_{1}}|, \, \overline{c_{0}} + |\overline{c_{1}}| + \overline{d_{0}} + |\overline{d_{1}}| \, \big] = \\
        & = & \Big[ \, (\underline{c_{0}} + \underline{d_{0}}) - \big( |\underline{c_{1}}| + |\underline{d_{1}}| \big), \, (\overline{c_{0}} + \overline{d_{0}}) + \big( |\overline{c_{1}}| + |\overline{d_{1}}| \big) \, \Big].
    \end{array}
\end{equation*}

Таким образом,
\begin{equation*}
    \begin{array}{r c l}
        \text{ran} (\mbf{a} + \mbf{b}) & = & \Big[ \, (\underline{c_{0}} + \underline{d_{0}}) - \big( |\underline{c_{1}}| + |\underline{d_{1}}| \big), \, (\overline{c_{0}} + \overline{d_{0}}) + \big( |\overline{c_{1}}| + |\overline{d_{1}}| \big) \, \Big], \\
        \text{ran} (\mbf{c} + \mbf{d}) & = & \big[ \, (\underline{c_{0}} + \underline{d_{0}}) - |\underline{c_{1}} + \underline{d_{1}}|, \, (\overline{c_{0}} + \overline{d_{0}}) + |\overline{c_{1}} + \overline{d_{1}}| \, \big].
    \end{array}
\end{equation*}

Рассмотрим выражение $\underline{\text{ran}(\mbf{c} + \mbf{d})} - \underline{\text{ran}(\mbf{a} + \mbf{b})}$:
\begin{equation*}
    \begin{array}{l}
        \big( (\underline{c_{0}} + \underline{d_{0}}) - |\underline{c_{1}} + \underline{d_{1}}| \big) - \Big( (\underline{c_{0}} + \underline{d_{0}}) - \big( |\underline{c_{1}}| + |\underline{d_{1}}| \big)  \Big) = \\
        = (\underline{c_{0}} + \underline{d_{0}}) - |\underline{c_{1}} + \underline{d_{1}}| - (\underline{c_{0}} + \underline{d_{0}}) + |\underline{c_{1}}| + |\underline{d_{1}}| = \\
        = |\underline{c_{1}}| + |\underline{d_{1}}| - |\underline{c_{1}} + \underline{d_{1}}| \geq 0, \\
        \hspace{1cm}\text{причём }|\underline{c_{1}}| + |\underline{d_{1}}| = |\underline{c_{1}} + \underline{d_{1}}| \Leftrightarrow \underline{c_{1}} \, \underline{d_{1}} \geq 0.
    \end{array}
\end{equation*}
То есть, $\underline{\text{ran}(\mbf{c} + \mbf{d})} \geq \underline{\text{ran}(\mbf{a} + \mbf{b})}$, а равенство достигается при условии $\underline{\mbf{c}}' \, \underline{\mbf{d}}' \geq 0$.

Аналогично рассмотрим выражение $\overline{\text{ran}(\mbf{a} + \mbf{b})} - \overline{\text{ran}(\mbf{c} + \mbf{d})}$:
\begin{equation*}
    \begin{array}{l}
        \Big( (\overline{c_{0}} + \overline{d_{0}}) + \big( |\overline{c_{1}}| + |\overline{d_{1}}| \big) \Big) - \big( (\overline{c_{0}} + \overline{d_{0}}) + |\overline{c_{1}} + \overline{d_{1}}| \big) = \\
        = (\overline{c_{0}} + \overline{d_{0}}) + |\overline{c_{1}}| + |\overline{d_{1}}| - (\overline{c_{0}} + \overline{d_{0}}) - |\overline{c_{1}} + \overline{d_{1}}| = \\
        = |\overline{c_{1}}| + |\overline{d_{1}}| - |\overline{c_{1}} + \overline{d_{1}}| \geq 0, \\
        \hspace{1cm} \text{причём }|\overline{c_{1}}| + |\overline{d_{1}}| = |\overline{c_{1}} + \overline{d_{1}}| \Leftrightarrow \overline{c_{1}} \, \overline{d_{1}} \geq 0. 
    \end{array}
\end{equation*}
То есть, $\overline{\text{ran}(\mbf{a} + \mbf{b})} \geq \overline{\text{ran}(\mbf{c} + \mbf{d})}$, а равенство достигается при условии $\overline{\mbf{c}}' \, \overline{\mbf{d}}' \geq 0$.

Таким образом, 
\begin{equation*}
    \text{ran} (\mbf{c} + \mbf{d}) \subseteq \text{ran} (\mbf{a} + \mbf{b}),
\end{equation*}
причём 
\begin{equation*}
    \text{ran} (\mbf{c} + \mbf{d}) = \text{ran} (\mbf{a} + \mbf{b}) \Leftrightarrow \underline{\mbf{c}}' \, \underline{\mbf{d}}' \geq 0 \text{ и } \overline{\mbf{c}}' \, \overline{\mbf{d}}' \geq 0.
\end{equation*} \hfill $\blacksquare$

\subsection{Унарный минус}

\begin{equation*}
	\begin{array}{r c l}
		 - \mbf{x} & = & -\big[ \, \underline{a} x + \underline{b}, \, \overline{a} x + \overline{b} \, \big] = \big[ \, -\overline{a} x - \overline{b}, \, -\underline{a} x - \underline{b} \, \big].
	\end{array}
\end{equation*}	

Результат операции унарного минуса также был получен из интерпретации линейного функционального интервала, как параметрического семейства классических интервалов.

\subsection{Вычитание}

Операция вычитания получается путём комбинирования рассмотренных выше операций унарного минуса и сложения.

\begin{equation*}
	\begin{array}{r c l}
		 \mbf{x} - \mbf{y} & = & \big[ \, \underline{a} x + \underline{b}, \, \overline{a} x + \overline{b} \, \big] - \big[ \, \underline{c} x + \underline{d}, \, \overline{c} x + \overline{d} \big] = \\
		  & = & \big[ \, \underline{a} x + \underline{b}, \, \overline{a} x + \overline{b} \, \big] + \big(- [ \, \underline{c} x + \underline{d}, \, \overline{c} x + \overline{d} ] \big) = \\
		 & = & \big[ \, (\underline{a} - \overline{c}) x\,+ (\underline{b} - \overline{d}), \, (\overline{a} - \underline{c}) x + (\overline{b} - \underline{d}) \, \big].
	\end{array}	
\end{equation*}

\subsection{Умножение}

По базовому принципу построения операций имеем

\begin{equation*}
	\begin{array}{r c l}
		 \mbf{x} \cdot \mbf{y} & = & \big[ \, \underline{a} x + \underline{b}, \, \overline{a} x + \overline{b} \, \big]\,\big[ \, \underline{c} x + \underline{d}, \, \overline{c} x + \overline{d} \, \big] = \\ [1mm]
		 & = & \Big[ \, F_{\downarrow} \big\{ \, (\underline{a} x + \underline{b})\,(\underline{c} x + \underline{d}), \, (\underline{a} x + \underline{b})\,(\overline{c} x + \overline{d}), \\
		 & & \hspace{1cm} (\overline{a} x + \overline{b})\,(\underline{c} x + \underline{d}), \, (\overline{a} x + \overline{b})\,(\overline{c} x + \overline{d}) \, \big\}, \\
		 & & \hspace{0.25cm} F_{\uparrow} \big\{ \, (\underline{a} x + \underline{b})\,(\underline{c} x + \underline{d}), \, (\underline{a} x + \underline{b})\,(\overline{c} x + \overline{d}), \\
		 & & \hspace{1.25cm} (\overline{a} x + \underline{b})\,(\underline{c} x + \underline{d}), \, (\overline{a} x + \overline{b})\,(\overline{c} x + \overline{d}) \, \big\} \, \Big].
	\end{array}	
\end{equation*}

Видно, что все аргументы функций $F_{\downarrow}$ и $F_{\uparrow}$ представляют собой квадратичные полиномы. По Леммам 1 и 2 нам известны способы построения этих функций для класса линейных функций. Рассмотрим, каким образом эти аргументы можно привести к этому классу.

Рассмотрим квадратичный полином
\begin{equation*}
	f(x) = a x^{2} + b x + c, \qquad x \in \unitint.
\end{equation*}

\noindentВ силу того, что вторая производная является знакопостоянной $\big( f''(x) = 2 a = \text{const} \big)$, функция является выпуклой или вогнутой.

Для начала предположим, что
\begin{equation*}
	f''(x) = 2 a \geq 0.
\end{equation*}

\noindentТак как $x \in \unitint$, можем воспользоваться Леммами 3 и 5. Тогда получим, что
\begin{equation*}
	\begin{array}{l}
		f(x) \geq f(x^{*}) + (x - x^{*})\,f'(x^{*}),\text{ где }\min\limits_{x \in \unitint} \, \big| f'(x) \big| = \big| f'(x^{*}) \big|, \\
		f(x) \leq \frac{x}{2}\,\big( f(1) - f(-1) \big) + \frac{1}{2}\,\big( f(1) + f(-1) \big). 
	\end{array}
\end{equation*}

Последовательно для каждого аргумента, рассмотрим процедуру $F_{\downarrow}$, входящую в операцию умножения.

\begin{equation*}
	\begin{array}{l}
		F_{\downarrow} \big\{ \, \dots, \, (A x + B)\,(C x + D), \, \dots \, \big\} \geq \\
		\hspace{2cm} \geq F_{\downarrow} \Big\{ \, \dots, \, \big( A x^{*} + B \big)\,\big( C x^{*}\,+ D \big)\\
		\hspace{4cm} + \, \big( x - x^{*} \big)\,\big( A\,(2 C x^{*} + D) + B C \big) , \, \dots \Big\} = \\
		= F_{\downarrow} \big\{ \, \dots, \, AC x^{*2} + ADx^{*} +BCx^{*} + BD + 2ACxx^{*} + ADx \\
		\hspace{2cm} + \, BCx - 2ACx^{*2} - ADx^{*} - BCx^{*} \big\} = \\
		= F_{\downarrow} \big\{ \, \dots, \, -AC x^{*2} + BD + 2ACxx^{*} + ADx + BCx, \, \dots \big\} =  \\
		= F_{\downarrow} \big\{ \, \dots, \, x \, (2ACx^{*} + AD + BC) + (BD -AC x^{*2}), \, \dots \big\},
	\end{array}
\end{equation*}
где
\begin{equation*}
	(A, \, B) \in \big\{ \, ( \underline{a}, \, \underline{b} ), \, ( \overline{a}, \, \overline{b} ) \, \big\},
	\qquad 
	(C, \, D) \in \big\{ \, (\underline{c}, \, \underline{d} ), \, ( \overline{c}, \, \overline{d} ) \, \big\}.
\end{equation*}

Аналогично рассмотрим процедуру $F_{\uparrow}$.
\begin{equation*}
	\begin{array}{l}
		F_{\uparrow} \big\{ \, \dots, \, (A x + B)\,(C x + D), \, \dots \, \big\} \leq \\
		\hspace{1.5cm} \leq F_{\uparrow} \Big\{ \, \dots, \, \frac{x}{2}\,\big( (A + B)\,(C + D) - (-A + B)\,(-C + D) \big) \\
		\hspace{3cm} + \, \frac{1}{2}\,\big( (A + B)\,(C + D) + (-A + B)\,(-C + D) \big), \, \dots \Big\} = \\
		= F_{\uparrow} \big\{ \, \dots, \, \frac{x}{2} \, ( AC + AD + BC + BD - AC + AD + BC - BD ) \\
		\hspace{1.5cm} + \, \frac{1}{2} \, ( AC + AD + BC + BD + AC - AD - BC + BD ), \, \dots \big\} = \\
		= F_{\uparrow} \big\{ \, \dots, \, \frac{x}{2} \, ( 2AD + 2BC \big) + \frac{1}{2} \, \big( 2AC + 2BD ), \, \dots \big\} \\
		= F_{\uparrow} \big\{ \, \dots, \, x \, ( AD + BC ) + ( AC + BD ), \, \dots \big\},
	\end{array}
\end{equation*}
где
\begin{equation*}
	\begin{array}{l}
		(A, \, B) \in \big\{ \, (\underline{a}, \, \underline{b} ), \, (\overline{a}, \, \overline{b}) \, \big\}, \qquad
		(C, \, D) \in \big\{ \, (\underline{c}, \, \underline{d}), \, (\overline{c}, \, \overline{d}) \, \big\}.
	\end{array}
\end{equation*}

Для случая
\begin{equation*}
	f''(x) = 2 a \leq 0
\end{equation*}
с помощью Лемм 4 и 6 получим аналогично
\begin{equation*}
	\begin{array}{l}
		F_{\downarrow} \big\{ \, \dots, \, (A x + B)\,(C x + D), \, \dots \, \big\} \geq \\
		\hspace{2cm} \geq F_{\downarrow} \big\{ \, \dots, \, x \, ( AD + BC ) + ( AC + BD ), \, \dots \big\}, \\
		F_{\uparrow} \big\{ \, \dots, \, (A x + B)\,(C x + D), \, \dots \, \big\} \leq \\
		\hspace{2cm} \leq F_{\uparrow} \big\{ \, \dots, \, x (2ACx^{*} + AD + BC) + (BD -AC x^{*2}), \, \dots \big\},
	\end{array}
\end{equation*}
где
\begin{equation*}
	(A, \, B) \in \big\{ \, ( \underline{a}, \, \underline{b} ), \, ( \overline{a}, \, \overline{b} ) \, \big\},
	\qquad (C, \, D) \in \big\{ \, (\underline{c}, \, \underline{d} ), \, ( \overline{c}, \, \overline{d} ) \, \big\}.
\end{equation*} 

После приведения аргументов процедур $F_{\downarrow}$ и $F_{\uparrow}$ от квадратичных полиномов к линейным, можем получить результат действия этих функций на них.

Для этого воспользуемся Леммой 1 о верхней оценке линейных функций на $\unitint$. Применяя это построение сначала для аргументов 1 и 2, 3 и 4, а затем для полученных результатов, получим общую огибающую для всех аргументов сверху.

Аналогично можем построить нижние огибающие, используя Лемму 2 о нижней оценке линейных функций на $\unitint$.

\begin{lemma}[об эффективности операции умножения в функциональной интервальной арифметике]
\end{lemma}

\noindentПусть
\begin{equation*}
    \begin{array}{c}
        \mbf{a}, \mbf{b} \in \mathbb{I}\mathbb{R}, \qquad \mbf{c}, \mbf{d} \in \mathbb{L}\mathbb{F}\mathbb{R}(x),\\
        \text{ran} \, \mbf{a} = \text{ran} \, \mbf{c}, \qquad \text{ran} \, \mbf{b} = \text{ran} \, \mbf{d}.
    \end{array}
\end{equation*}
Тогда
\begin{equation*}
    \text{ran} (\mbf{a} \cdot \mbf{b}) = \text{ran} (\mbf{c} \cdot \mbf{d}).
\end{equation*}
\textbf{Доказательство}
По базовому принципу построения операций имеем
\begin{equation*}
	\begin{array}{r c l}
		 \mbf{x} \cdot \mbf{y} & = & \big[ \, \underline{a} x + \underline{b}, \, \overline{a} x + \overline{b} \, \big]\,\big[ \, \underline{c} x + \underline{d}, \, \overline{c} x + \overline{d} \, \big] = \\
		 & = & \Big[ \, F_{\downarrow} \big\{ \, (\underline{a} x + \underline{b})\,(\underline{c} x + \underline{d}), \, (\underline{a} x + \underline{b})\,(\overline{c} x + \overline{d}), \\
		 & & \hspace{1cm} (\overline{a} x + \overline{b})\,(\underline{c} x + \underline{d}), \, (\overline{a} x + \overline{b})\,(\overline{c} x + \overline{d}) \, \big\}, \\
		 & & \hspace{0.25cm} F_{\uparrow} \big\{ \, (\underline{a} x + \underline{b})\,(\underline{c} x + \underline{d}), \, (\underline{a} x + \underline{b})\,(\overline{c} x + \overline{d}), \\
		 & & \hspace{1.25cm} (\overline{a} x + \underline{b})\,(\underline{c} x + \underline{d}), \, (\overline{a} x + \overline{b})\,(\overline{c} x + \overline{d}) \, \big\} \, \Big].
	\end{array}	
\end{equation*}

Каждый аргумент функций $F_{\downarrow}$ и $F_\uparrow$ представляет собой квадратичный полином. 

Начнём переводить аргументы в класс линейных функций. Если очередной аргумент представляет собой квадратичный полином с $f''(x) > 0$, то воспользовавшись следствием из Лемм 3 и 5 получим, что область значений не увеличится. Аналогично, если $f''(x) < 0$, то по следствиям из Лемм 4 и 6 также получим, что область значений не увеличится.

Осталось заметить, что если мы начнём строить результаты функций $F_{\downarrow}$ и $F_{\uparrow}$ для линейных функций с помощью Лемм 1 и 2, то мы также не будем увеличивать область значений.

Исходя из приведённых выше рассуждений получим, что $\text{ran} (\mbf{a} \cdot \mbf{b}) = \text{ran} (\mbf{c} \cdot \mbf{d})$.

\hfill $\blacksquare$

\subsection{Деление}

По базовому принципу построения операций имеем

\begin{equation*}
	\begin{array}{r c l}
	    \mbf{x} \, / \, \mbf{y} & = &\big[ \, \underline{a} x + \underline{b}, \, \overline{a} x + \overline{b} \, \big] \, / \, \big[ \, \underline{c} x + \underline{d}, \, \overline{c} x + \overline{d} \big] = \\
		 & = & \Big[ \, F_{\downarrow} \big\{ \, (\underline{a} x + \underline{b}) \, / \, (\underline{c} x + \underline{d}), \, (\underline{a} x + \underline{b}) \, / \, (\overline{c} x + \overline{d}), \\
		 & & \hspace{1cm} (\overline{a} x + \overline{b}) \, / \, (\underline{c} x + \underline{d}), \, (\overline{a} x + \overline{b}) \, / \, (\overline{c} x + \overline{d}) \, \big\}, \\
		 & & \hspace{0.25cm} F_{\uparrow} \big\{ \, (\underline{a} x + \underline{b}) \, / \, (\underline{c} x + \underline{d}), \, (\underline{a} x + \underline{b}) \, / \, (\overline{c} x + \overline{d}), \\
		 & & \hspace{1.25cm} (\overline{a} x + \underline{b}) \, / \, (\underline{c} x + \underline{d}), \, (\overline{a} x + \overline{b}) \, / \, (\overline{c} x + \overline{d}) \, \big\} \, \Big].
	\end{array}	
\end{equation*}

В данной операции все аргументы процедур $F_{\downarrow}$ и $F_{\uparrow}$ представляют собой выражения вида
\begin{equation*}
	f(x) = (a x + b) \, / \, (c x + d), \qquad x \notin [ \, d - c, \, d + c \, ].
\end{equation*}

Заметим, что в силу того что
\begin{equation*}
	f''(x) = 2 c \, (b c - a d) \, / \, (c x + d) ^ {3} = \text{const} \, / \, (c x + d) ^ {3}
\end{equation*}
и в силу условия $x \notin [ \, d - c, \, d + c \, ]$ получаем, что $f''(x)$ является знакопостоянной. Значит, $f(x)$ является выпуклой или вогнутой функцией.

Предположим, что 
\begin{equation*}
	f''(x) = 2 c \, (b c - a d) \geq 0.
\end{equation*}
Тогда для функции $f(x)$ также применимы Леммы 3 и 5 и верно
\begin{equation*}
	\begin{array}{l}
		f(x) \geq f(x^{*}) + (x - x^{*})\,f'(x^{*}),\text{ где }\text{min} \, |f'(x)| \, _{x \in \unitint} = |f'(x^{*})|, \\
		f(x) \leq \frac{x}{2}\,\big( f(1) - f(-1) \big) + \frac{1}{2}\,\big( f(1) + f(-1) \big). 
	\end{array}
\end{equation*}

Рассмотрим результат применения функции для каждого аргумента $F_{\uparrow}$, входящих в операцию деления:

\begin{equation*}
	\begin{array}{l}
		F_{\uparrow} \big\{ \, \dots, \, (A x + B) \, / \, (C x + D), \, \dots \, \big\} \leq \\
		\hspace{1cm} \leq F_{\uparrow} \Big\{ \, \dots, \, \frac{x}{2}\,\big( (A + B) \, / \, (C + D) - (-A + B) \, / \, (-C + D) \big) \\
		\hspace{2cm} + \, \frac{1}{2}\,\big( (A + B) \, / \, (C + D) + (-A + B) \, / \, (-C + D) \big), \, \dots \Big\} = \\
		= F_{\uparrow} \Big\{ \, \dots, \, \frac{x}{2} \, \big( (A + B) (C - D) - (A - B) (C + D) \big) \, / \, (C ^ {2} - D ^ {2}) \\
		\hspace{1cm} + \, \frac{1}{2}\,\big( (A + B) (C - D) + (A - B) \, (C + D) \big) \, / \, (C ^ {2} - D ^ {2}), \, \dots \Big\} = \\ 
		= F_{\uparrow} \Big\{ \, \dots, \, x \, (BC - AD) \, / \, (C ^ {2} - D ^ {2}) + \, (AC - BD) \, / \, (C ^ {2} - D ^ {2}), \, \dots \, \Big\},
	\end{array}
\end{equation*}
где
\begin{equation*}
	(A, \, B) \in \{ \, ( \underline{a}, \, \underline{b} ), \, ( \overline{a}, \, \overline{b} ) \, \},
	\qquad (C, \, D) \in \{ \, (\underline{c}, \, \underline{d} ), \, ( \overline{c}, \, \overline{d} ) \, \}.
\end{equation*} 

\noindentАналогично результат применения функции для каждого аргумента $F_{\downarrow}$, входящих в операцию деления:

\begin{equation*}
    \begin{array}{l}
		F_{\downarrow} \big\{ \, \dots, \, (A x + B) \, / \, (C x + D), \, \dots \, \big\} \geq \\
		\hspace{2cm} \geq F_{\downarrow} \Big\{ \, \dots, \, \big( A x^{*} + B \big) \, / \, \big( C x^{*}\,+ D \big)\\
		\hspace{4cm} + \, \big( x - x^{*} \big)\,\big( A D - B C \big) \, / \, (C x ^ {*} + D) ^ {2}, \, \dots \Big\} = \\
		= F_{\downarrow} \big\{ \, \dots, \, x \, (AD - BC) \, / \, (Cx^{*} + D) ^ {2} \\
		\hspace{2cm} + \, (ACx ^ {*2} + 2BCx ^ {*} + BD) \, / \, (Cx ^ {*} + D) ^ {2}, \, \dots \, \big\},
	\end{array}
\end{equation*}
где
\begin{equation*}
		(A, \, B) \in \big\{ \, ( \underline{a}, \, \underline{b} ), \, ( \overline{a}, \, \overline{b} ) \, \big\},
		\qquad (C, \, D) \in \big\{ \, (\underline{c}, \, \underline{d} ), \, ( \overline{c}, \, \overline{d} ) \, \big\}.
\end{equation*} 

В предположении, что
\begin{equation*}
	f''(x) = 2 c \, (b c - a d) \leq 0,
\end{equation*}
получим аналогичные результаты после применения лемм 4 и 6
\begin{equation*}
	\begin{array}{l}
		F_{\uparrow} \big\{ \, \dots, \, (A x + B) \, / \, (C x + D), \, \dots \, \big\} \geq  \\
		\hspace{2cm}  \geq F_{\uparrow} \Big\{ \, \dots, \, x \, (BC - AD) \, / \, (C ^ {2} - D ^ {2}) \\
		\hspace{4cm} + \, (AC - BD) \, / \, (C ^ {2} - D ^ {2}), \, \dots, \, \Big\}, \\
		F_{\downarrow} \big\{ \, \dots, \, (A x + B) \, / \, (C x + D), \, \dots \, \big\} \leq \\
		\hspace{2cm} \leq F_{\downarrow} \big\{ \, \dots, \, ( A x^{*} + B ) \, / \, ( C x^{*} + D )\\
		\hspace{4cm} + \, ( x - x^{*} ) \, ( A D - B C ) \, / \, (C x ^ {*} + D) ^ {2}, \, \dots \big\},
	\end{array}
\end{equation*}
где
\begin{equation*}
	(A, \, B) \in \big\{ \, ( \underline{a}, \, \underline{b} ), \, ( \overline{a}, \, \overline{b} ) \, \big\},
	\qquad (C, \, D) \in \big\{ \, (\underline{c}, \, \underline{d} ), \, ( \overline{c}, \, \overline{d} ) \, \big\}.
\end{equation*} 

\section{Приложения линейной \\функциональной арифметики}

\subsection{Задача решения уравнения}

\subsubsection{Определение интервала поиска корней}

Изначально при решении задачи доказательного численного нахождения корней уравнения на вещественной оси $\mathbb{I}\mathbb{R}$ необходимо выделить исходный интервал поиска. 

Для задания начального интервала поиска корней, можно воспользоваться следствием, приведённом в книге \cite{Prasolov}, из одной из известных теорем теории функций --- теоремы Руше \cite{Rouche}, \cite{Privalov}.\\
\textbf{Следствие из теоремы Руше.}\\
Пусть 
\begin{equation*}
	f(z) = z ^ {n} + a_{n - 1} z ^ {n - 1} + \dots + a_{1} z + a_{0}, \qquad a_{i} \in \mathbb{C}, \qquad z \in \mathbb{C}.
\end{equation*}
Тогда внутри круга комплексной плоскости, определяемого выражением
\begin{equation*}
	|z| \leq 1 + \max\limits_{i = 1,\dots, n} \, |a_{i}|
\end{equation*}
расположено ровно $n$ корней полинома $f$ (с учетом их кратностей).

Из этого, в частности, следует, что если мы рассмотрим полином
\begin{equation*}
	f(x) = x ^ {n} + a_{n - 1} x ^ {n - 1} + \dots + a_{1} x + a_{0}, \qquad a_{i} \in \mathbb{R},
\end{equation*} 
тогда получим, что все вещественные корни полинома по модулю не превосходят $1 + \max\limits_{i = 1, \dots, n} \, {|a_{i}|}$.

\subsubsection{Сжатие функционального интервала к нулевой части}

Предположим, что у нас в распоряжении имеется интервал $\mbf{x} \in \mathbb{L}\mathbb{F}\mathbb{R}(x)$. Этот интервал можно записать в виде
\begin{equation*}
    \mbf{x} = [ \, \underline{a} x + \underline{b}, \, \overline{a} x + \overline{b} \, ], \qquad x \in \unitint.
\end{equation*}

\begin{definition}
\textit{Нулевой частью} интервала $\mbf{y} \in \mathbb{L}\mathbb{F}\mathbb{R}(x)$ будем называть такой интервал $\mbf{r} \in \mathbb{I}\mathbb{R}$, что 
\begin{equation}
\label{eq:interval_root}
    \left\{\begin{array}{ll}
        \forall x \in \mbf{r} & 0 \in \mbf{y}(x), \\
        \forall x \in \unitint \setminus \mbf{r} & 0 \notin \mbf{y}(x).
    \end{array}\right.
\end{equation}
\end{definition} 

Нахождение нулевой части интервала $\mathbb{L}\mathbb{F}\mathbb{R}(x)$ соответствует решению системы неравенств
\begin{equation*}
    \left\{\begin{array}{l}
        \underline{a} x + \underline{b} \leq 0, \\
        \overline{a} x + \overline{b} \geq 0.
    \end{array}\right.
\end{equation*}

\noindent Решением первого неравенства системы будет
\begin{equation*}
    \left\{\begin{array}{cl}
        x \leq -\underline{b} \, / \, \underline{a}, & \text{если } \underline{a} > 0, \\
        x \geq -\underline{b} \, / \, \underline{a}, & \text{иначе }(\text{если } \underline{a} < 0).
    \end{array}\right.
\end{equation*}

Так как $x \in \unitint$, то итоговое решение первого неравенства системы будет выглядеть следующим образом
\begin{equation*}
    \mathbb{I}\mathbb{R} \ni \mbf{S}_{1} = 
    \left\{\begin{array}{cl}
        \varnothing, & \text{если } \underline{a} > 0 \text{ и } -\underline{b} \, / \, \underline{a} < -1, \\
        \big[ \, -1, \, \min \, \{ \, -\underline{b} \, / \underline{a}, \, 1 \, \} \, \big], & \text{если } \underline{a} > 0 \text{ и } -\underline{b} \, / \underline{a} \geq -1, \\
        \varnothing, & \text{если } \underline{a} < 0 \text{ и } -\underline{b} \, / \, \underline{a} > 1, \\
        \big[ \, \max \, \{ \, -\underline{b} \, / \underline{a}, \, -1 \, \}, \, 1 \, \big], & \text{если } \underline{a} < 0 \text{ и } -\underline{b} \, / \underline{a} \leq 1.
    \end{array}\right.
\end{equation*}
Решением второго неравенства системы будет
\begin{equation*}
    \left\{\begin{array}{cl}
        x \geq -\overline{b} \, / \, \overline{a}, & \text{если } \overline{a} > 0, \\
        x \leq -\overline{b} \, / \, \overline{a}, & \text{иначе }(\text{если } \overline{a} < 0).
    \end{array}\right.
\end{equation*}

Так как $x \in \unitint$, то итоговое решение первого неравенства системы будет выглядеть следующим образом
\begin{equation*}
    \mathbb{I}\mathbb{R} \ni \mbf{S}_{2} = \left\{\begin{array}{cl}
        \varnothing, & \text{если } \overline{a} > 0 \text{ и } -\overline{b} \, / \, \overline{a} > 1, \\
        \big[ \, \max \, \{ \, -\overline{b} \, / \, \overline{a}, \, -1 \, \}, \, 1 \, \big], & \text{если } \overline{a} > 0 \text{ и } -\overline{b} \, / \, \overline{a} \leq 1, \\
        \varnothing, & \text{если } \overline{a} < 0 \text{ и } -\overline{b} \, / \, \overline{a} < -1, \\
        \big[ \, -1, \, \min \, \{ \, -\overline{b} \, / \, \overline{a}, \, 1 \, \} \, \big], & \text{если } \overline{a} < 0 \text{ и } -\overline{b} \, / \, \overline{a} \geq -1.
    \end{array}\right.
\end{equation*}

Итоговым интервалом, в котором гарантированно содержатся вся нулевая часть интервала $\mbf{x}$, 
будет $\mbf{S}_{1} \cap \mbf{S}_{2}$. Если в результате такого пересечения мы получили 
$\varnothing$, то это означает, что на интервале $\unitint$ не содержится 
нулевая часть интервала $\mbf{x}$. 
  
  
\subsubsection{Алгоритм поиска корней}
  
В интервальном анализе для доказательного поиска вещественных нулей функций 
широко распространены алгоритма класса <<ветвлений-и-отсечений>> \cite{Shary}, 
\cite{HansenWalster}. Для описания алгоритма введём ряд вспомогательных понятий. 

\begin{definition} Рабочий список --- массив упорядоченных пар
вида $( \omega, \, \mbf{I} )$, отсортированный по возрастанию значения $\omega$, где\\
$\mbf{I} \in \mathbb{I}\mathbb{R}$ --- \textit{область определения функции},\\
$\omega \in \mathbb{R}$ --- \textit{характеристика интервала}.

\end{definition}

\noindentАвтором работы предлагается следующий алгоритм этого класса:
  
\bigskip   
\begin{small}
\noindent\textbf{Входные данные:}

\begin{singlespacing}
    \begin{enumerate}
    
        \item Функция для вычисления \textit{характеристики интервала} $\omega$.
        
        \item Исходный интервал для поиска нулей функции.
        
        \item \textit{Пороговое значение остановки алгоритма} $\varepsilon > 0$.
        
    \end{enumerate}
\end{singlespacing}

\noindent\textbf{Алгоритм:}

    \begin{singlespacing}
        \begin{enumerate}
        	
        	\item Вычислить характеристику $\omega$ для исходного интервала поиска нулей функции. Поместить в пустой рабочий список полученную пару. Перейти на шаг 2.
        
        	\item Если рабочий список пуст или условие остановки выполнено, перейти на шаг 11, иначе перейти на шаг 3.
        	
        	\item Извлечь из рабочего списка первую пару $A = (\mbf{a}, \, a)$. Перейти на шаг 4.
        
        	\item Найти интервальное расширение функции $f(\mbf{a})$. Перейти на шаг 5.
        	
        	\item Проверить, содержится ли $0$ в интервальном расширении $f(\mbf{a})$. Если не содержится, перейти на шаг 2, иначе перейти на шаг 6.
        	
        	\item Применить к $\mbf{a}$ процедуру сжатия функционального интервала к нулевой части. Перейти на шаг 7.
        	
        	\item Если после процедуры сжатия $\mbf{a} = \varnothing$, перейти на шаг 2, иначе перейти на шаг 8.
        	
        	\item Произвести дихотомию $\mbf{a} \rightarrow \big\{ \mbf{a}_{-} = [ \, \underline{\mbf{a}}, \, \text{mid} \, \mbf{a} \, ], \, \mbf{a}_{+} = [ \, \text{mid} \, \mbf{a}, \, \overline{\mbf{a}} \, ] \big\}$. Вставить рабочие интервалы $\big(\mbf{a}_{-}, \, \omega_{-} \big)$ и $\big( \mbf{a}_{+}, \, \omega_{+} \big)$ в рабочий список. Перейти на шаг 2.
        	
        	\item Вывести в качестве ответа список всех полученных интервалов. Алгоритм завершён.
        	
        \end{enumerate}
    \end{singlespacing}
\end{small}
\clearpage

\subsubsection{Применение теоремы Больцано-Коши\\к функциональному интервалу}

Для интервалов семейства $\mathbb{L}\mathbb{F}\mathbb{R}(x)$ есть возможность <<встроенного>> 
использования теоремы Больцано-Коши \cite{Zorich} для доказательства наличия нуля непрерывной 
функции на рассматриваемом интервале. 
  
\bigskip\noindent  
\textbf{Теорема Больцано-Коши.}

\noindentЕсли функция, непрерывная на отрезке, принимает на его концах значения разных знаков, то на отрезке есть точка, в которой функция обращается в ноль.
 
\begin{figure}[ht]
    \begin{center}
    
    \unitlength=1mm  
    
    \begin{picture}(70,70)		
    
        \put(0, 0){\includegraphics[width = 65mm]{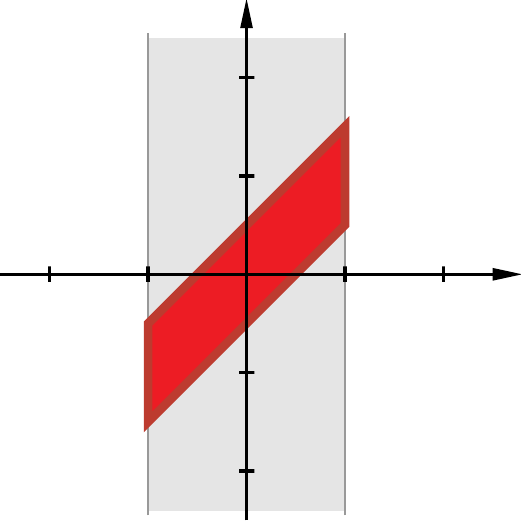}} 
        \put(46, 49){$y = x + 0.5$}
        \put(46, 36){$y = x - 0.5$}
        \put(34, 62){$y$}
        \put(62, 25){$x$}
        
        \put(34, 53){$2$}
        \put(34, 42){$1$}
        \put(34, 25){$0$}
        \put(34, 17){$-1$}
        \put(34, 5){$-2$}
        
        \put(1, 24){$-2$}
        \put(11, 24){$-1$}
        \put(44, 24){$1$}
        \put(54, 24){$2$}
        
    \end{picture} 	
	
	\caption{Пример интервала $\mathbb{L}\mathbb{F}\mathbb{R}(x)$, для которого \\ 
	применима теорема Больцано-Коши: $[ \, x - 0.5, \, x + 0.5 \, ].$}
	\label{fig:bolcano_example}
    	
    \end{center}
\end{figure}

\bigskip   
Предположим, что в результате нахождения интервального расширения функции $f(x)$ на некотором интервале $x \in [ \, l, \, r \, ]$ был получен интервал $\mbf{y}$ семейства $\mathbb{L}\mathbb{F}\mathbb{R}(x)$, у которого
\begin{equation*}
    \mbf{y}(-1) = [ \, l_{1}, \, l_{2} \, ], \qquad \mbf{y}(1) = [ \, r_{1}, \, r_{2} \, ], \qquad l_{1} \leq l_{2} \leq 0 \leq r_{1} \leq r_{2}.
\end{equation*}
Тогда по теореме Больцано-Коши на интервале $[ \, l, \, r \, ]$ существует корень уравнения $f(x) = 0$ (рис. \ref{fig:bolcano_example}).

Аналогично применима теорема для случая $l_{1} \geq l_{2} \geq 0 \geq r_{1} \geq r_{2}$.

\subsection{Задача нахождения минимума функции}

\subsubsection{Процедура сжатия интервала к минимуму}

Предположим, что у нас в распоряжении имеется интервал $\mbf{x} \in \mathbb{L}\mathbb{F}\mathbb{R}(x)$. Этот интервал можно записать в виде
\begin{equation*}
\mbf{x} = [\,\underline{a} x + \underline{b}, \,\overline{a} x + \overline{b}\,], 
   \qquad x \in \unitint. 
\end{equation*}

Пусть также известна оценка минимума сверху, которую обозначим $M$. Нужно убедиться, 
что $M \leq \min \Big\{ \, \overline{\mbf{x}(-1)}, \,  \overline{\mbf{x}(1)} \Big\}$, 
поскольку иначе именно $\min \Big\{\, \overline{\mbf{x}(-1)}, \,\overline{\mbf{x}(1)} 
\Big\}$ будет являться оценкой минимума сверху. 
    
Таким образом, нам необходимо найти такой интервал $\mbf{m} \in \mathbb{I}\mathbb{R}$, 
что 
\begin{equation*}
    \left\{\begin{array}{ll}
        \forall x \in \mbf{m} & M \in \mbf{x},  \\
        \forall x \in \unitint \setminus \mbf{m} & M \notin \mbf{x}.
    \end{array}\right.
\end{equation*}

Исходя из определения (\ref{eq:interval_root}) видно, что задачу нахождения минимума интервала $\mbf{x}$ можно свести к применению процедуры сжатия интервала к нулевой части интервала $(\mbf{x} - M)$.

\subsubsection{Алгоритм поиска минимума функции}
\label{min_alg}

В интервальном анализе для доказательного поиска глобального минимума функций широко распространены алгоритма класса <<ветвей-и-границ>> \cite{Shary}, \cite{HansenWalster}. Для описания алгоритма введём вспомогательные понятия.

\begin{definition} Рабочий список --- массив упорядоченных пар $( \mbf{I} , \, \omega )$, отсортированный по возрастанию значений $\omega$,
где \\
$\mbf{I} \in \mathbb{I}\mathbb{R}$ --- \textit{область определения функции},\\
$\omega \in \mathbb{R}$ --- \textit{характеристика интервала}.

\end{definition}

\noindentАвтором работы предлагается следующий алгоритм этого класса:

\begin{small}
\noindent\textbf{Входные данные:}

\begin{singlespacing}
    \begin{enumerate}
    
        \item Функция для вычисления \textit{характеристики интервала} $\omega$.
        
        \item Исходный интервал для поиска минимума функции.
        
        \item Положительное число $\varepsilon$ -  пороговое значение остановки алгоритма.
        
    \end{enumerate}
\end{singlespacing}

\noindent\textbf{Алгоритм:}

    \begin{singlespacing}
        \begin{enumerate}
        	
        	\item Инициализировать рабочий список рабочим интервалом, где предполагается поиск глобального минимума. Перейти на шаг 4.
        
        	\item Если рабочий список пуст или условие остановки выполнено, перейти на шаг 11, иначе перейти на шаг 5.
        	
        	\item Извлечь из рабочего списка первую пару $A = (\mbf{a}, \, a)$. Перейти на шаг 6.
        
        	\item Найти интервальное расширение функции $f(\mbf{a})$. Обновить гарантированное значение минимума числом $\overline{f(\mbf{a})}$. Перейти на шаг 7.
        	
        	\item Проверить, больше ли $f(\mbf{a})$, чем значение гарантированного значения минимума. Если больше, перейти на шаг 4, иначе перейти на шаг 8.
        	
        	\item Применить к $\mbf{a}$ процедуру сжатия функционального интервала к минимуму. Перейти на шаг 9.
        	
        	\item Если после процедуры сжатия $\mbf{a} = \varnothing$, перейти на шаг 4, иначе перейти на шаг 10.
        	
        	\item Произвести дихотомию $\mbf{a} \rightarrow \big\{ \mbf{a}_{-} = [ \, \underline{\mbf{a}}, \, \text{mid} \, \mbf{a} \, ], \, \mbf{a}_{+} = [ \, \text{mid} \, \mbf{a}, \, \overline{\mbf{a}} \, ] \big\}$. Вставить рабочие интервалы $\big(\mbf{a}_{-}, \, \omega_{-} \big)$ и $\big(\mbf{a}_{+}, \, \omega_{+} \big)$ в рабочий список. Перейти на шаг 4.
        	
        	\item Вывести в качестве ответа список всех полученных интервалов. Алгоритм завершён.
        	
        \end{enumerate}
    \end{singlespacing}
\end{small}

\section{Численные эксперименты}

Рассмотрим применение построенной функциональной интервальной арифметики и продемонстрируем её эффективность при решении численных задач.

\subsection{Задача решения уравнения}

Рассмотрим задачу доказательного численного нахождения корней уравнения на вещественной оси $\mathbb{R}$. Для примера рассмотрим функцию
\begin{equation*}
    f(x) = x^{3} - x^{2} - x = 0.
\end{equation*}

Используя следствие из теоремы Руше, возьмём в качестве начального интервала поиска вещественных корней интервал
\begin{equation*}
    \begin{array}{l}
        \big[ \, -1 - \max\limits_{i = 1, \dots, n} \, |a_{i}|, \, 1 + \max\limits_{i = 1, \dots, n} \, |a_{i}|\, \big] = \\
        \hspace{1cm} = \big[ \, -1 - \max\{ 2 \, / \, 1, \, 2 \, / \, 1, \, 2 \, / \, 1 \}, \, 1 + \max\{ 2 \, / \, 1, \, 2 \, / 1, \, 2 \, / \, 1 \} \, \big] = \\
        \hspace{2cm} = [ \, -2, \, 2 \, ].
    \end{array}
\end{equation*}

Реализуем алгоритм поиска корней, рассмотренный в пункте 4.4.2 с интервалами семейств $\mathbb{L}\mathbb{F}\mathbb{R}(x)$ и $\mathbb{I}\mathbb{R}$. Для случая классических интервалов также рассмотрим центрированные формы. 

В качестве функции вычисления характеристики рабочего интервала была выбрана
\begin{equation*}
    -\text{mag} \, f(\mbf{I}), \qquad \mbf{I} \in \mathbb{F}\mathbb{R}.
\end{equation*}

В качестве порогового значения остановки алгоритма был взят $\varepsilon = 10 ^ {-2}$.

Программа была написана на языке программирования \texttt{C Sharp} и запускалась на ЭВМ с 8-ми ядерным процессором \texttt{Intel Core i7-3770}.

\clearpage

\subsubsection{Иллюстрации итераций с интервалами $\mathbb{I}\mathbb{R}$}

Также приведём иллюстрации для некоторых итераций, которые проводятся с помощью интервалов семейства $\mathbb{I}\mathbb{R}$.

\begin{figure}[ht]
	\begin{center}
		
		\includegraphics[width = 0.45 \linewidth]{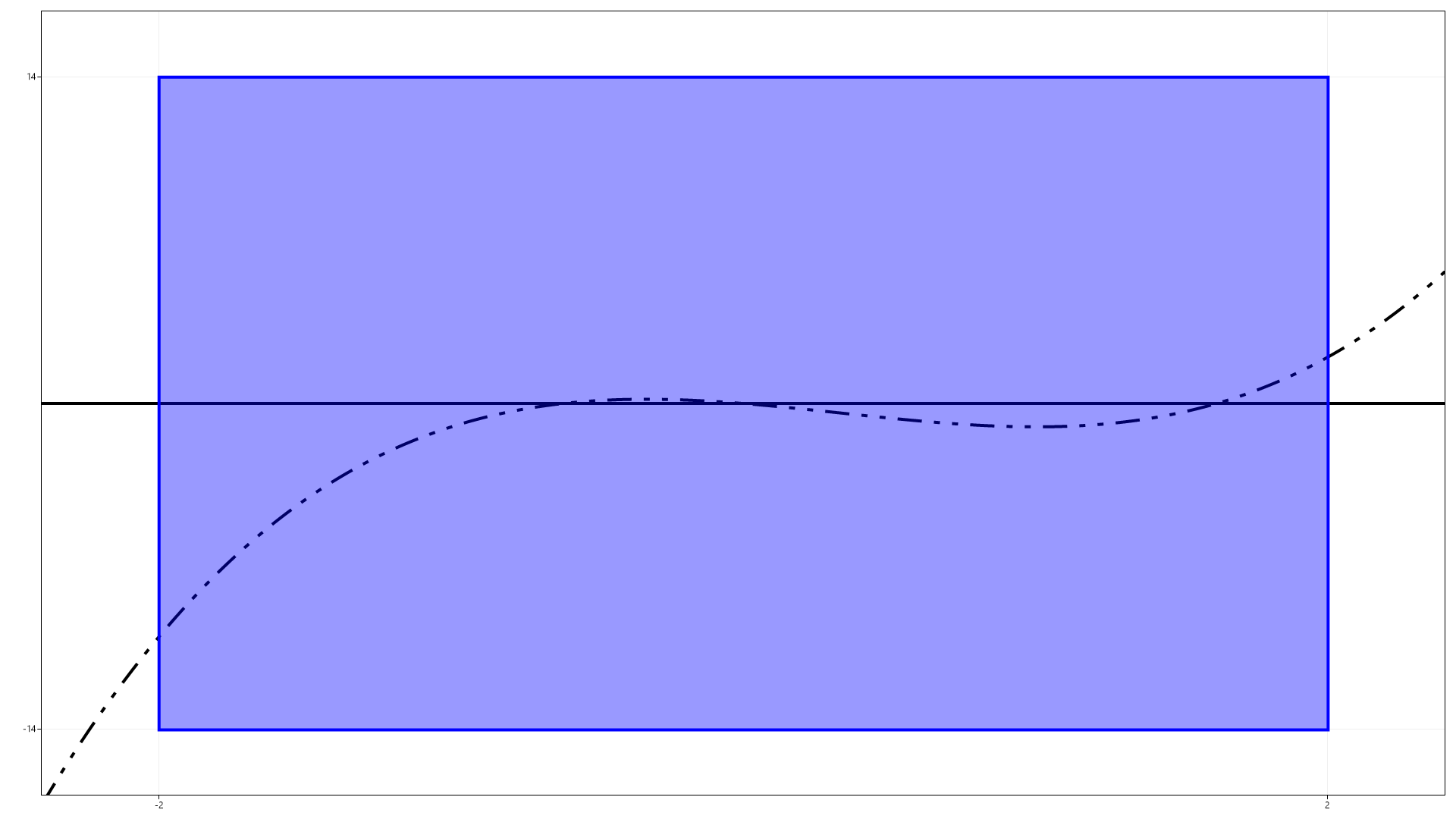}
		\includegraphics[width = 0.45 \linewidth]{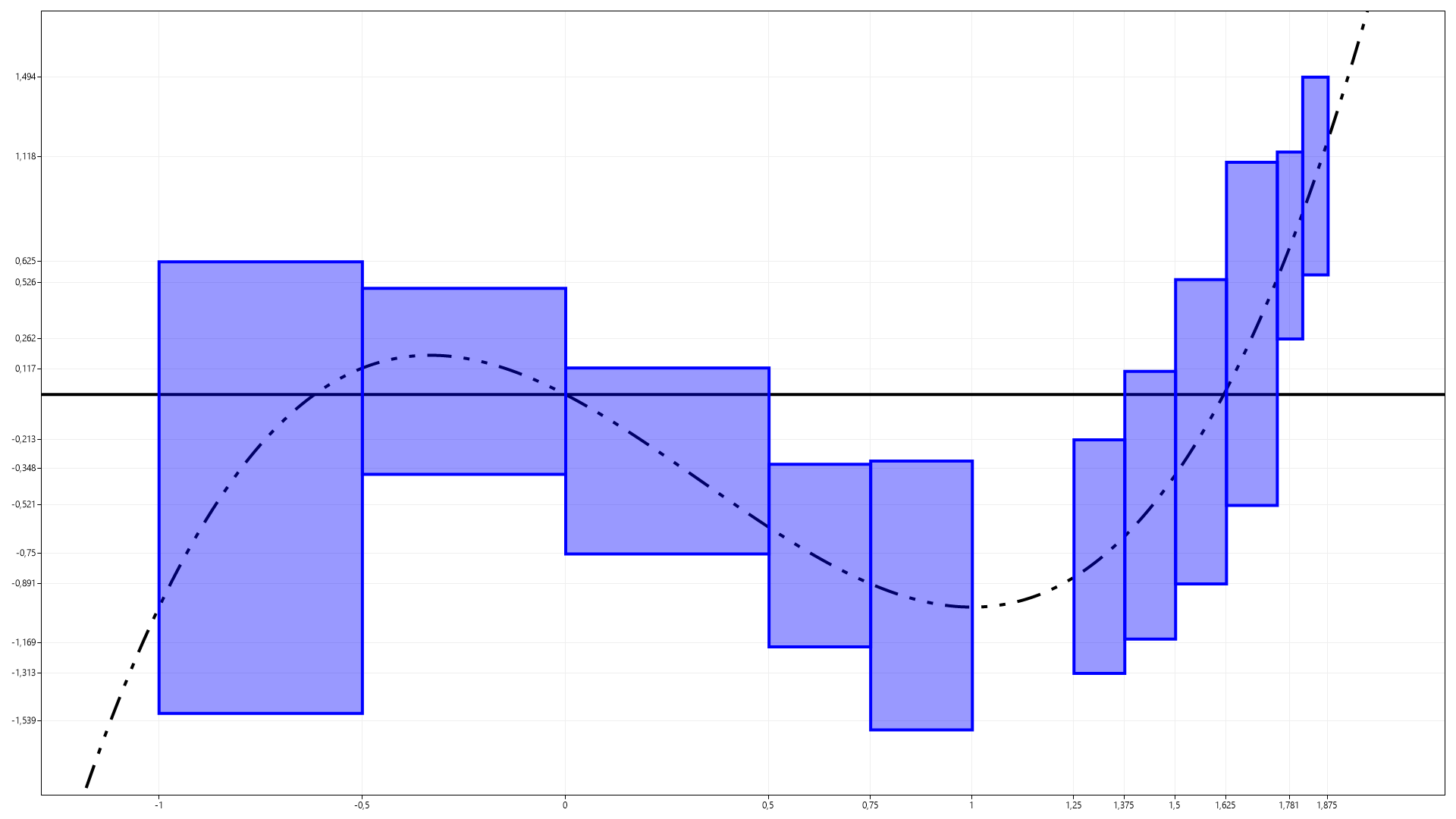}
		\caption{Итерации 0 и 14.}
		\label{fig:equation_classic_iteration_0_14}
	
	\end{center}
\end{figure}

\begin{figure}[ht]
	\begin{center}
		
		\includegraphics[width = 0.45 \linewidth]{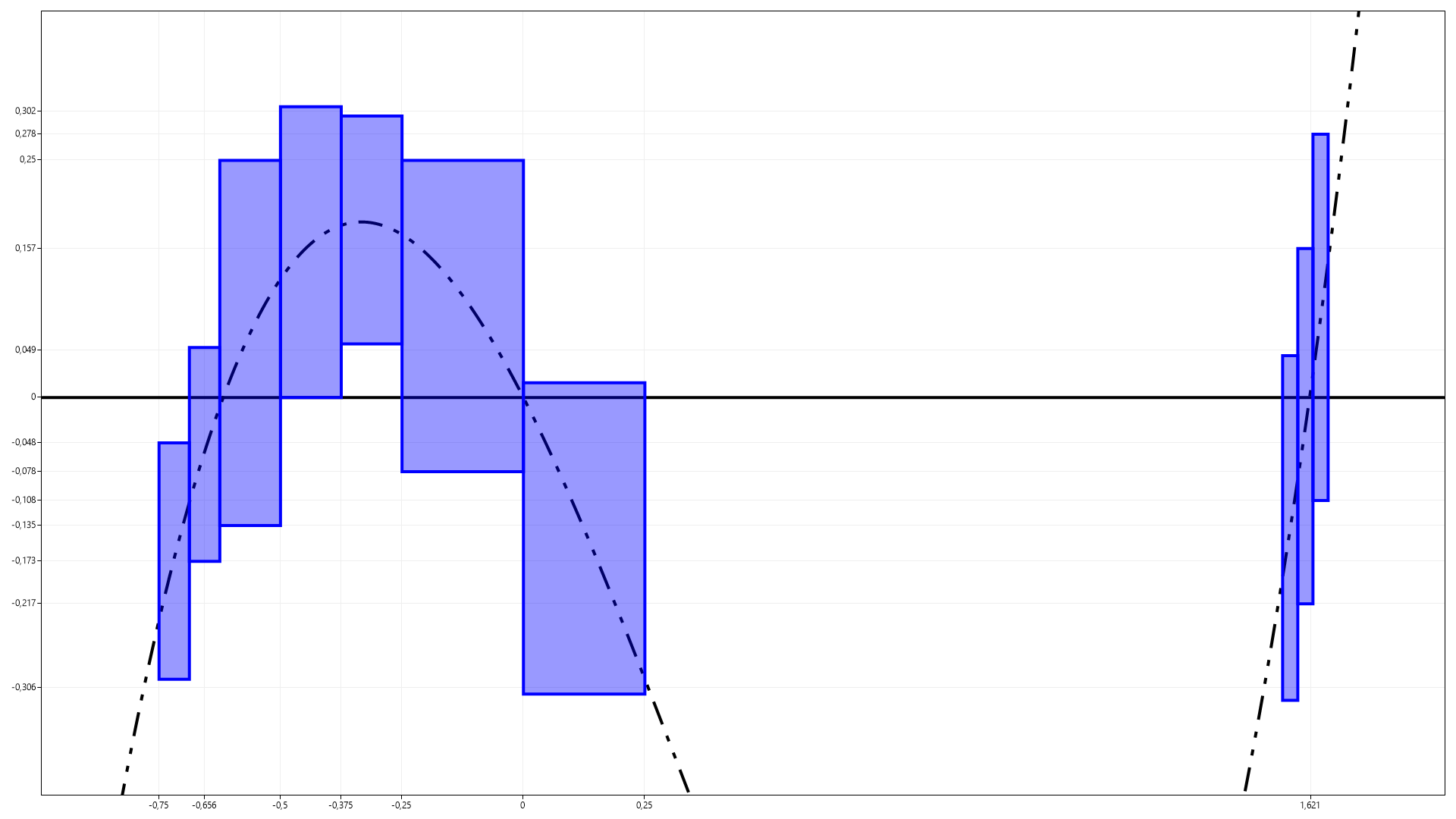}
		\includegraphics[width = 0.45 \linewidth]{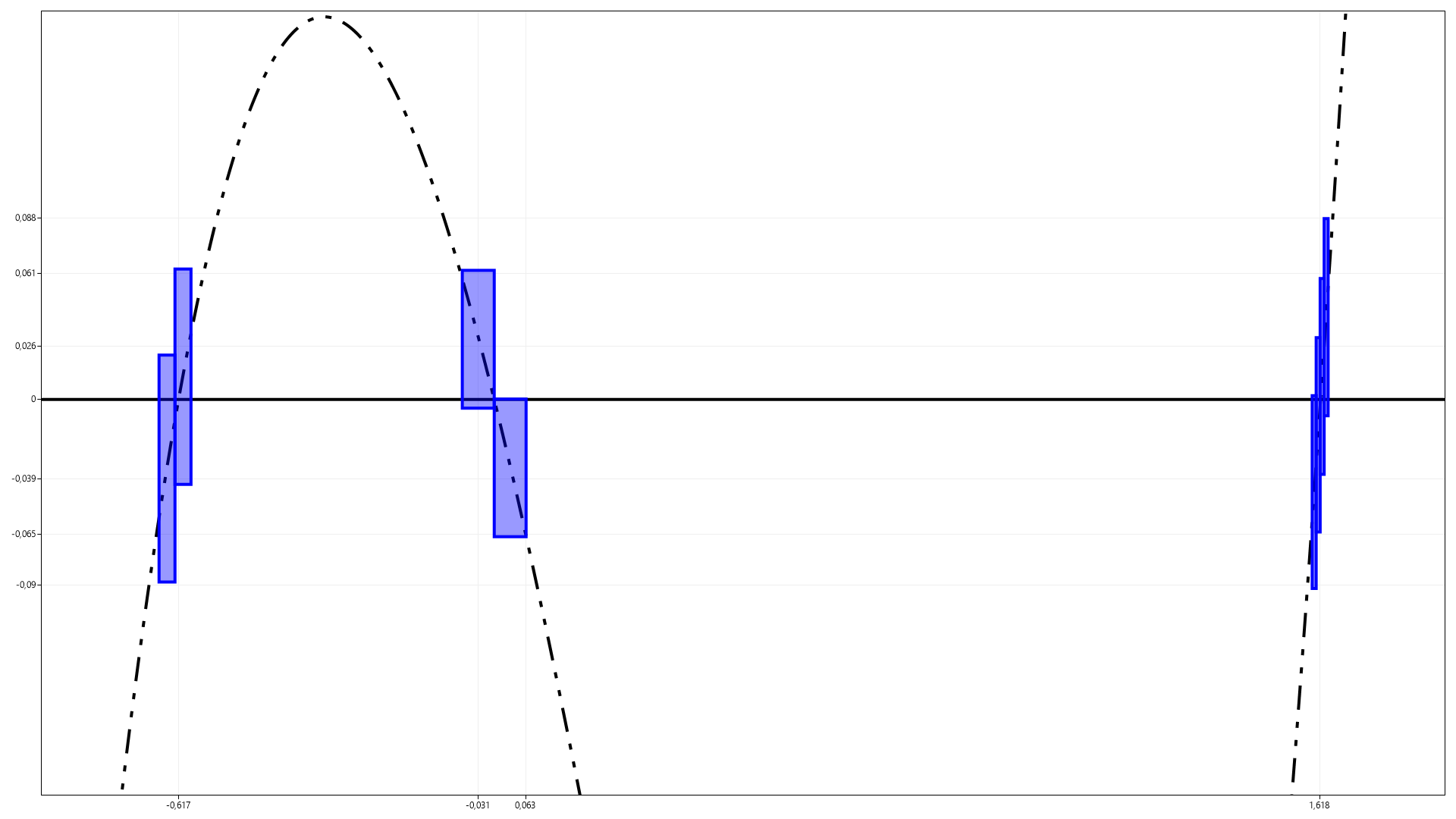}
		\caption{Итерации 28 и 42.}
		\label{fig:equation_classic_iteration_28_42}
	
	\end{center}
\end{figure}

\begin{figure}[ht]
	\begin{center}
		
		\includegraphics[width = 0.45 \linewidth]{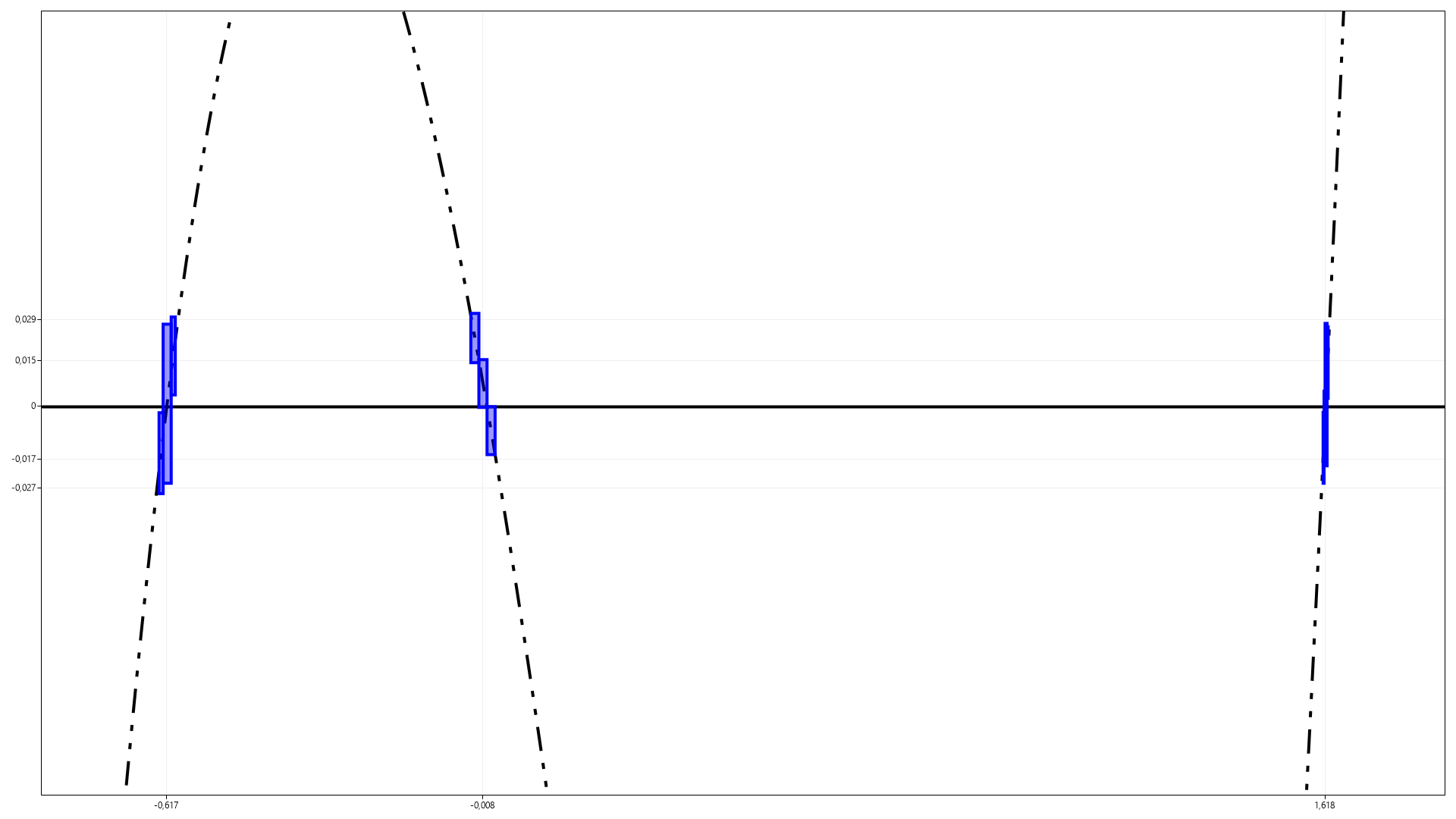}
		\includegraphics[width = 0.45 \linewidth]{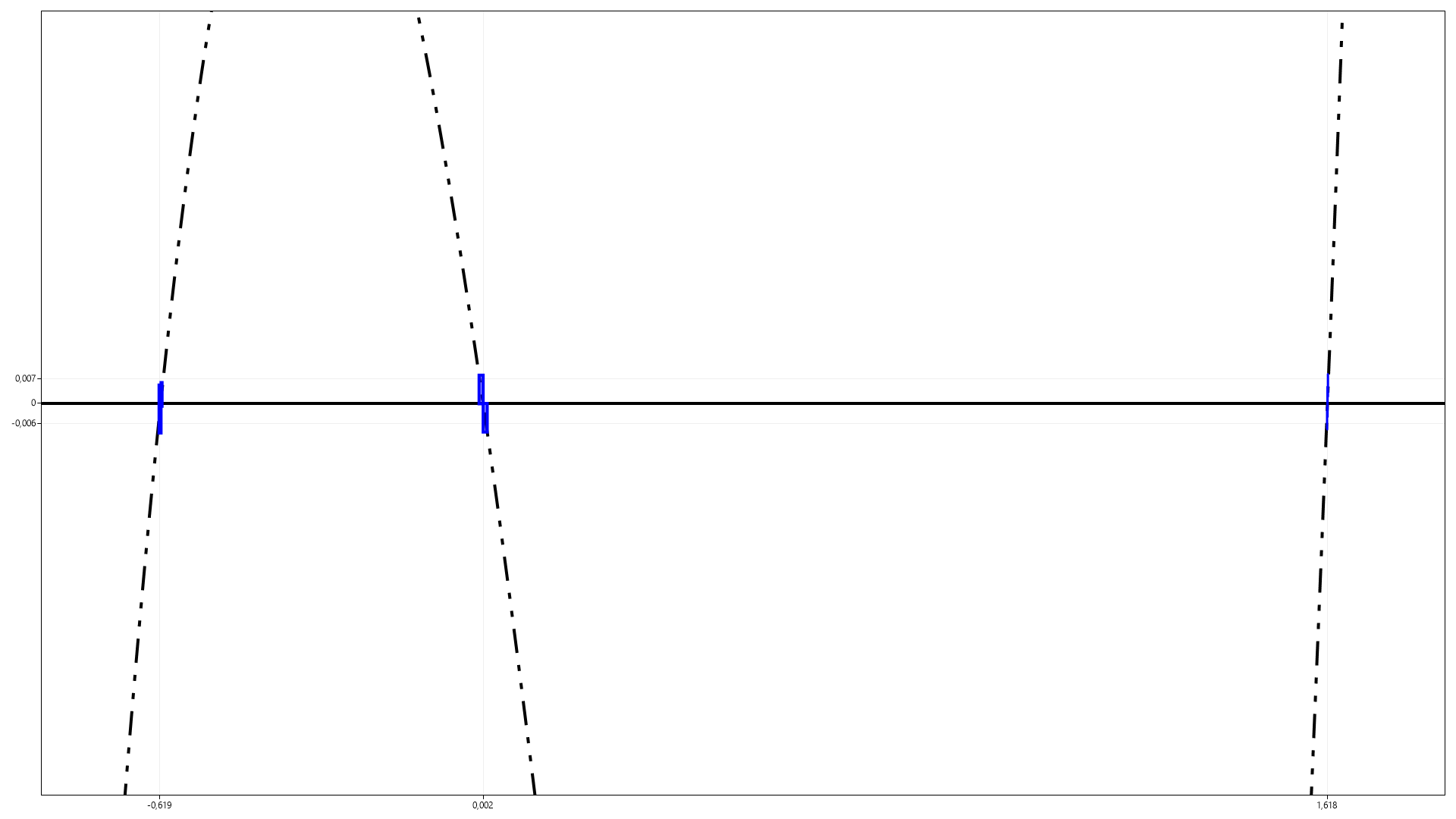}
		\caption{Итерации 56 и 69 --- последняя итерация алгоритма.}
		\label{fig:equation_classic_iteration_56_69}
	
	\end{center}
\end{figure}

\clearpage

\subsubsection[Иллюстрации итераций с центрированными формами]%
              {Иллюстрации итераций с центрированными формами интервалов $\mathbb{I}\mathbb{R}$}

Также приведём иллюстрации для некоторых итераций, которые проводятся с помощью интервалов семейства $\mathbb{I}\mathbb{R}$ с помощью центрированных форм.

\begin{figure}[ht]
	\begin{center}
		
		\includegraphics[width = 0.45 \linewidth]{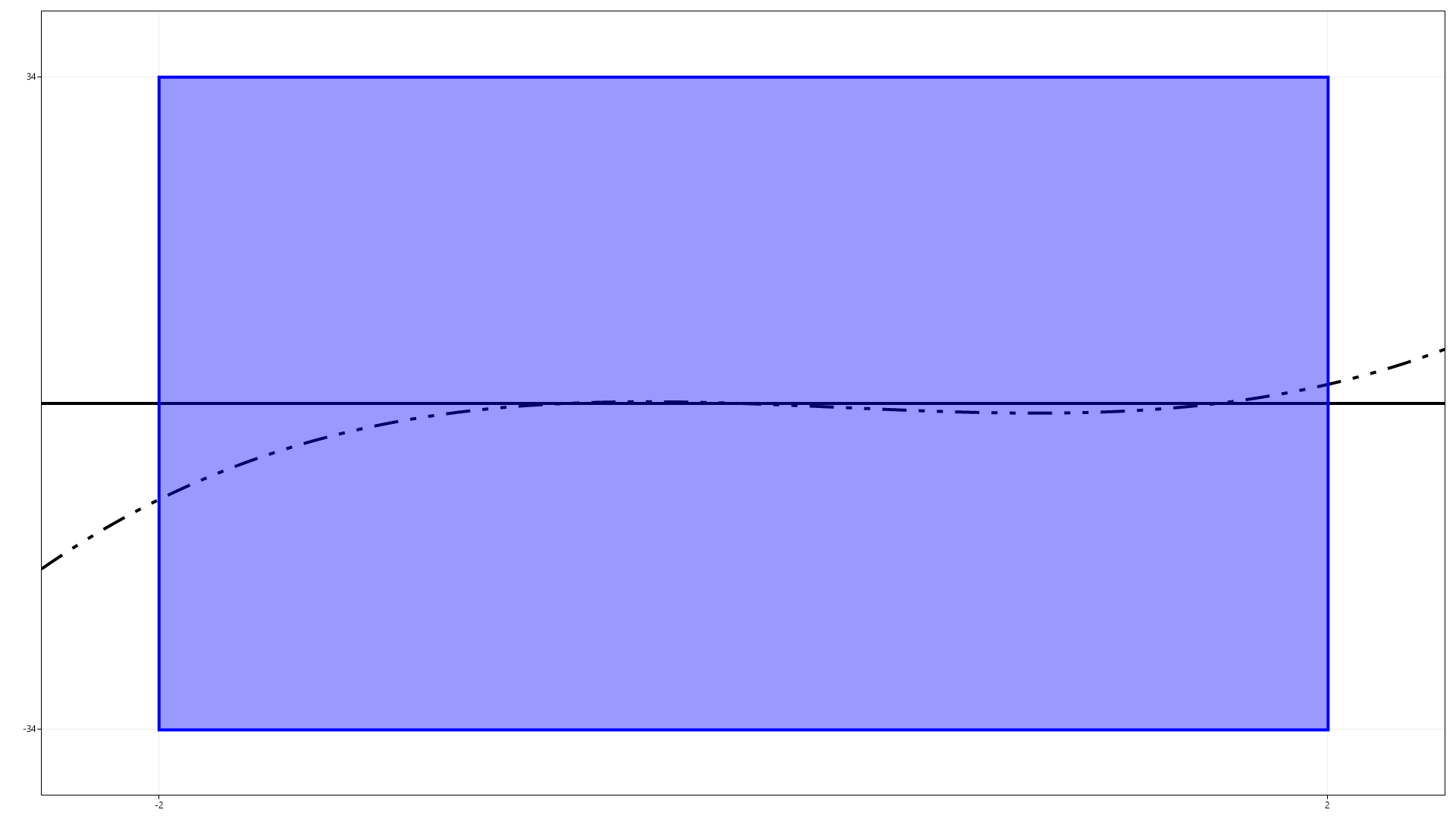}
		\includegraphics[width = 0.45 \linewidth]{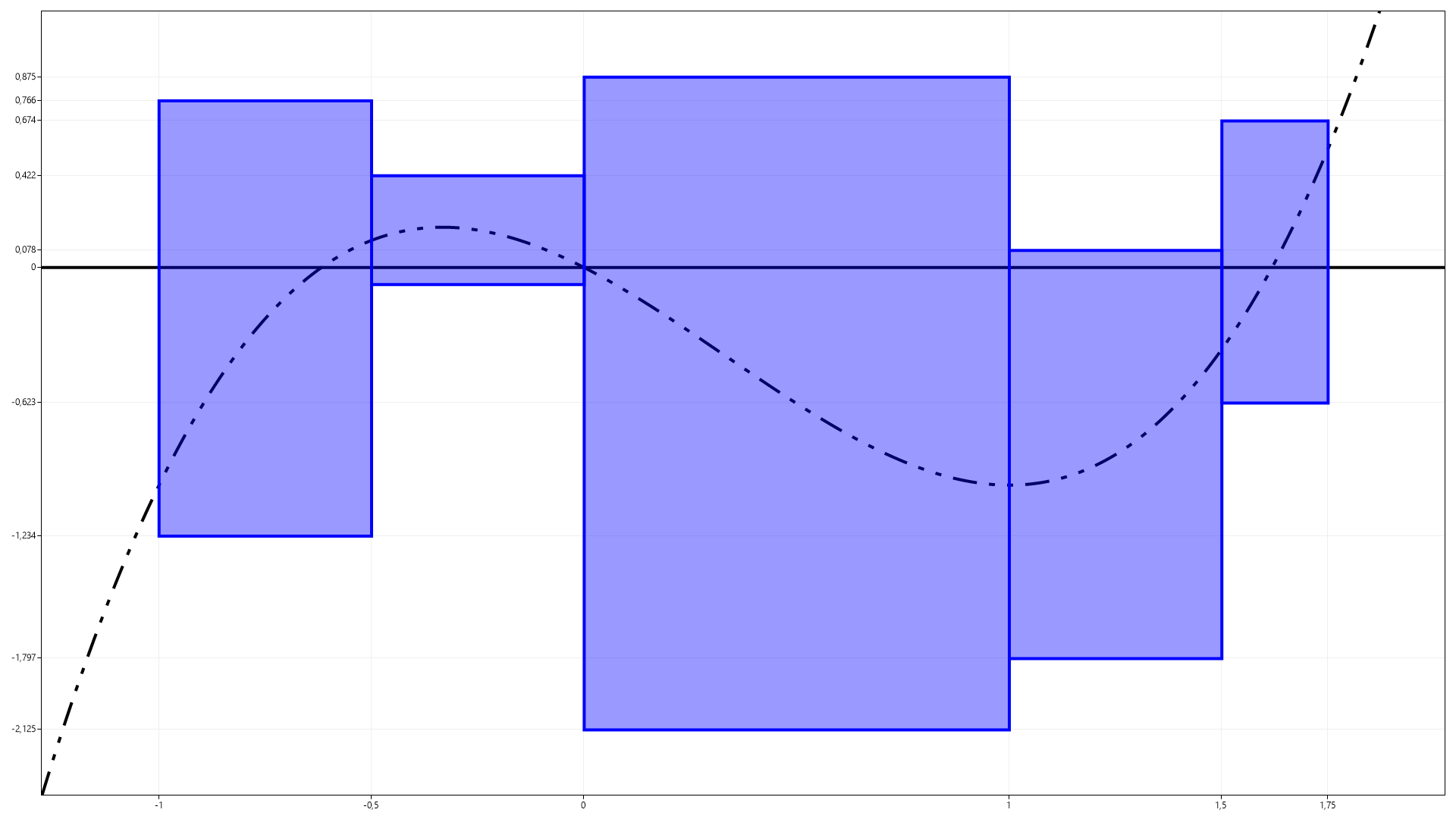}
		\caption{Итерации 0 и 7.}
		\label{fig:equation_center_iteration_0_7}
	
	\end{center}
\end{figure}

\begin{figure}[ht]
	\begin{center}
		
		\includegraphics[width = 0.45 \linewidth]{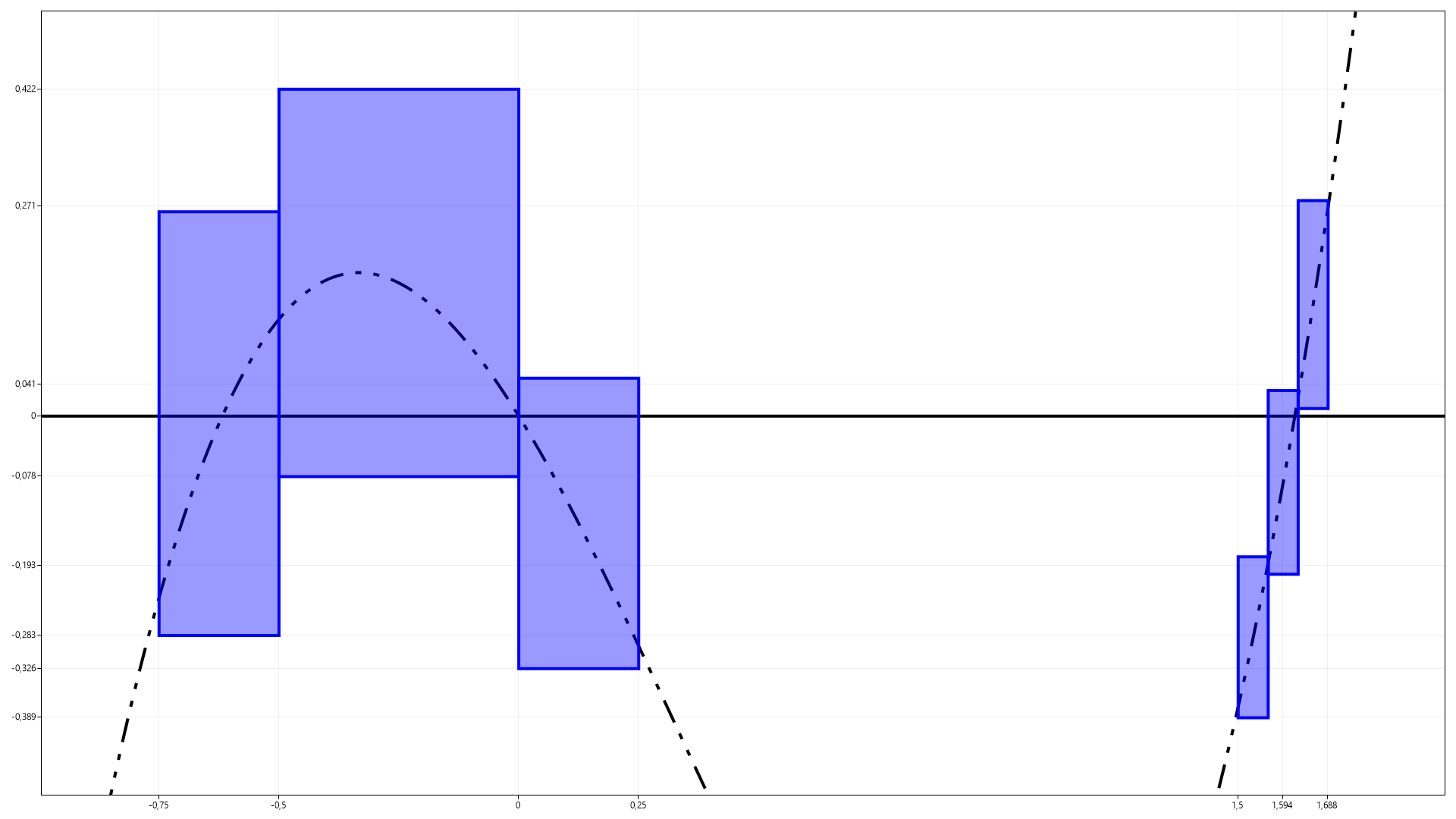}
		\includegraphics[width = 0.45 \linewidth]{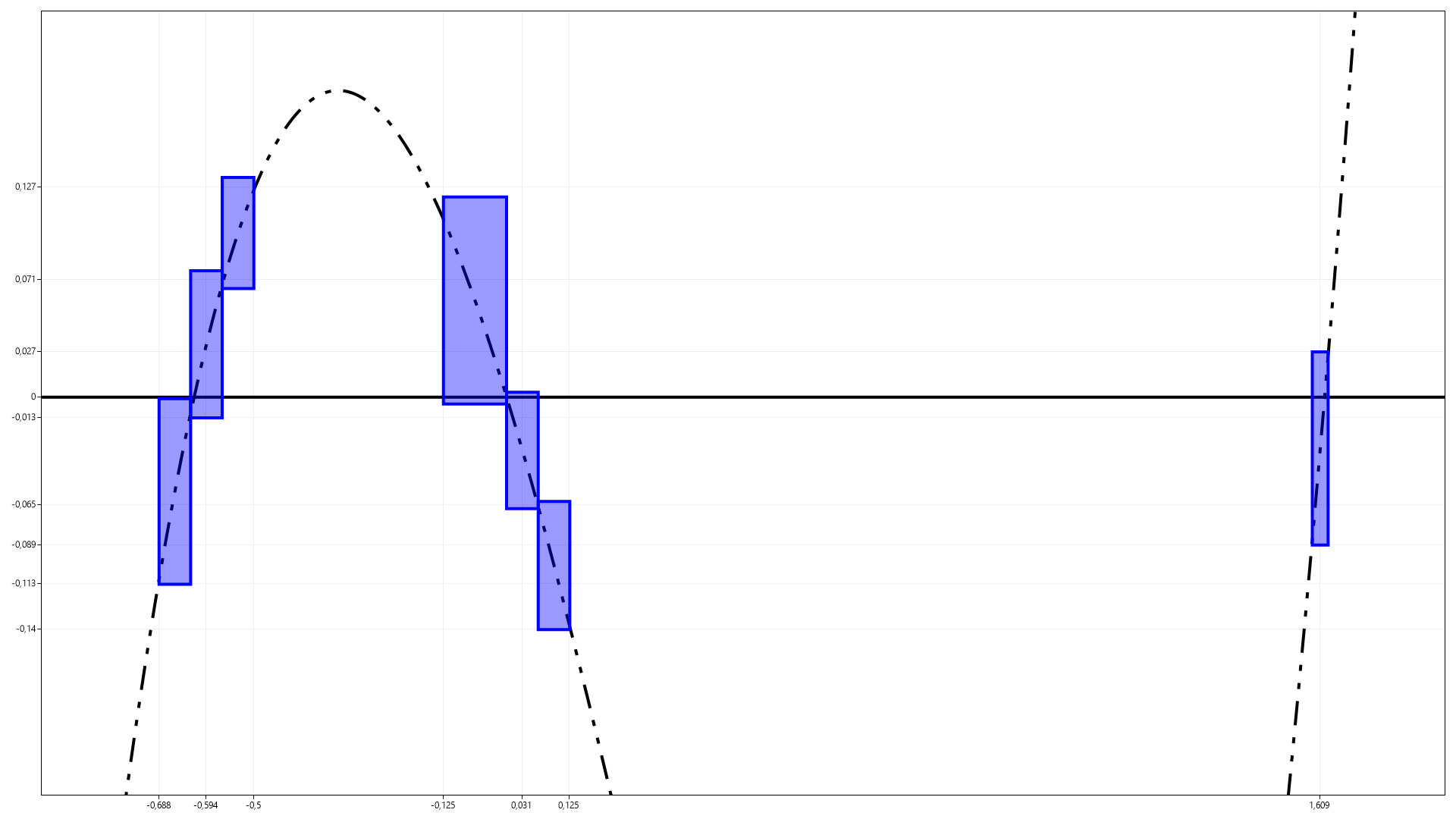}
		\caption{Итерации 14 и 22.}
		\label{fig:equation_center_iteration_14_22}
	
	\end{center}
\end{figure}

\begin{figure}[ht]
	\begin{center}
		
		\includegraphics[width = 0.45 \linewidth]{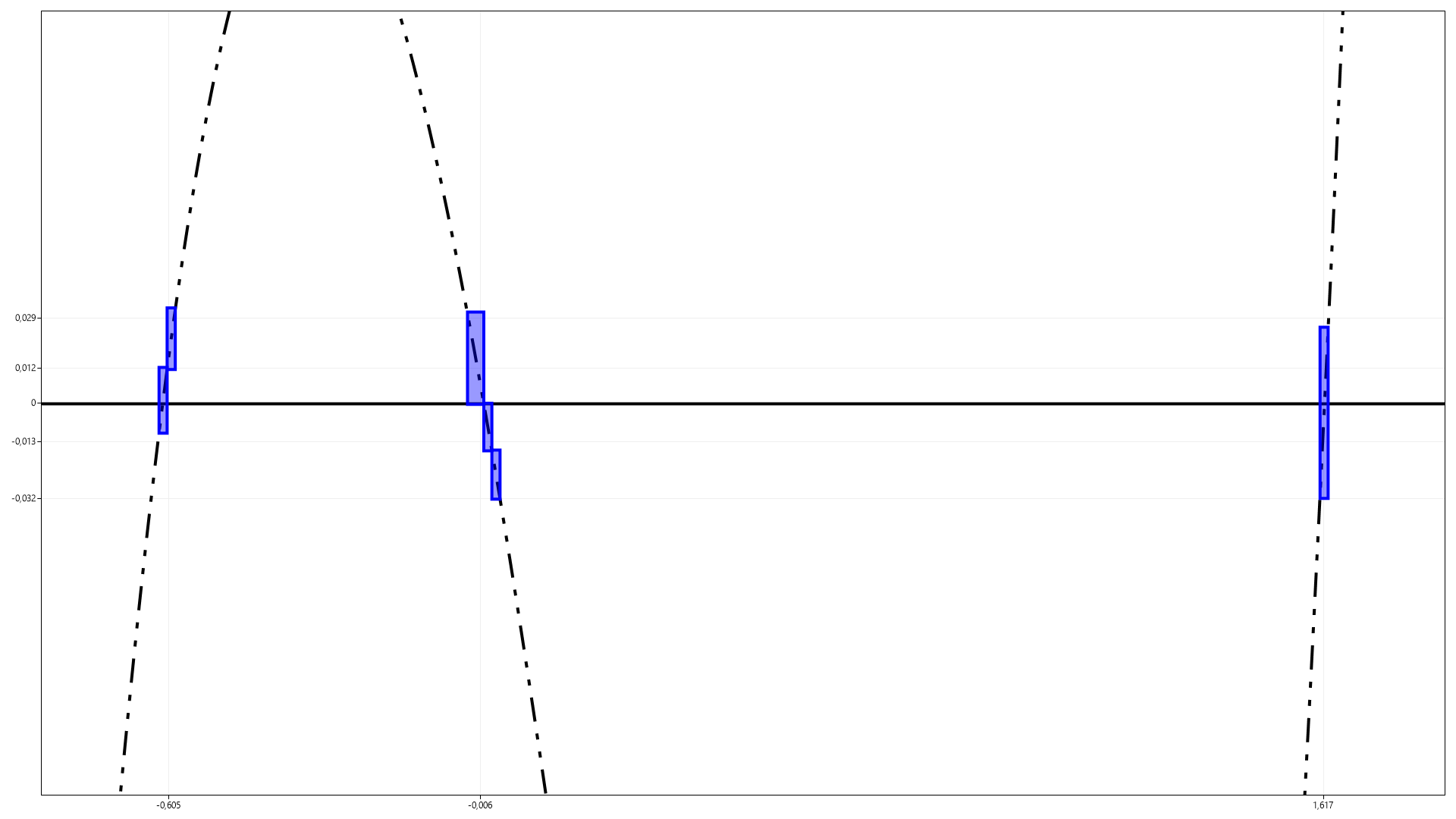}
		\includegraphics[width = 0.45 \linewidth]{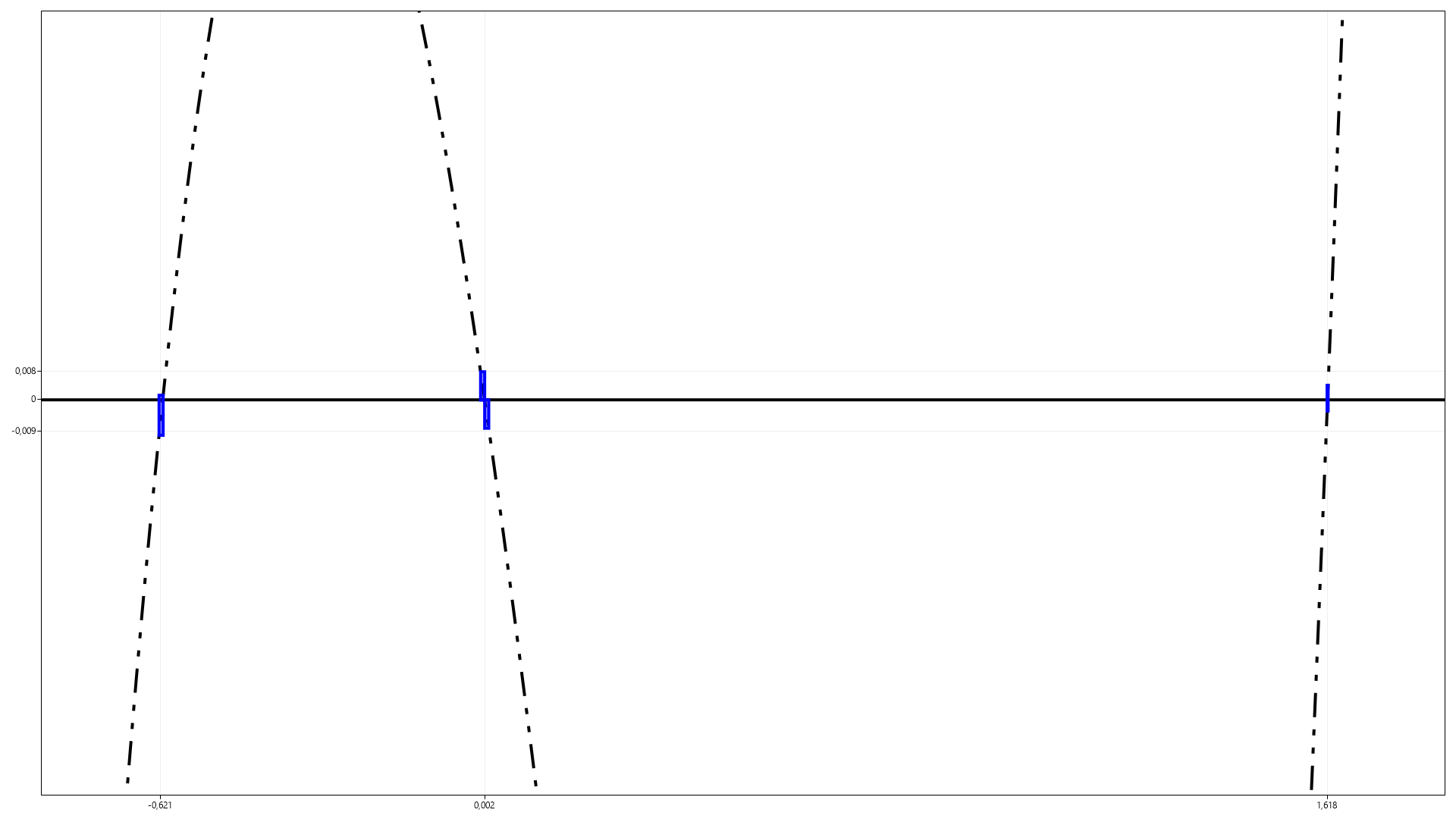}
		\caption{Итерации 29 и 37 --- последняя итерация алгоритма.}
		\label{fig:equation_center_iteration_29_37}
	
	\end{center}
\end{figure}

\clearpage

\subsubsection{Иллюстрации итераций с интервалами $\mathbb{L}\mathbb{F}\mathbb{R}$}

Далее приведены иллюстрации для некоторых итераций при использовании интервалов семейства $\mathbb{L}\mathbb{F}\mathbb{R}(x)$.

\begin{figure}[ht]
	\begin{center}
		
		\includegraphics[width = 0.45 \linewidth]{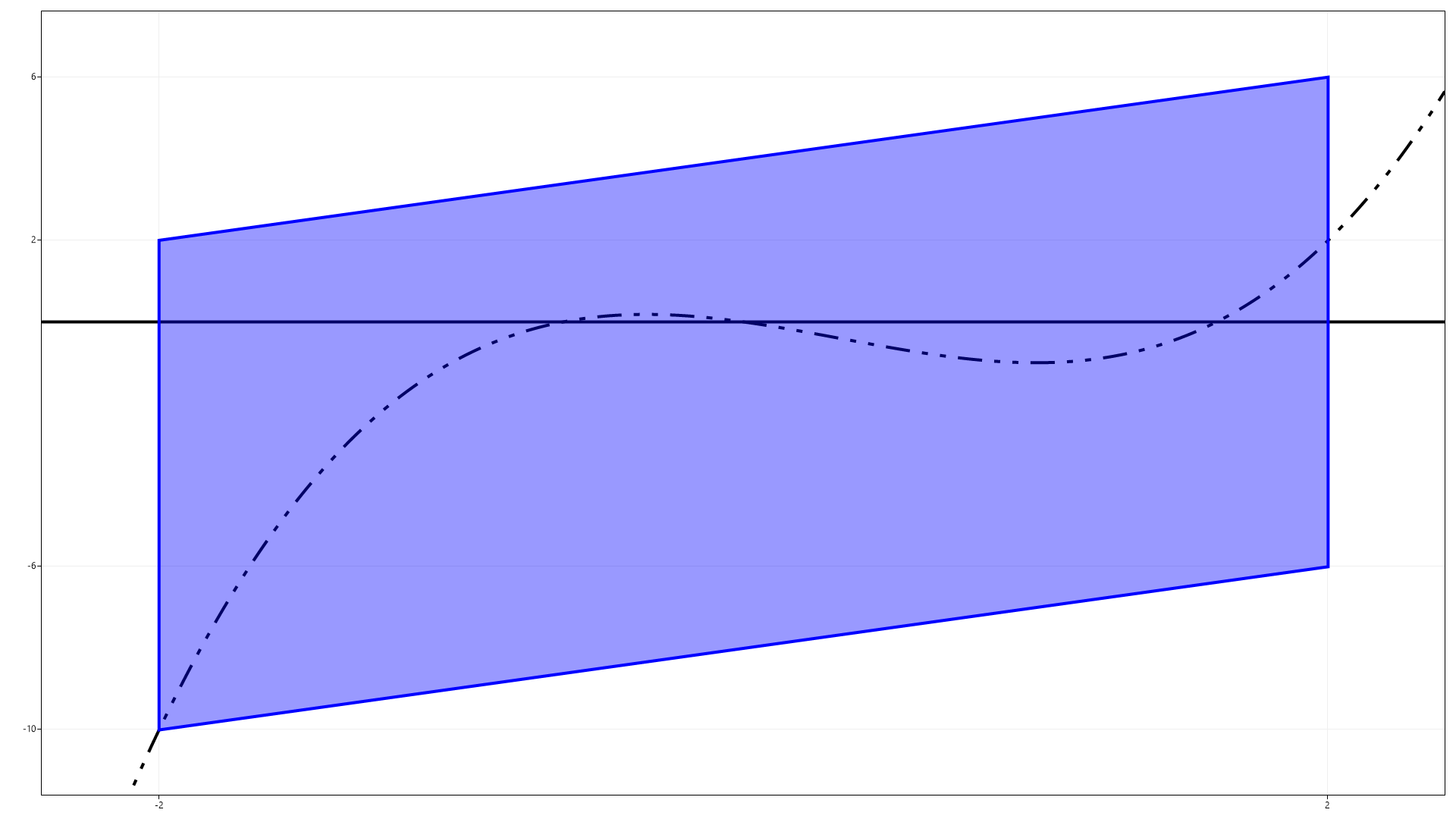}
		\includegraphics[width = 0.45 \linewidth]{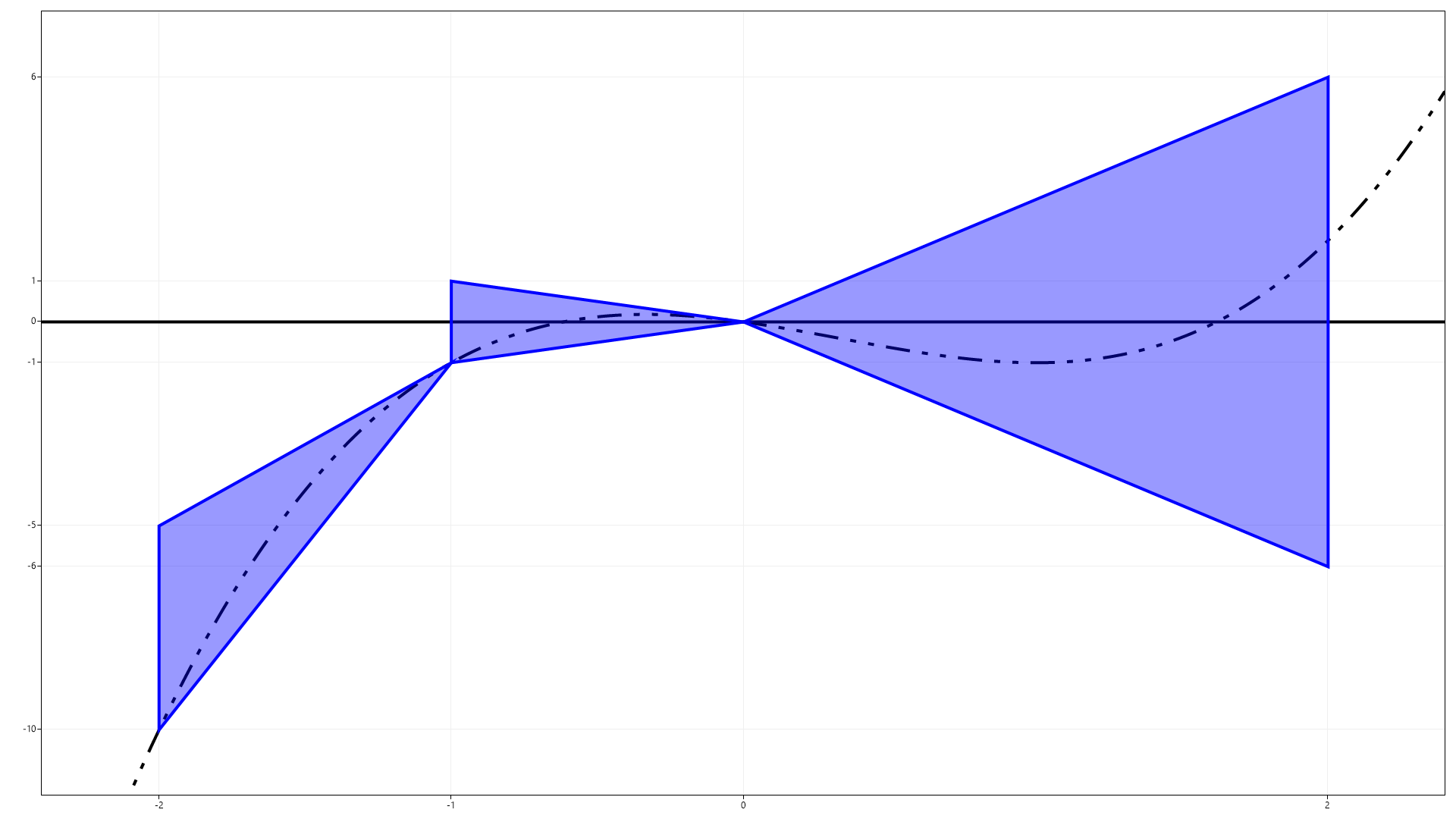}
		\caption{Итерации 0 и 2.}
		\label{fig:equation_iteration_0_2}
		
	\end{center}
\end{figure}

\begin{figure}[ht]
	\begin{center}
	
		\includegraphics[width = 0.45 \linewidth]{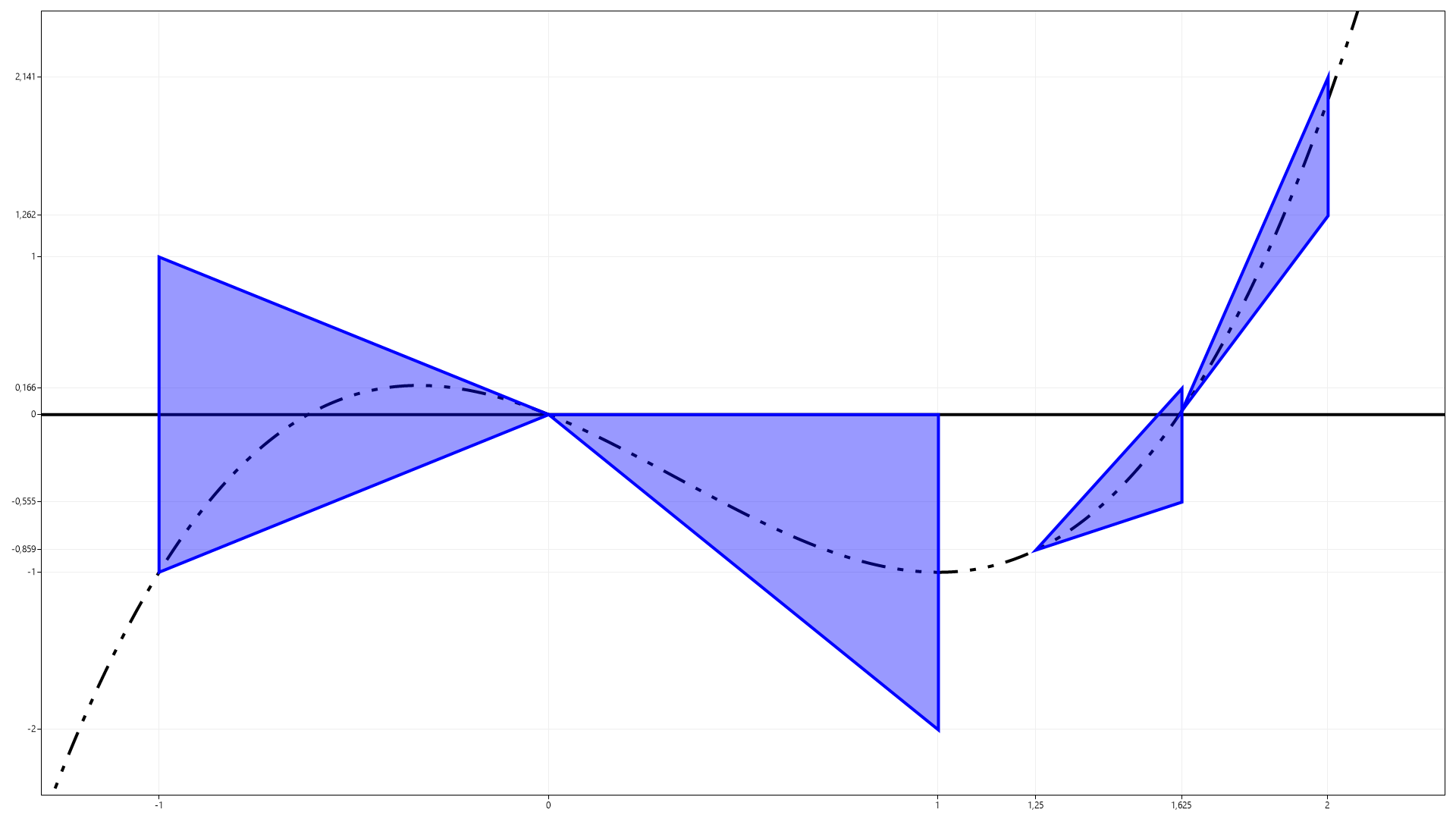}
		\includegraphics[width = 0.45 \linewidth]{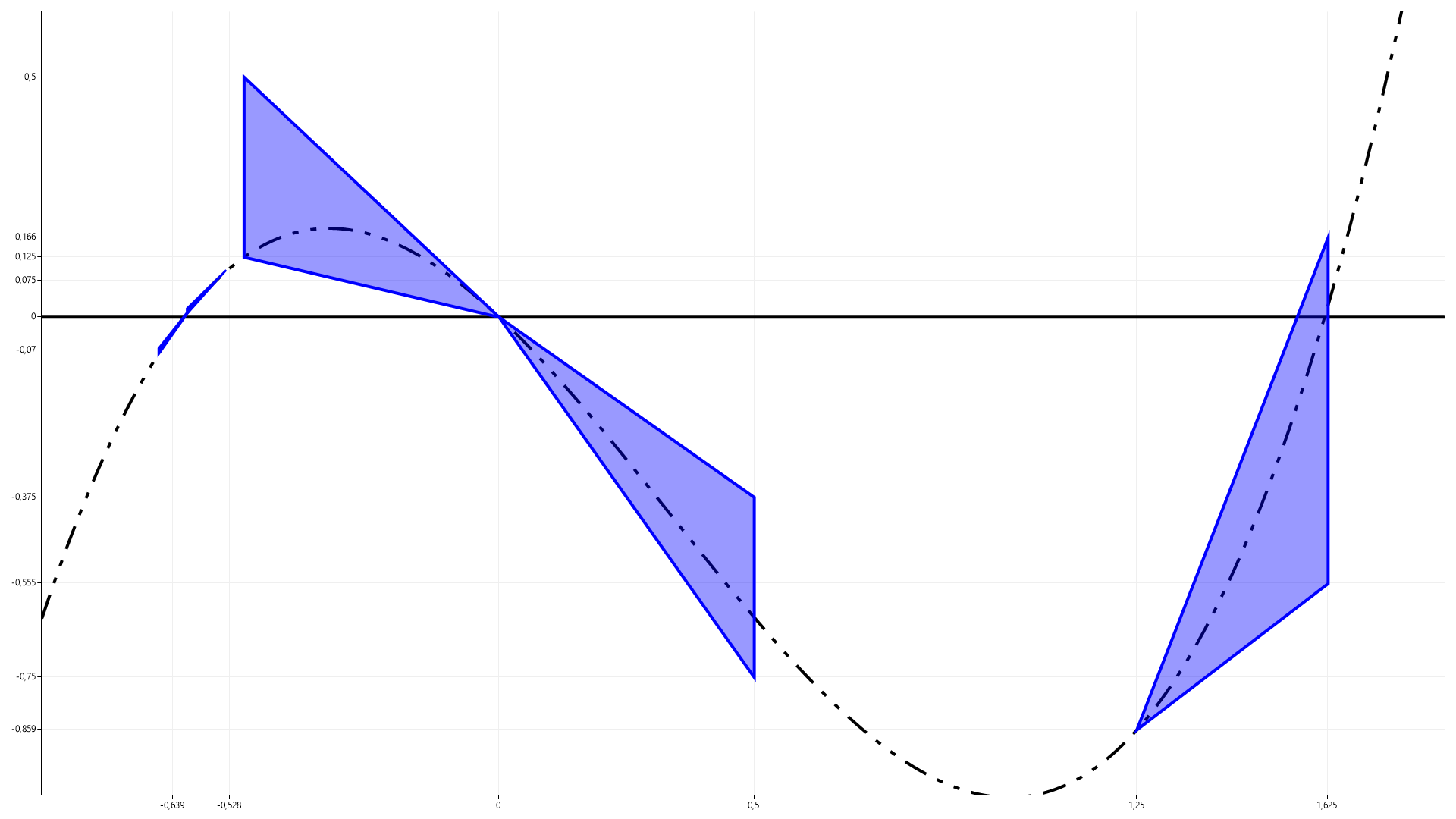}
		\caption{Итерации 4 и 7.}
		\label{fig:equation_iteration_4_7}
	
	\end{center}
\end{figure}

\begin{figure}[ht]
	\begin{center}
	
		\includegraphics[width = 0.45 \linewidth]{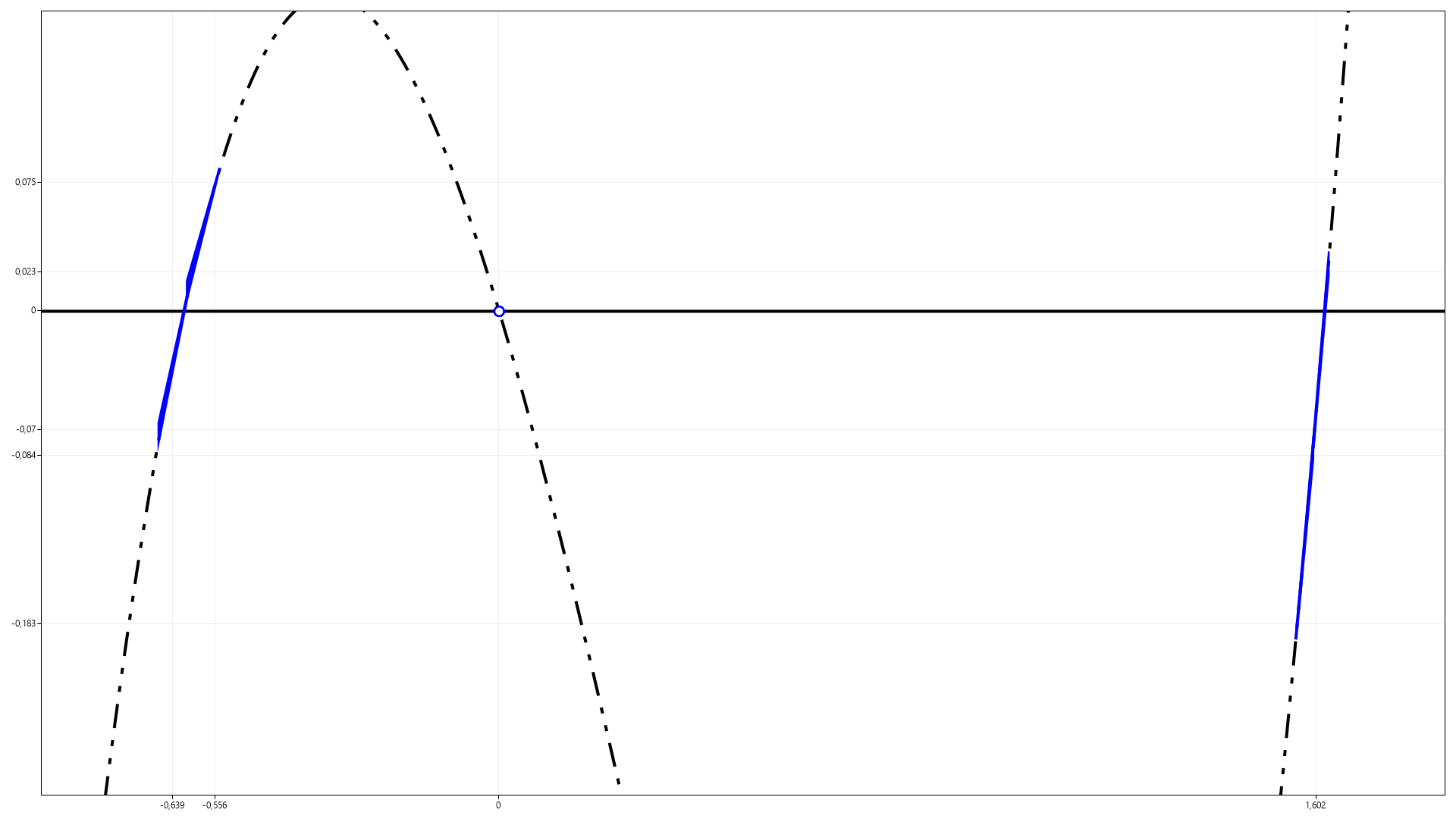}
		\includegraphics[width = 0.45 \linewidth]{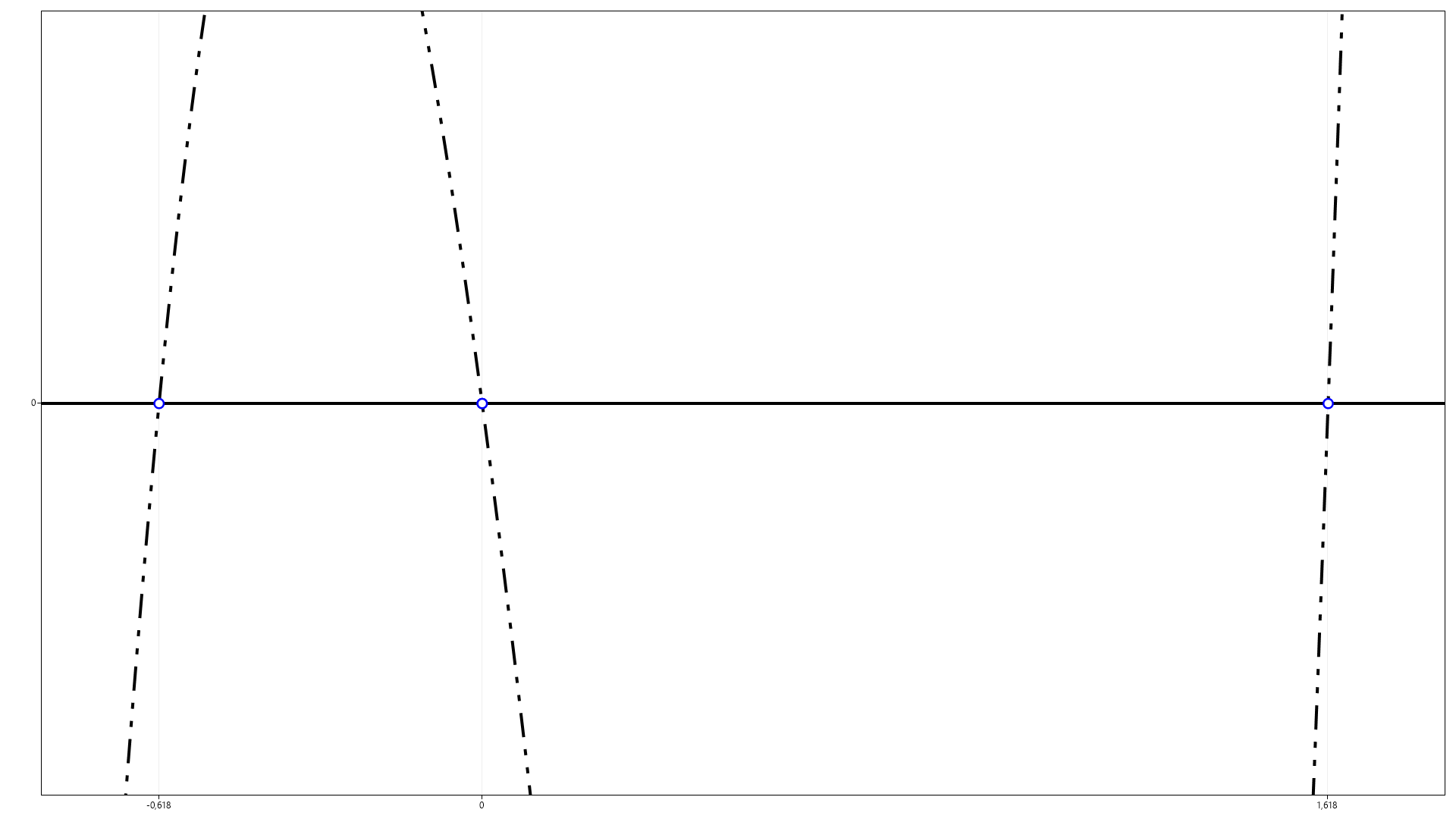}
		\caption{Итерации 10 и 13 --- последняя итерация алгоритма.}
		\label{fig:equation_iteration_10_13}
		
	\end{center}
\end{figure}	

\clearpage

\subsubsection{Сравнение результатов}

Далее приведён график по по результатам численного эксперимента.

\begin{figure}[ht]
	\begin{center}
		
		\includegraphics[width = 0.75 \linewidth]{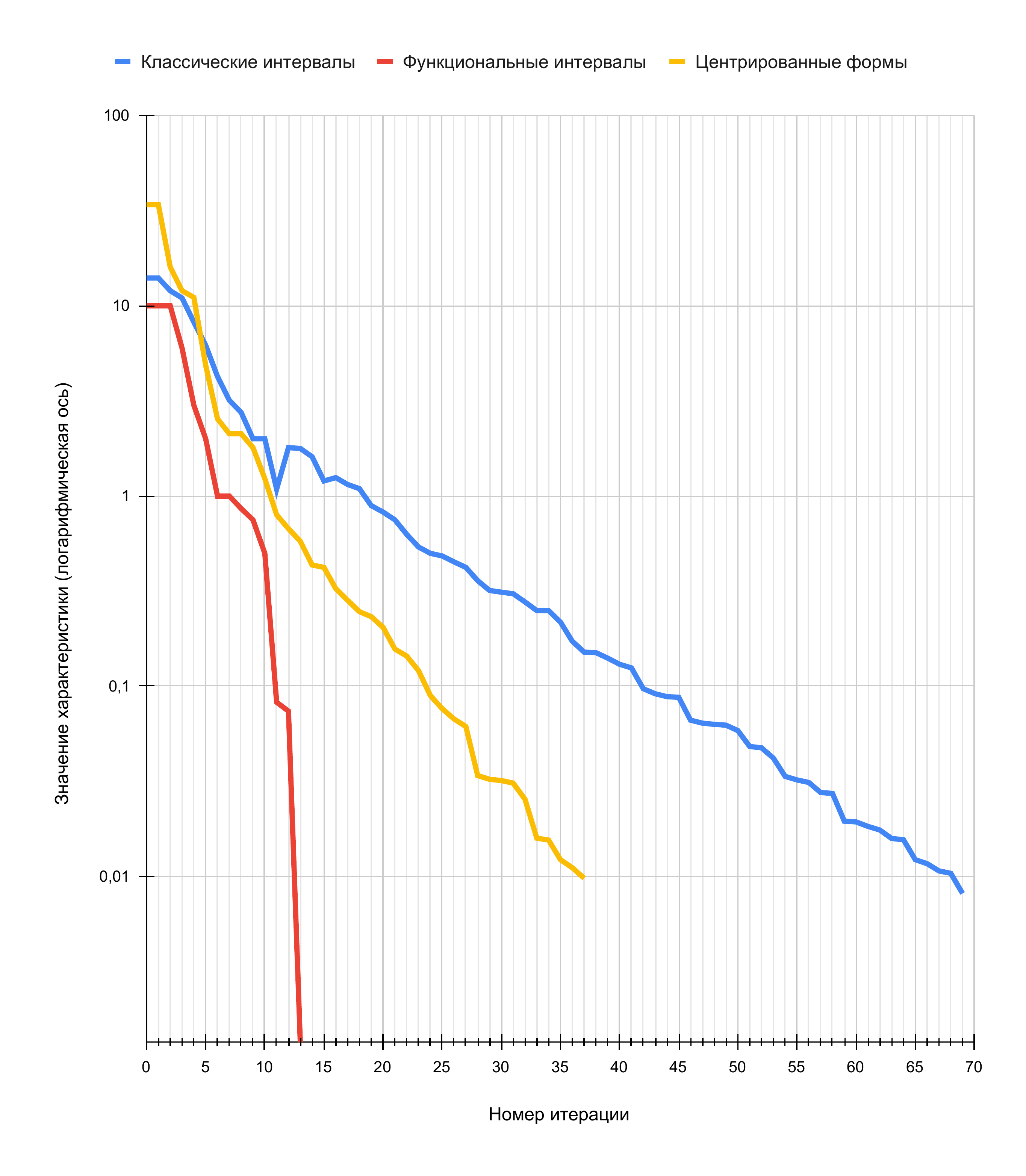}
		\caption{График сравнения максимальной характеристики для разных видов интервалов.}
		
	\end{center}
\end{figure}

\clearpage

\subsubsection{Дополнительные результаты}

Для проведения дополнительных исследований алгоритм был запущен с функцией вычисления характеристики рабочего интервала
\begin{equation*}
    -\text{wid} \, f(\mbf{I}), \qquad \mbf{I} \in \mathbb{F}\mathbb{R}.
\end{equation*}

В качестве \textit{порогового значения остановки алгоритма} было выбрано $\varepsilon = 10 ^ {-12}$. Критерий остановки был дополнен условием на максимальное количество совершённых итераций, равное $10 ^ {4}$.

\begin{figure}[ht]
	\begin{center}
		
		\includegraphics[width = 0.75 \linewidth]{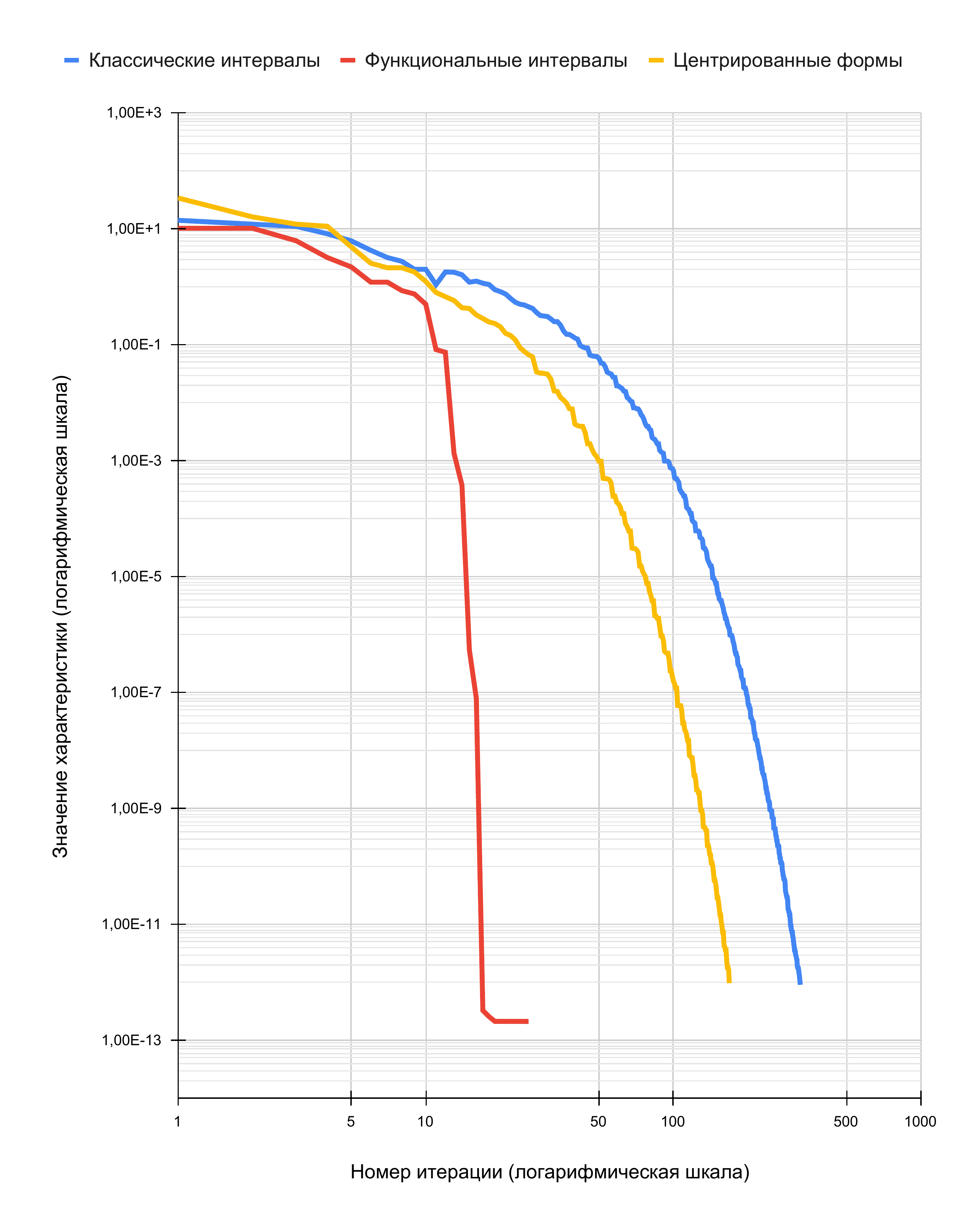}
		\caption{График сравнения максимальной характеристики рабочего интервала для разных видов интервалов на каждой итерации алгоритма.}
		
	\end{center}
\end{figure}

\begin{figure}[ht]
	\begin{center}
		
		\includegraphics[width = 0.75 \linewidth]{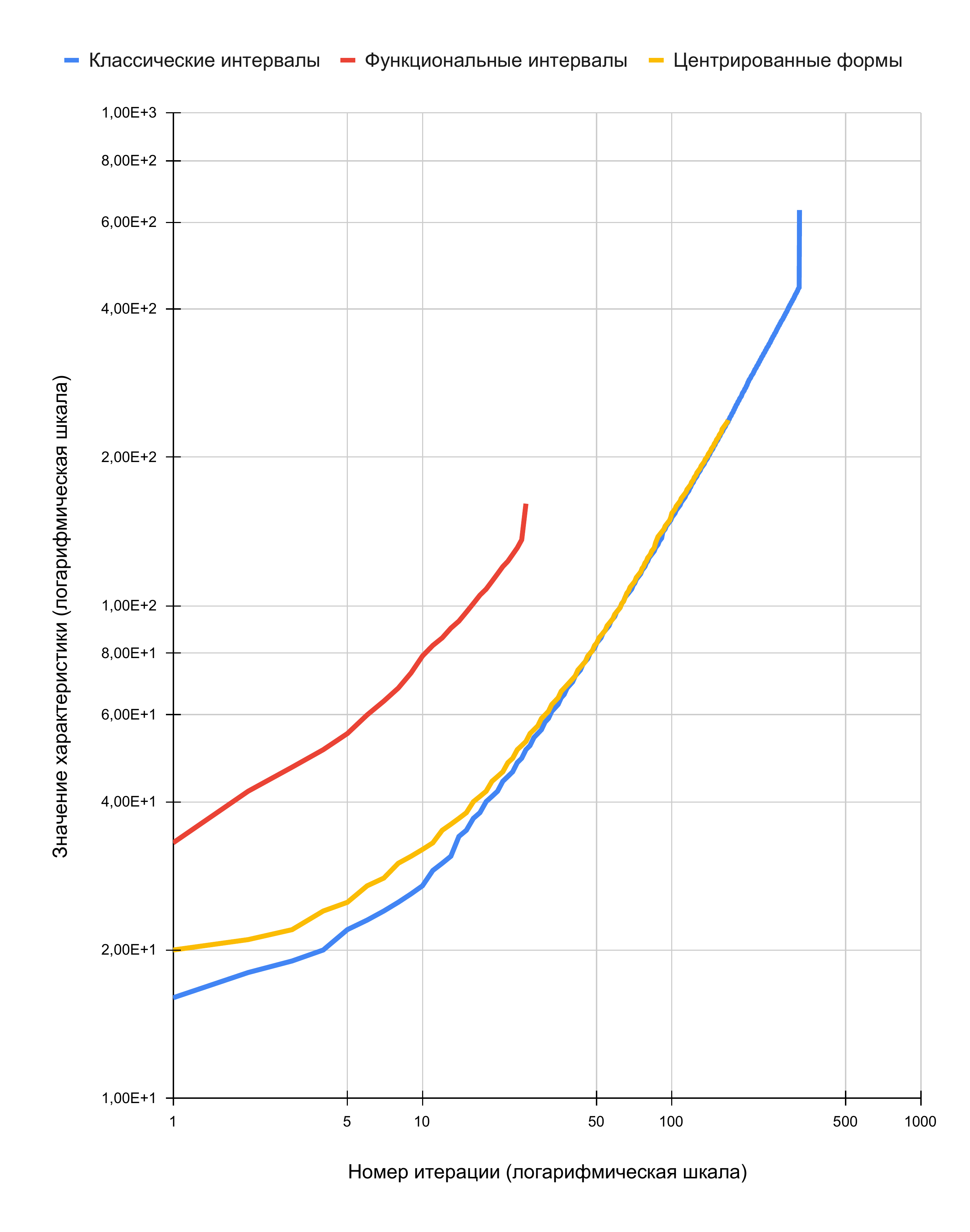}
		\caption{Суммарное время работы в миллисекундах, достигнутое на каждой итерации алгоритма.}
		
	\end{center}
\end{figure}

\begin{figure}[ht]
	\begin{center}
		
		\includegraphics[width = 0.75 \linewidth]{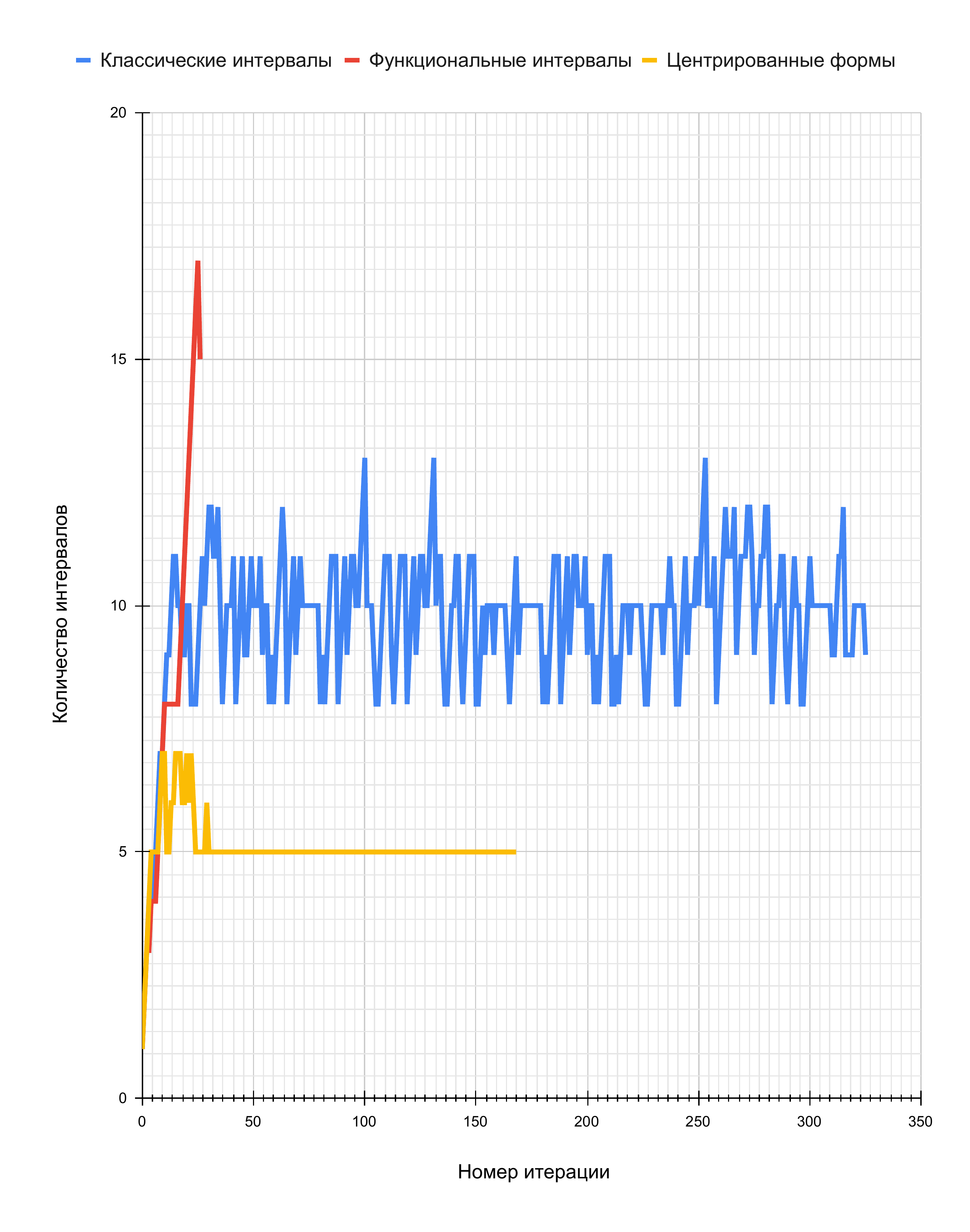}
		\caption{Количество рабочих интервалов в рабочем списке на каждой итерации алгоритма.}
		
	\end{center}
\end{figure}

\clearpage

\subsection{Задача нахождения минимума функции}

Рассмотрим задачу доказательного нахождения безусловного глобального минимума функции на всей вещественной оси $\mathbb{R}$. В качестве примера возьмём целевую функцию

\begin{equation*}
	f(x) = x ^ {6} + x ^ {5} - 10 x ^ {4} - 8 x ^ {3} + 15 x ^ {2} + 11.
\end{equation*}

Для применения алгоритмов численного решения задачи, ограничим область поиска интервалом, который гарантированно содержит в себе глобальный минимум. Для этого заметим, что старший член полинома чётный, а коэффициент перед ним положительный. Следовательно, данный полином ограничен снизу, так как при неограниченной возрастании модуля аргумента, значение полинома также увеличивается. Значит, такой интервал существует, а глобальный минимум будет являться одним из локальных минимумов исходной функции.

Каждый локальный минимум является точкой экстремума, поэтому необходимо рассмотреть нули производной исходного полинома

\begin{equation*}
	f'(x) = 6 x ^ {5} + 5 x ^ {4} - 40 x ^ {3} - 24 x ^ {2} + 30 x.
\end{equation*}

Для этого полинома можно с помощью следствия из теоремы Руше \cite{Prasolov} получить интервал, который гарантированно содержит в себе все его корни. Это интервал
\begin{equation*}
	\begin{array}{l}
		[ \, - 1 - \text{max} \{ \, 5 \, / \, 6, \, 40 \, /\,\, 6, \, 24 \, / \, 6, \, 30 \, / \, 6, \, 11 \, / \, 6 \, \}, \\
		\hspace{5cm} 1 + \text{max} \{ \, 5 \, / \, 6, \, 40 \, / \, 6, \, 24 \, / \, 6, \, 30 \, / \, 6, \, 11 \, /\,\, 6 \, \} \, ] = \\
		= [ \, -1 - 40 \, / \, 6, \, 1 + 40 \, / \, 6 \, ] = [ \, - 46 \, / \, 6, \, 46 \, / \, 6 \, ] \subset [ \, -7.67, \, 7.67 \, ].
	\end{array}
\end{equation*}

\clearpage

Реализуем данный алгоритм с интервалами семейств $\mathbb{L}\mathbb{F}\mathbb{R}(x)$ и $\mathbb{I}\mathbb{R}$. Для случая классических интервалов также рассмотрим центрированные формы. 

В качестве функции для вычисления характеристики рабочего интервала была выбрана
\begin{equation*}
    -\text{wid} \, f(\mbf{I}), \qquad \mbf{I} \in \mathbb{F}\mathbb{R}.
\end{equation*}

В качестве порогового значения остановки алгоритма было взято $\varepsilon = 10 ^ {-2}$.

Программа была написана на языке программирования \texttt{C Sharp} и запускалась на ЭВМ с 8-ми ядерным процессором \texttt{Intel Core i7-3770}.

\clearpage

\subsubsection{Иллюстрации итераций с интервалами $\mathbb{I}\mathbb{R}$}

Приведём иллюстрации для итераций алгоритма, в котором использовались интервалы семейства $\mathbb{I}\mathbb{R}$.

\begin{figure}[ht]
	\begin{center}
		
		\includegraphics[width = 0.45 \linewidth]{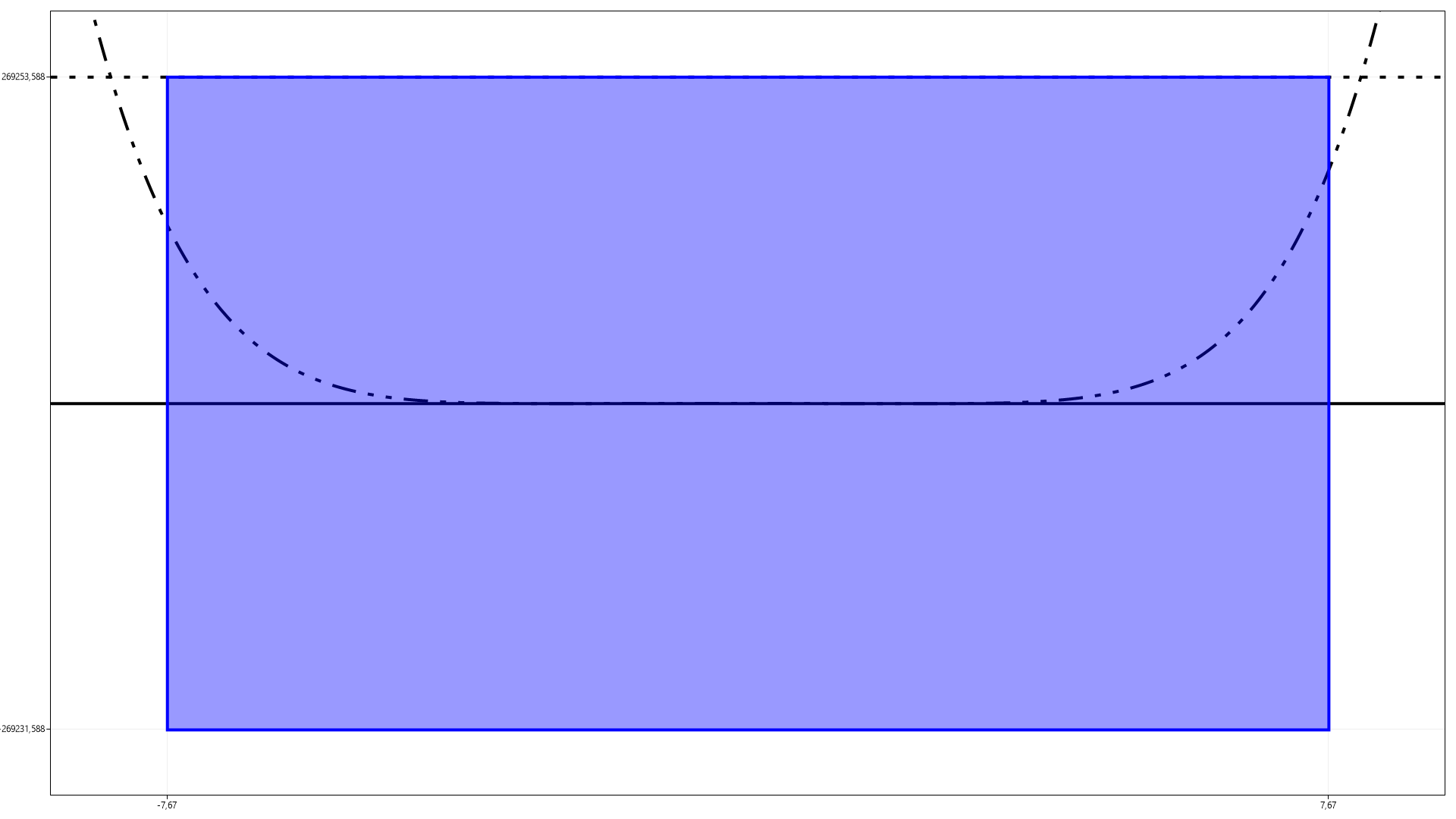}
		\includegraphics[width = 0.45 \linewidth]{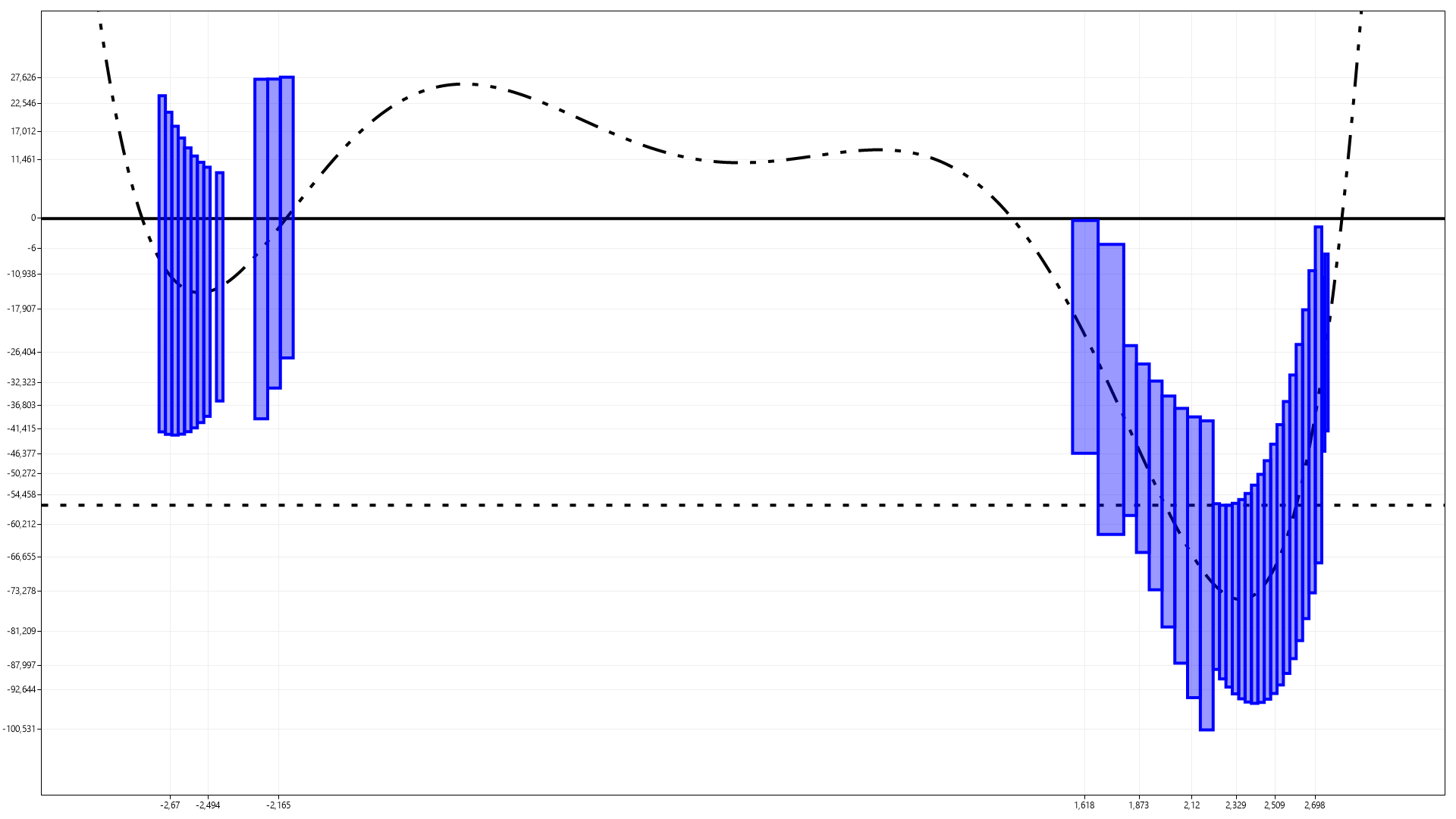}
		\caption{Итерации 0 и 112.}
		\label{fig:minimize_classical_iteration_0_112}
	
	\end{center}
\end{figure}

\begin{figure}[ht]
	\begin{center}
		
		\includegraphics[width = 0.45 \linewidth]{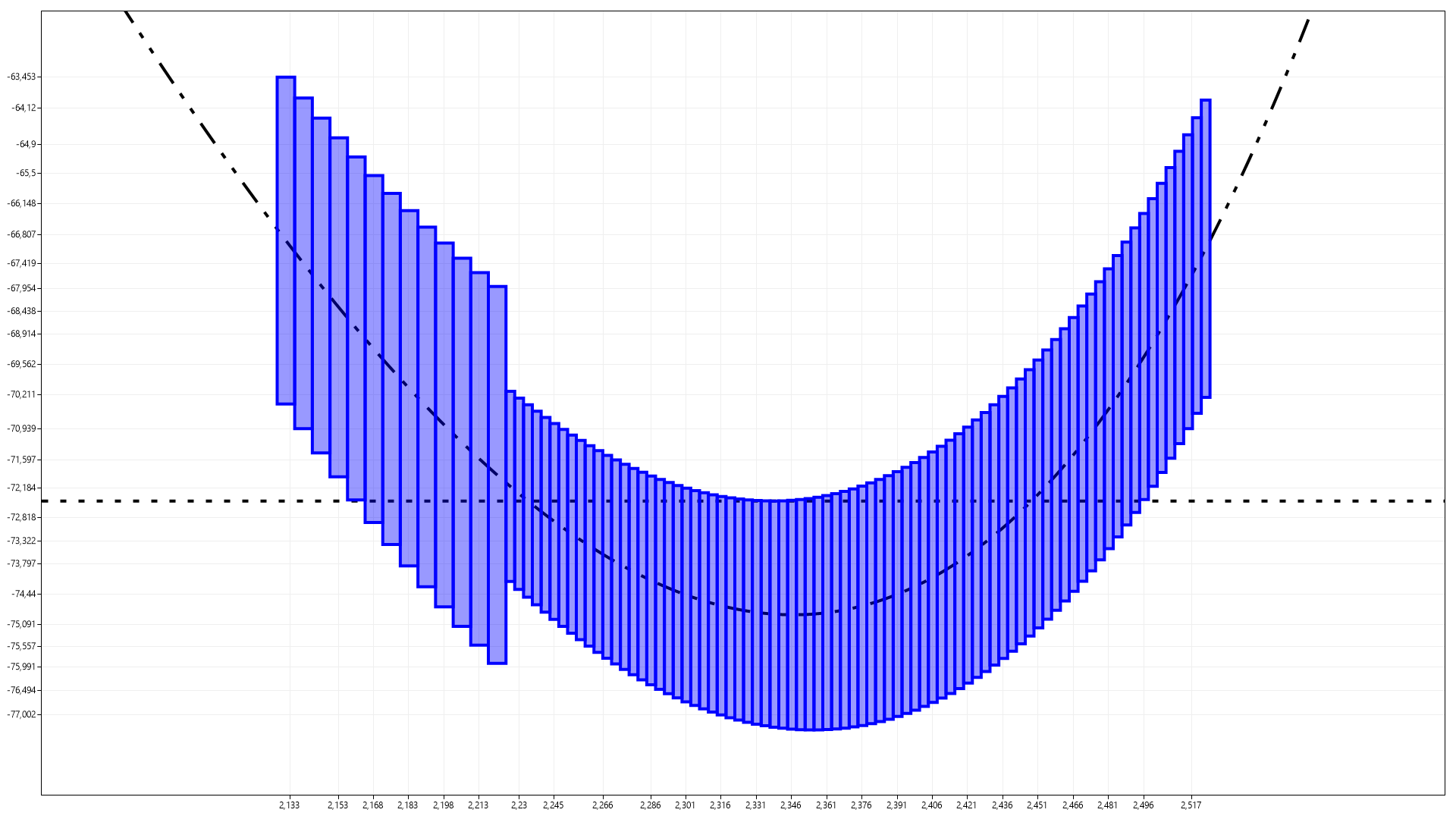}
		\includegraphics[width = 0.45 \linewidth]{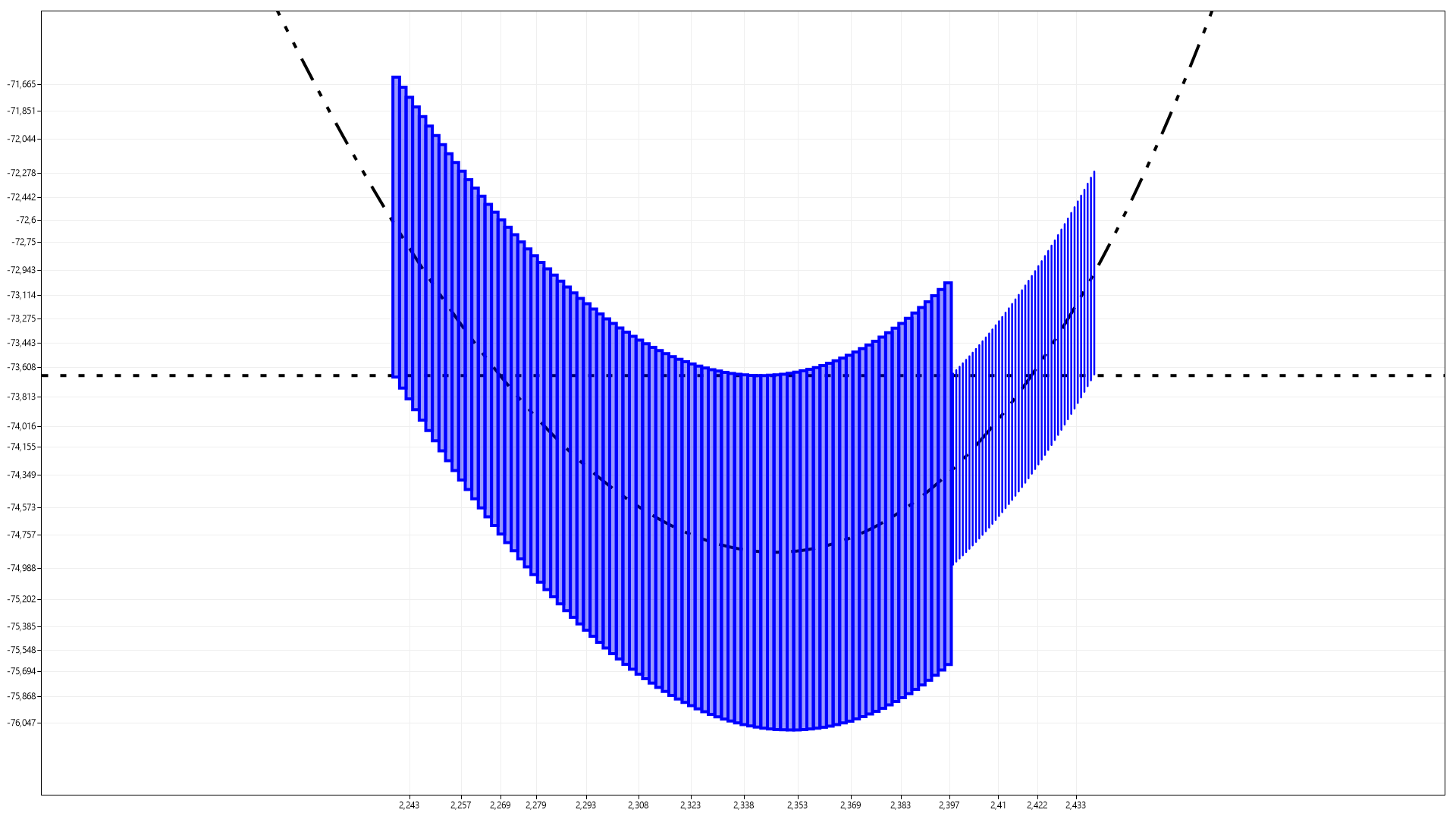}
		\caption{Итерации 224 и 336.}
		\label{fig:minimize_classical_iteration_224_336}
	
	\end{center}
\end{figure}

\begin{figure}[ht]
	\begin{center}
		
		\includegraphics[width = 0.45 \linewidth]{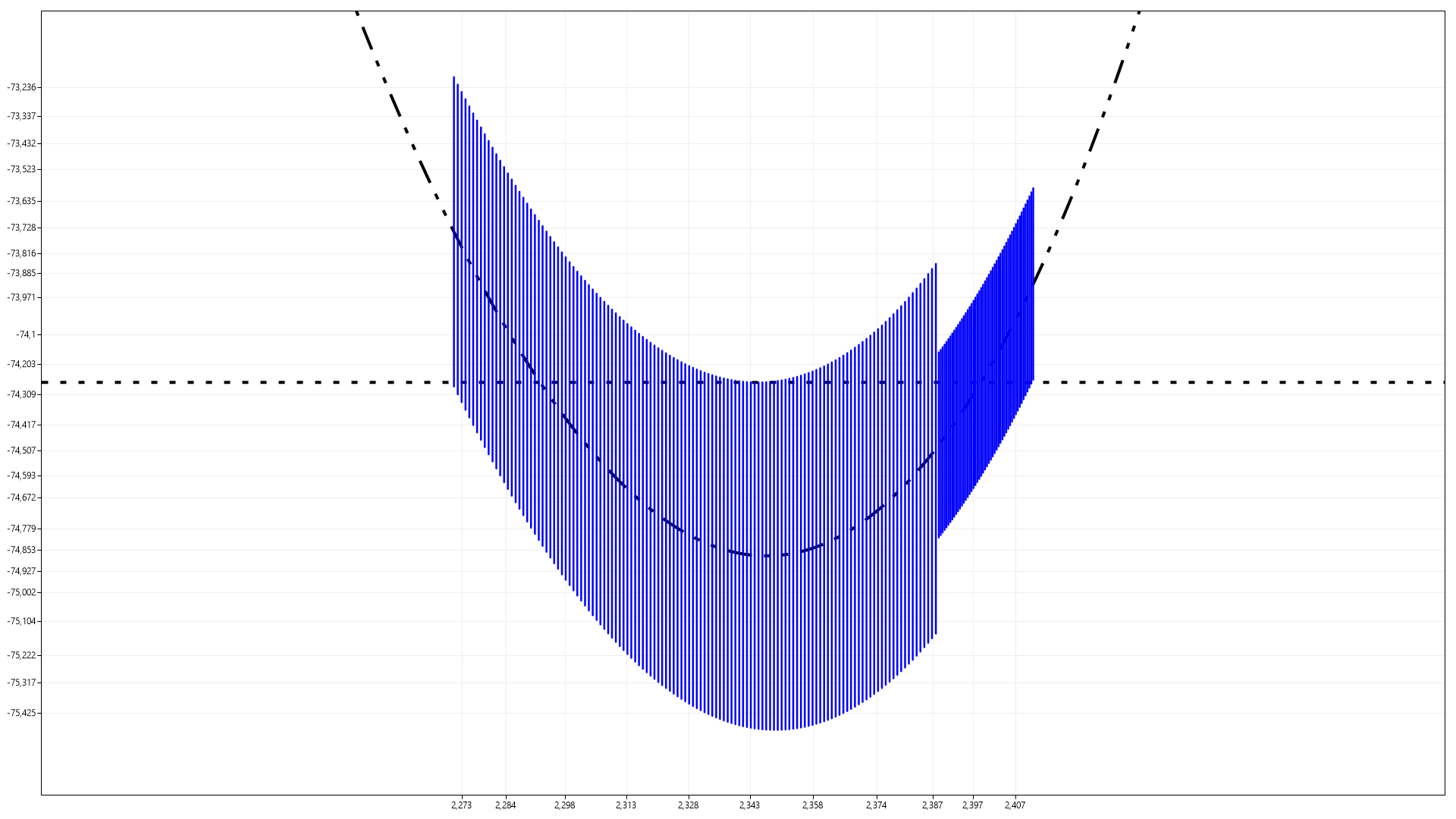}
		\includegraphics[width = 0.45 \linewidth]{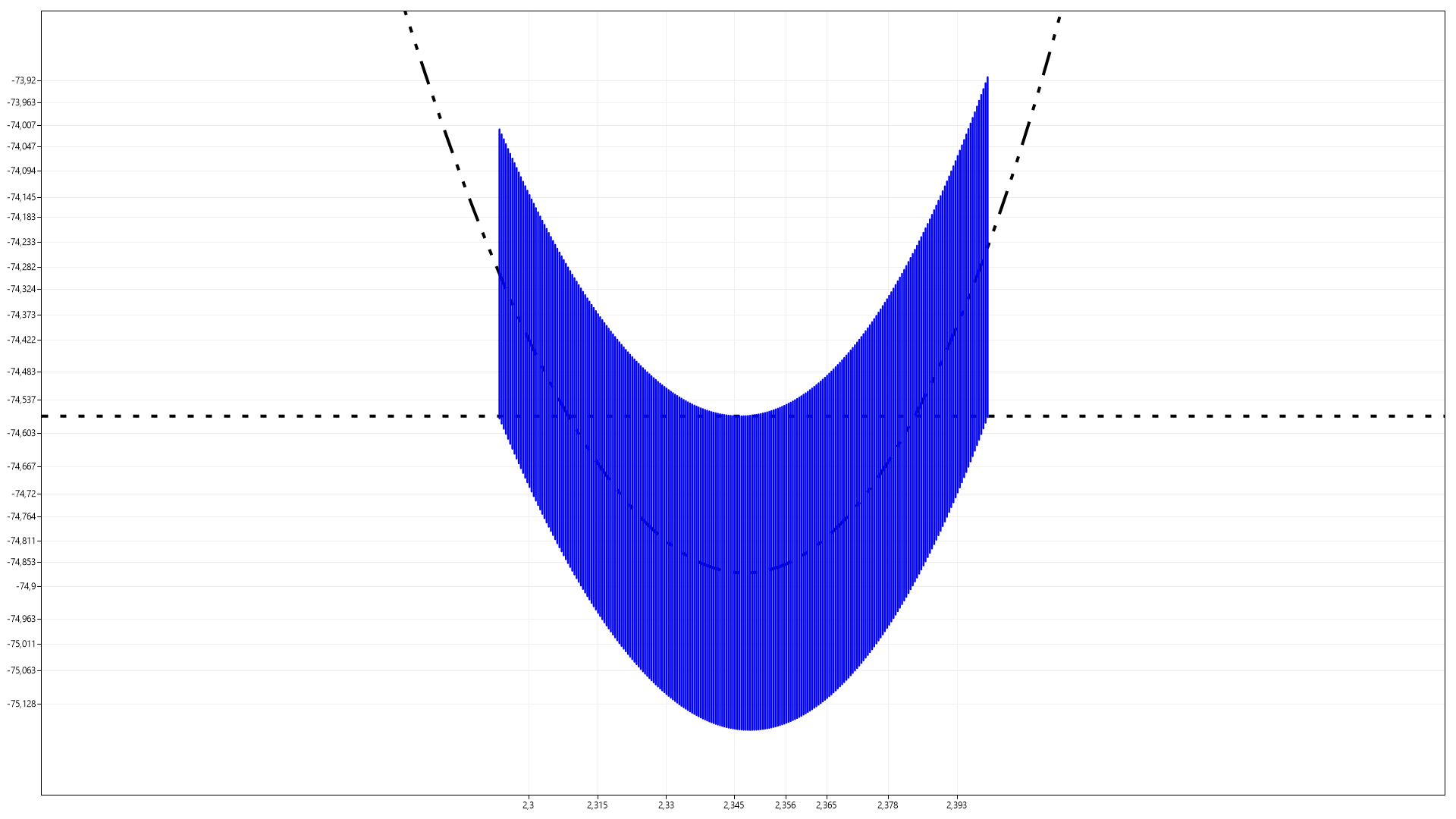}
		\caption{Итерации 448 и 563 --- последняя итерация алгоритма.}
		\label{fig:minimize_classical_iteration_448_563}
	
	\end{center}
\end{figure}

\clearpage

\subsubsection{Иллюстрации итераций с центрированными формами\\ интервалов $\mathbb{I}\mathbb{R}$}

Также приведём иллюстрации итераций алгоритма, в котором использовались центрированные формы в $\mathbb{I}\mathbb{R}$.

\begin{figure}[ht]
	\begin{center}
		
		\includegraphics[width = 0.45 \linewidth]{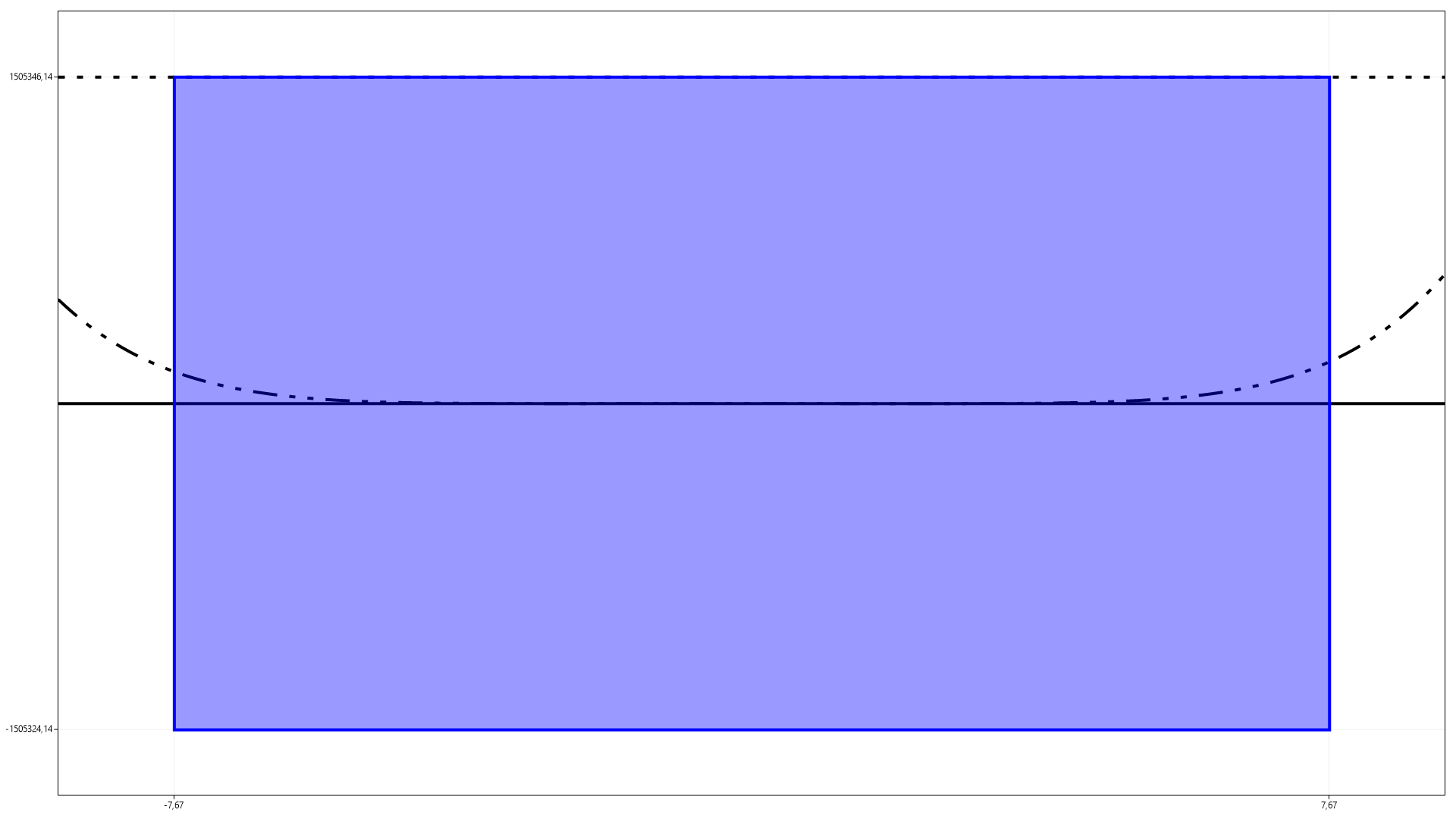}
		\includegraphics[width = 0.45 \linewidth]{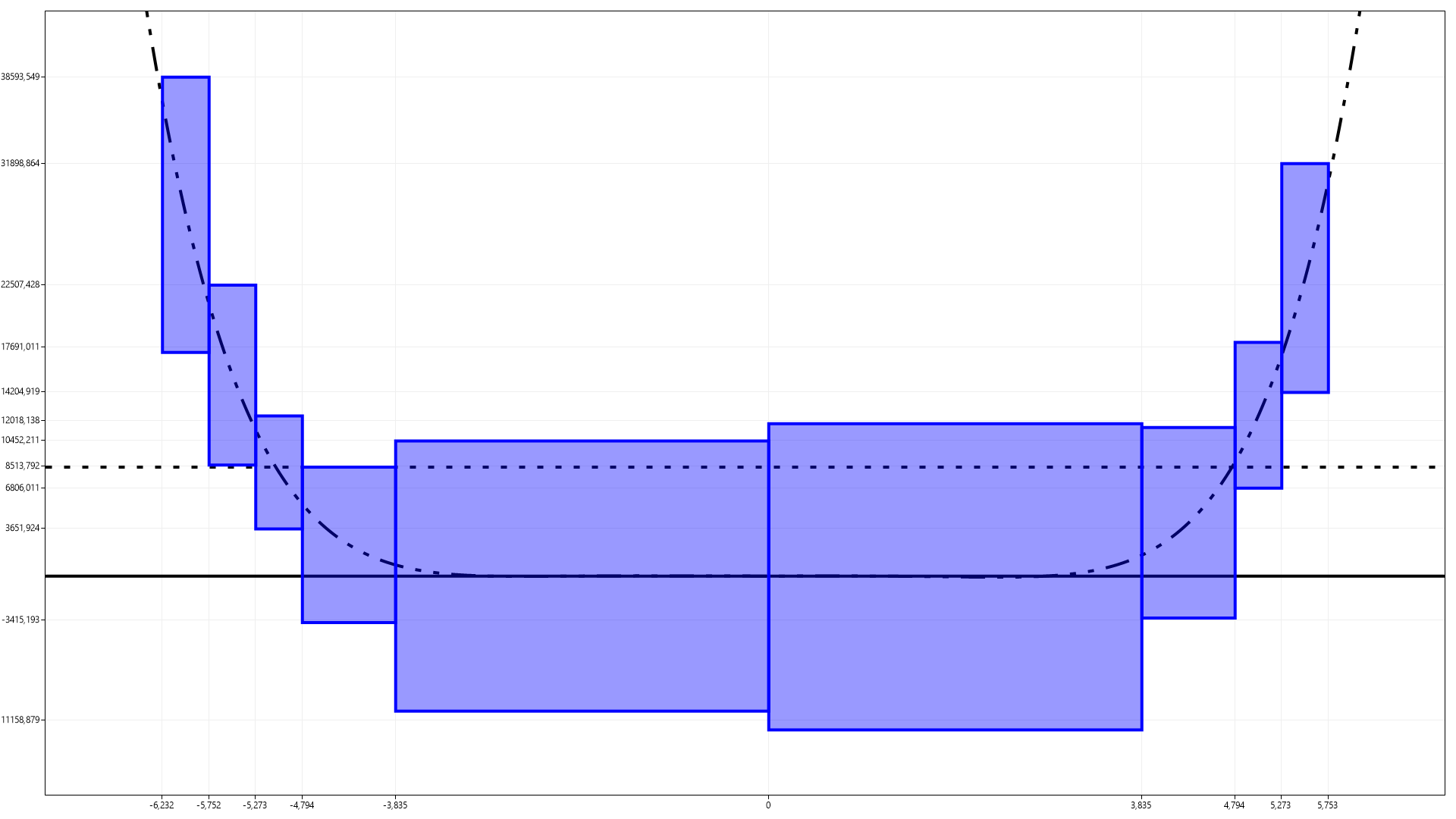}
		\caption{Итерации 0 и 12.}
		\label{fig:minimize_center_iteration_0_12}
	
	\end{center}
\end{figure}

\begin{figure}[ht]
	\begin{center}
		
		\includegraphics[width = 0.45 \linewidth]{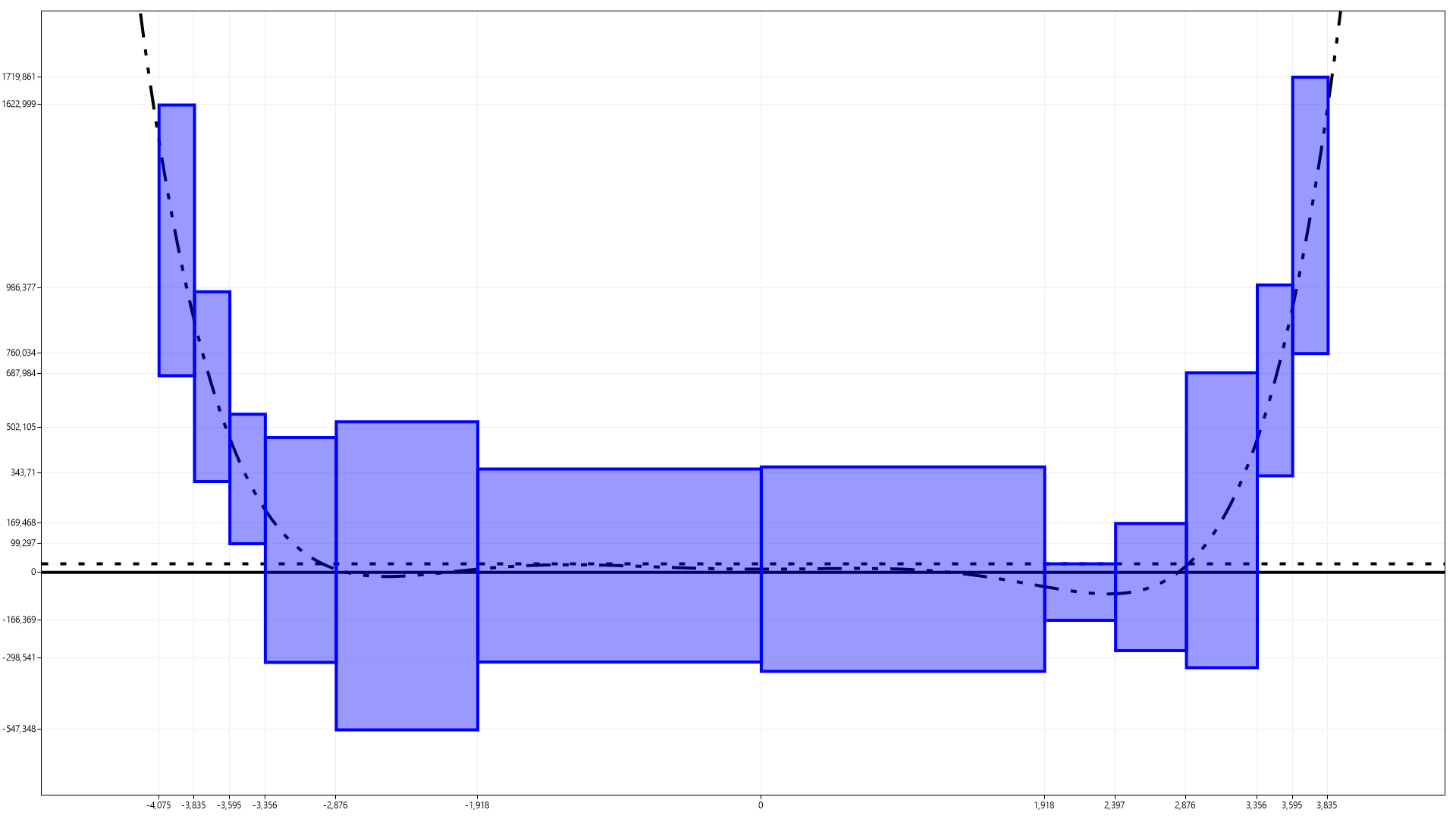}
		\includegraphics[width = 0.45 \linewidth]{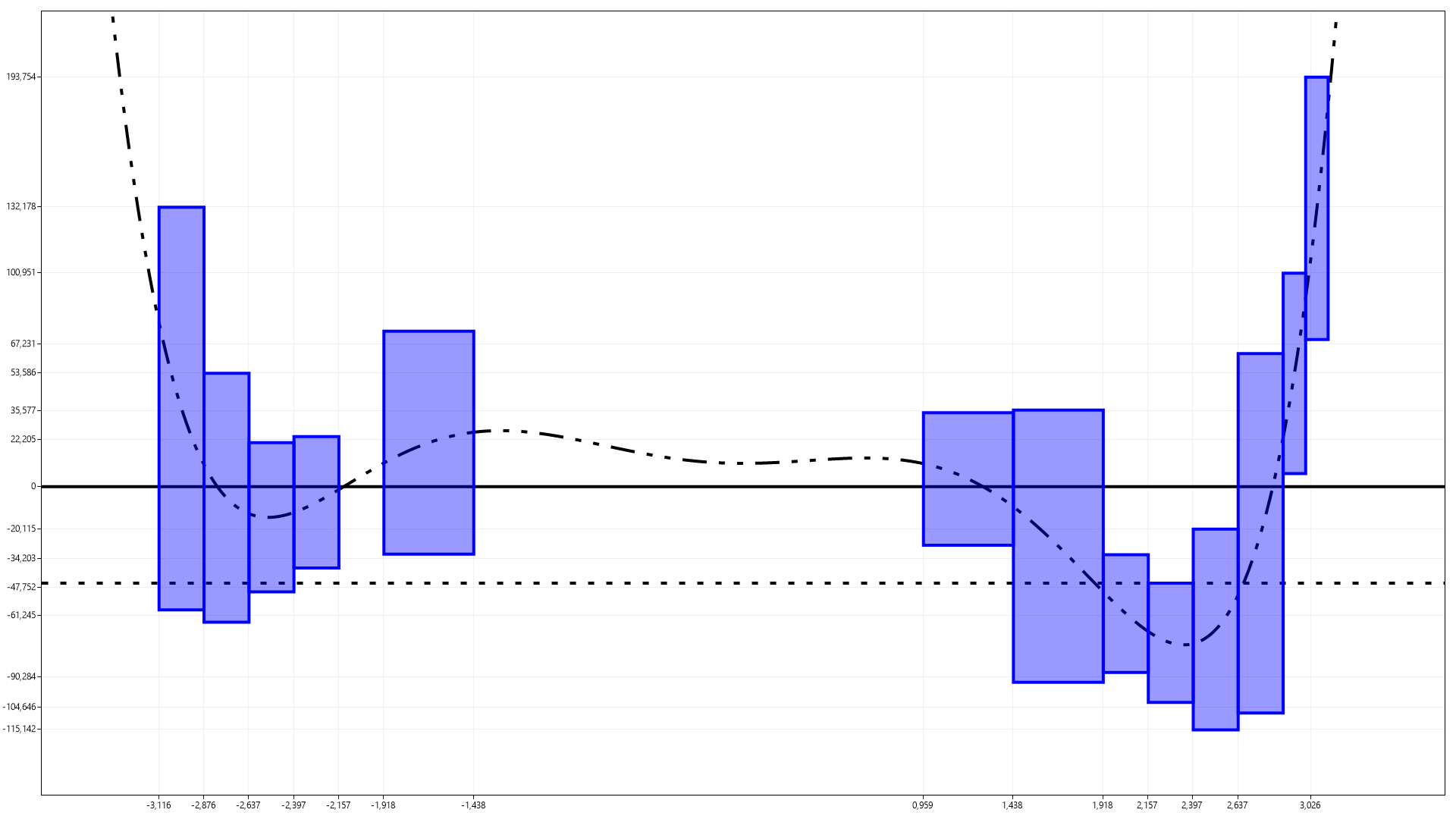}
		\caption{Итерации 24 и 36.}
		\label{fig:minimize_center_iteration_24_36}
	
	\end{center}
\end{figure}

\begin{figure}[ht]
	\begin{center}
		
		\includegraphics[width = 0.45 \linewidth]{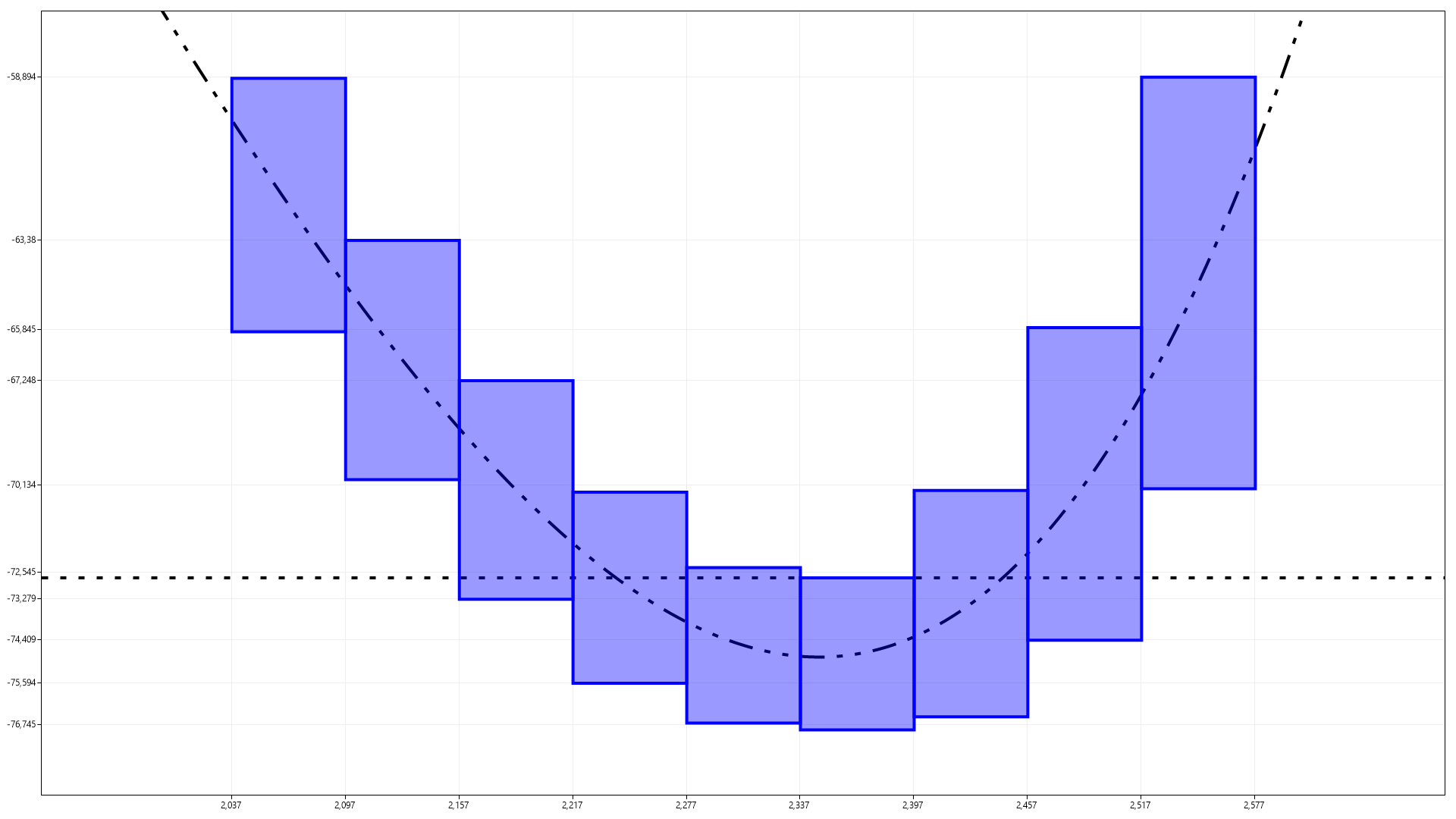}
		\includegraphics[width = 0.45 \linewidth]{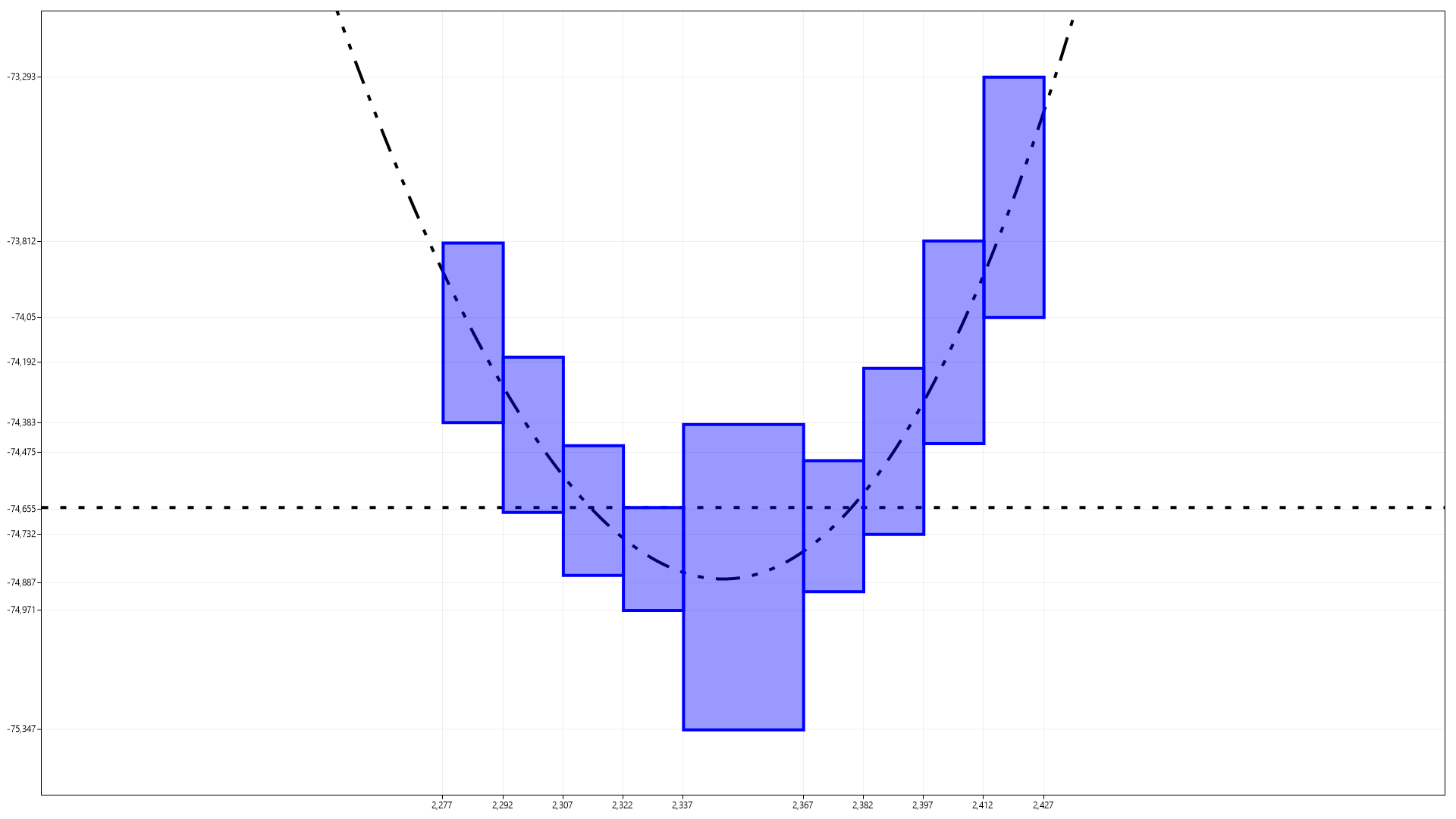}
		\caption{Итерации 48 и 59 --- последняя итерация алгоритма.}
		\label{fig:minimize_center_iteration_48_59}
	
	\end{center}
\end{figure}

\clearpage

\subsubsection{Иллюстрации итераций с интервалами $\mathbb{L}\mathbb{F}\mathbb{R}$}

Далее приведены иллюстрации для некоторых итераций алгоритма минимизации функции с использованием интервалов семейства $\mathbb{L}\mathbb{F}\mathbb{R}$.

\begin{figure}[ht]
	\begin{center}
		
		\includegraphics[width = 0.45 \linewidth]{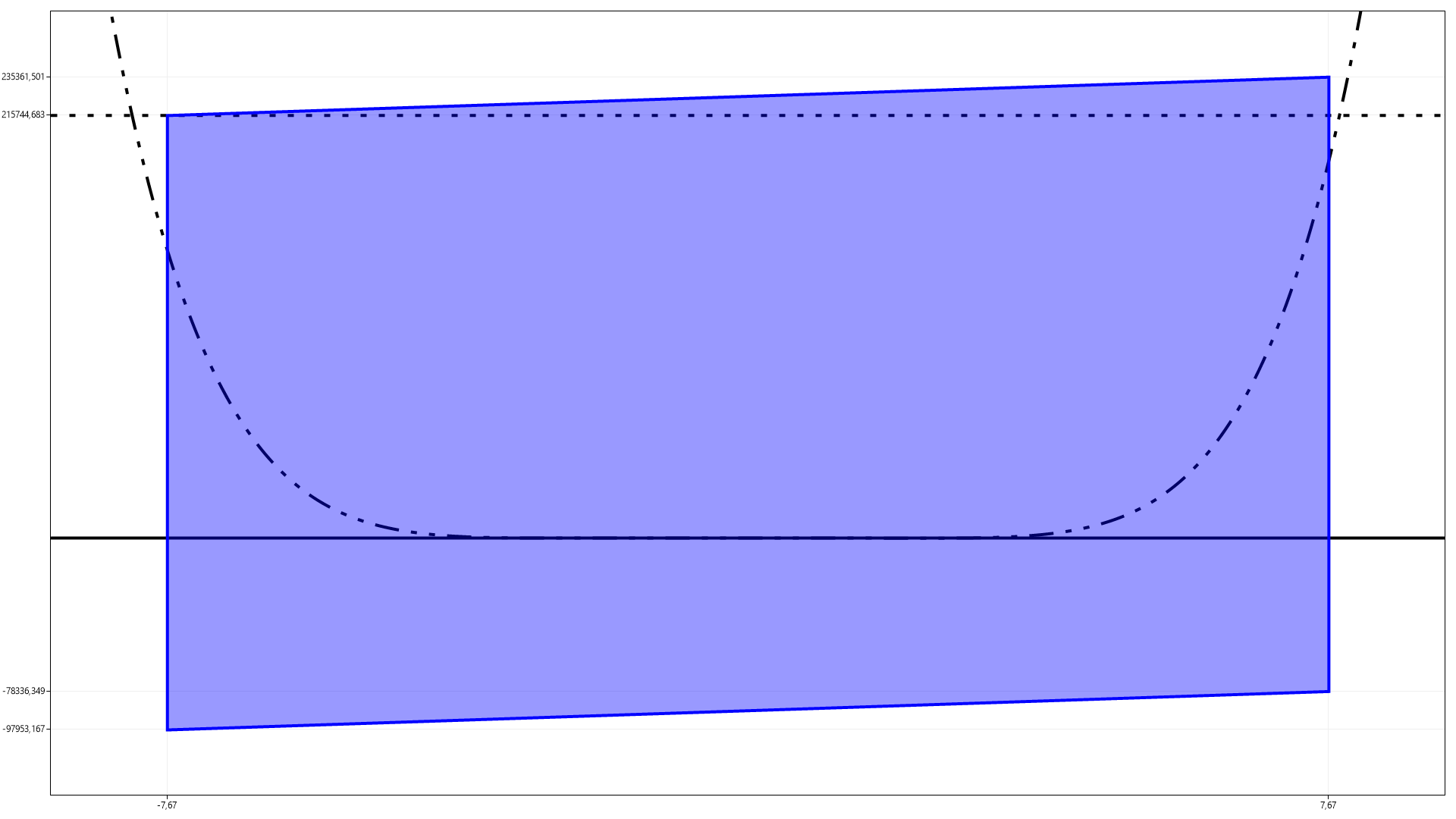}
		\includegraphics[width = 0.45 \linewidth]{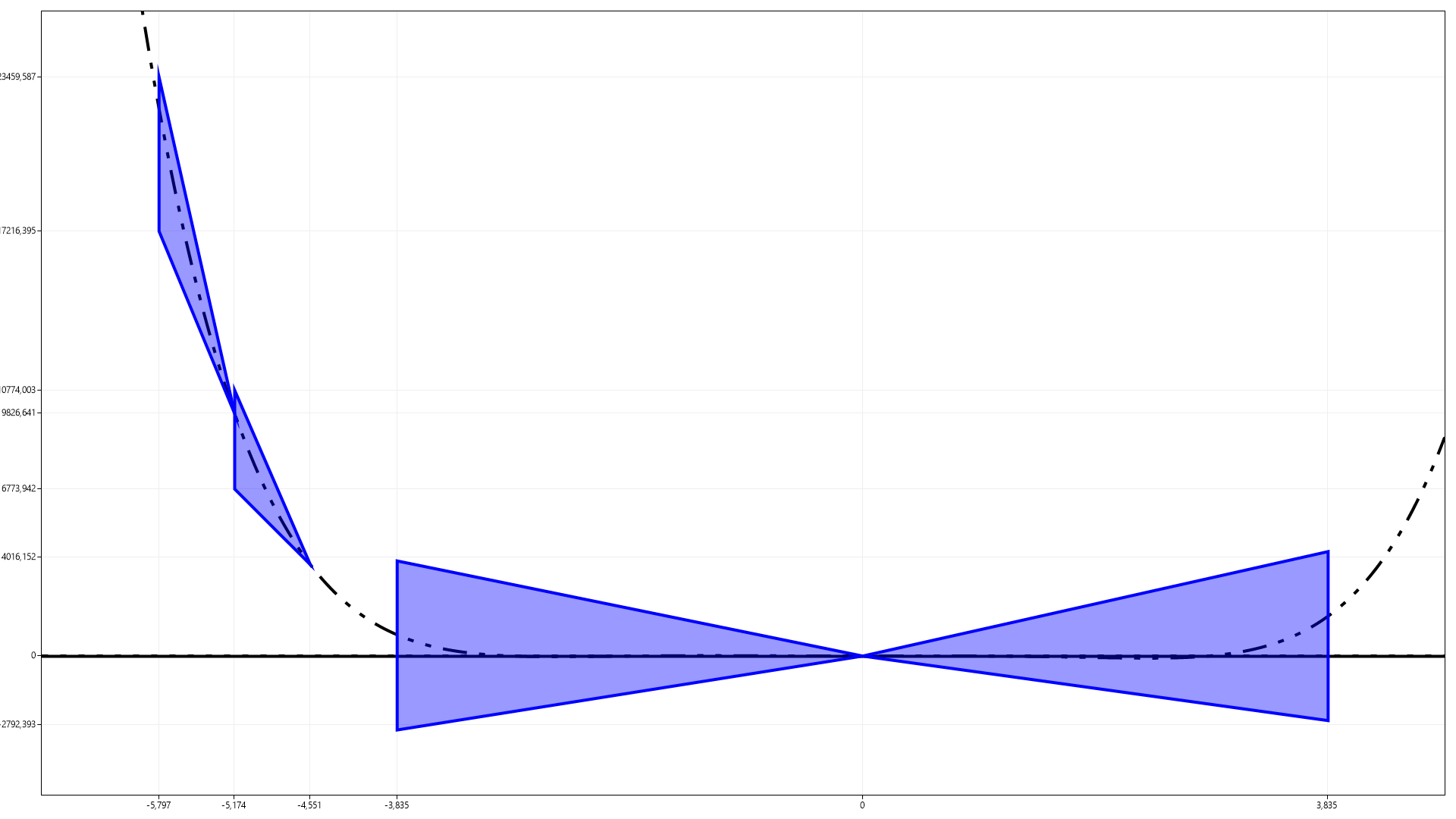}
		\caption{Итерации 0 и 6.}
		\label{fig:minimize_functional_iteration_0_6}
	
	\end{center}
\end{figure}

\begin{figure}[ht]
	\begin{center}
		
		\includegraphics[width = 0.45 \linewidth]{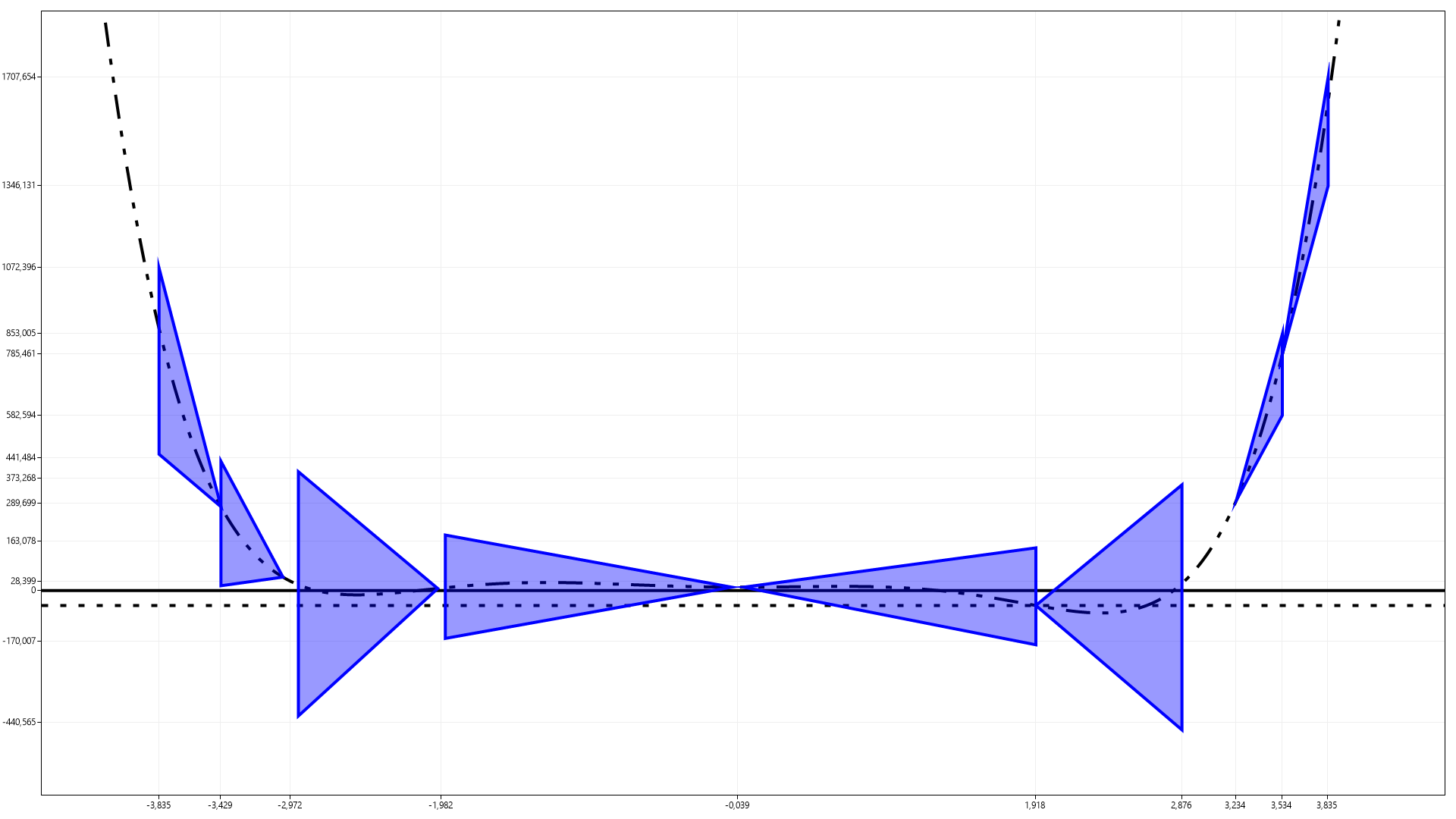}
		\includegraphics[width = 0.45 \linewidth]{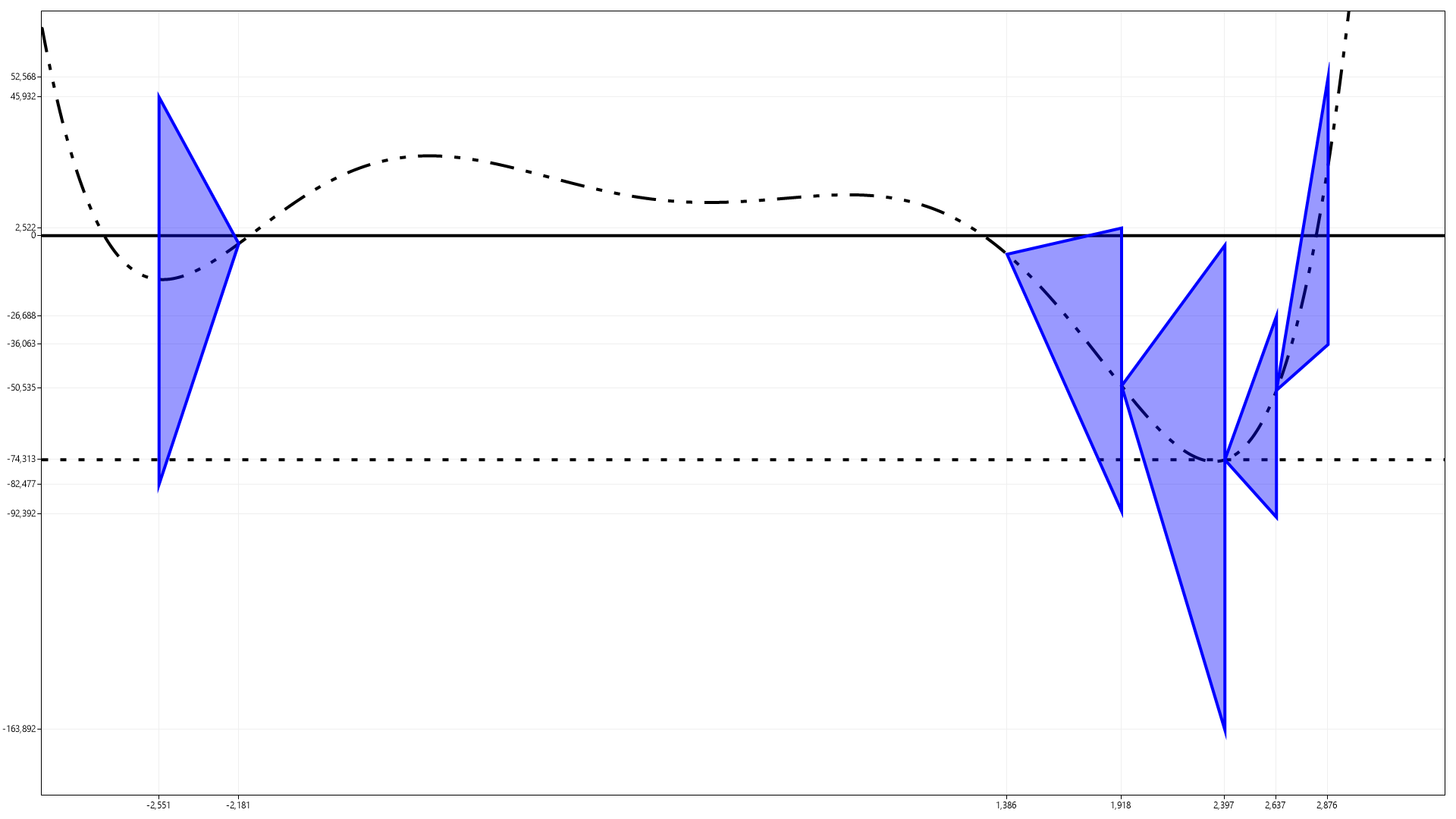}
		\caption{Итерации 12 и 18.}
		\label{fig:minimize_functional_iteration_12_18}
	
	\end{center}
\end{figure}

\begin{figure}[ht]
	\begin{center}
		
		\includegraphics[width = 0.45 \linewidth]{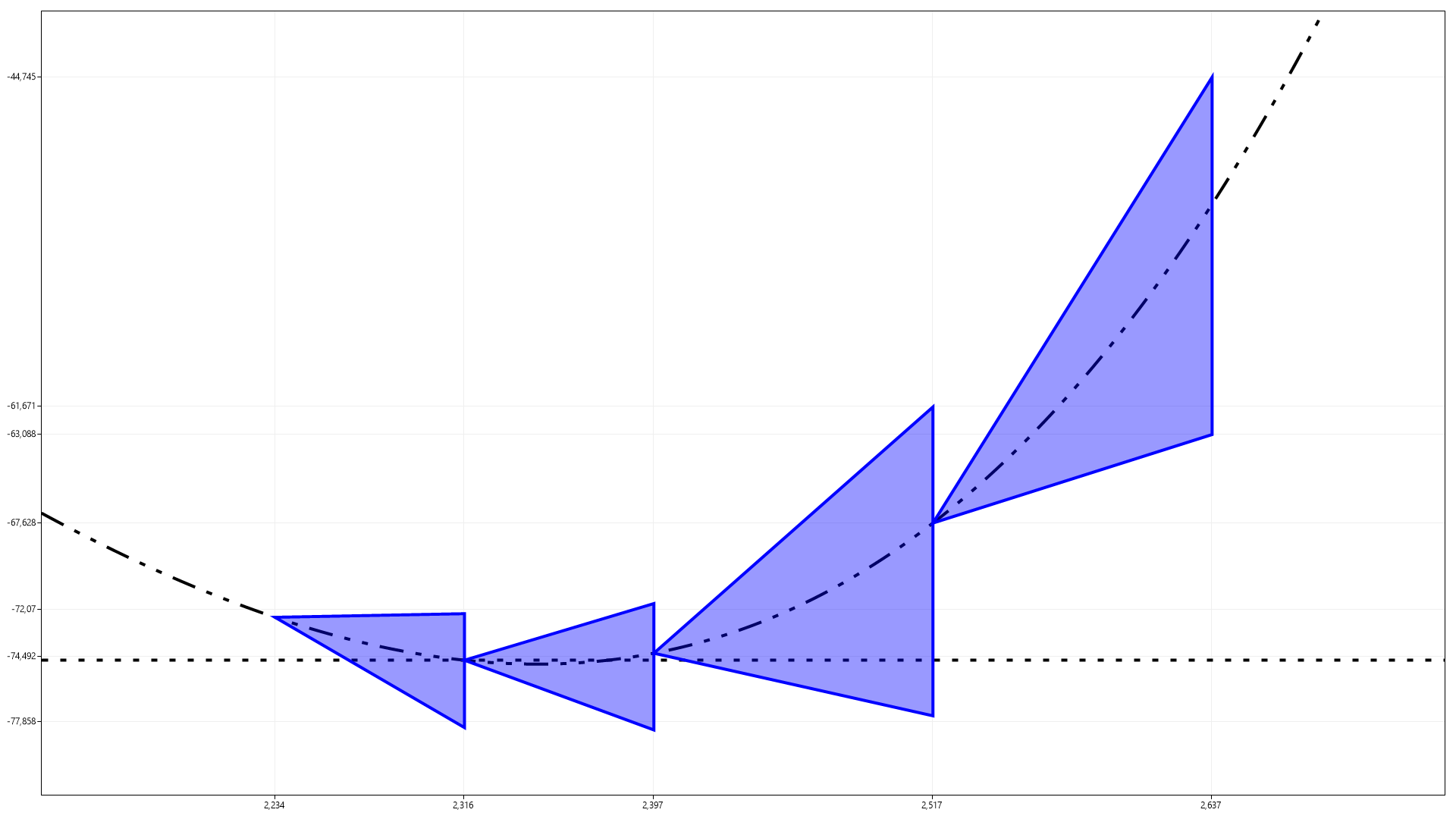}
		\includegraphics[width = 0.45 \linewidth]{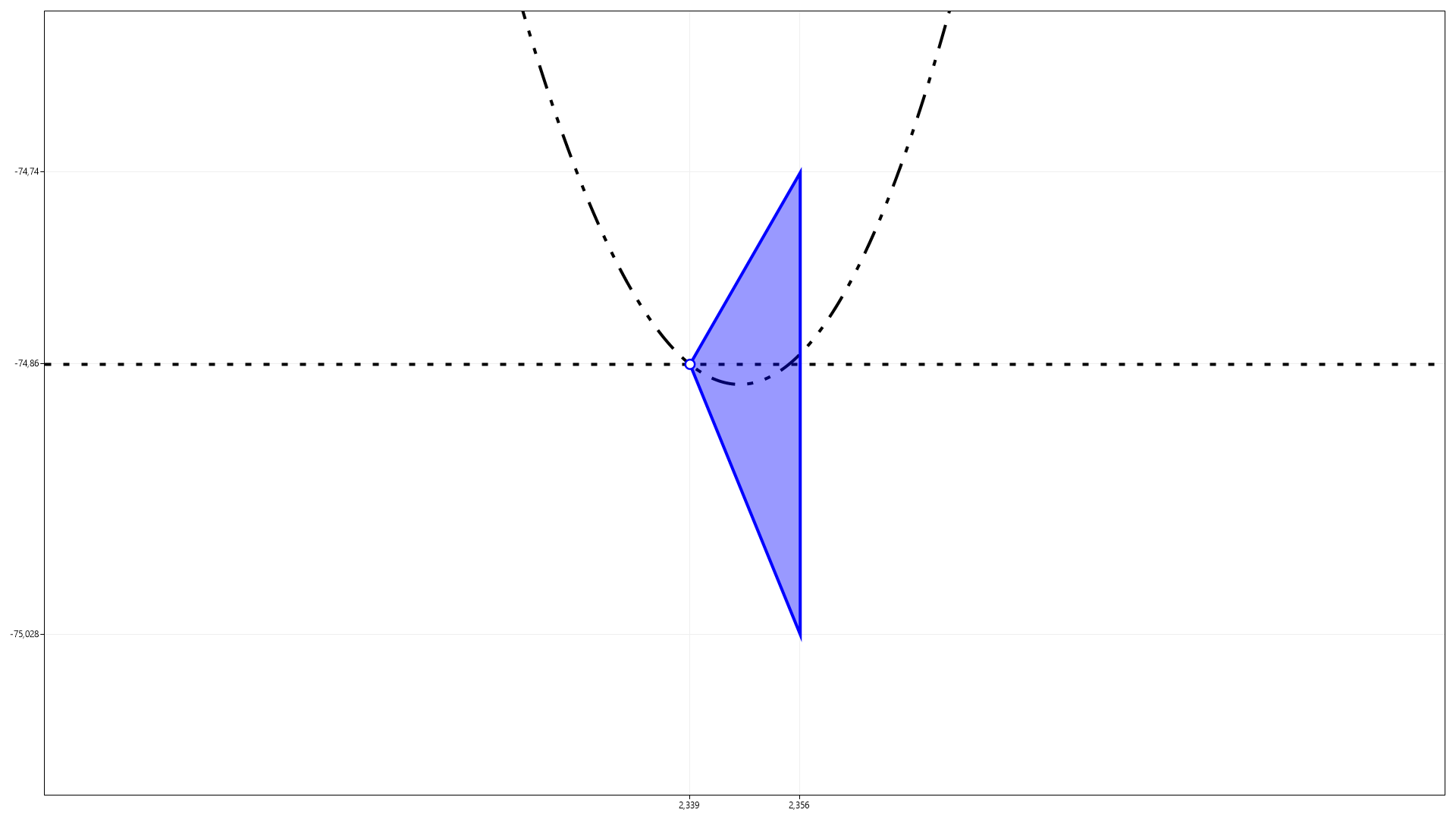}
		\caption{Итерации 24 и 30 --- последняя итерация алгоритма.}
		\label{fig:minimize_functional_iteration_24_30}
	
	\end{center}
\end{figure}

\clearpage

\subsubsection{Сравнение результатов}

Далее приведён график по по результатам численного эксперимента.

\begin{figure}[ht]
	\begin{center}
		
		\includegraphics[width = 0.8 \linewidth]{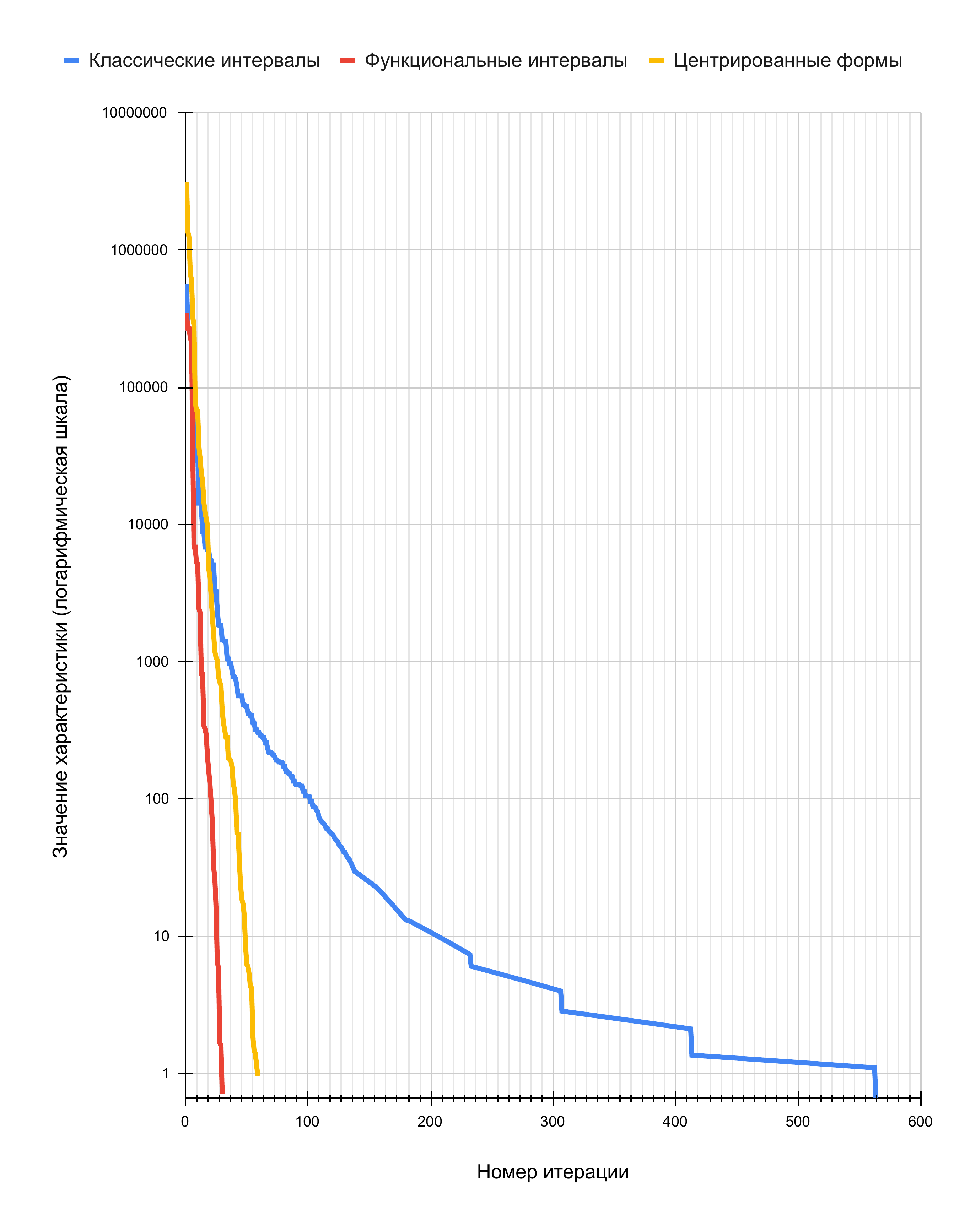}
		\caption{График сравнения максимальной характеристики рабочего интервала в рабочем списке для разных видов интервалов на каждой итерации.}
		
	\end{center}
\end{figure}

\clearpage

\subsubsection{Дополнительные результаты}

Для проведения дополнительных исследований алгоритм был запущен с функцией вычисления характеристики рабочего интервала
\begin{equation*}
    -\text{wid} \, f(\mbf{I}), \qquad \mbf{I} \in \mathbb{F}\mathbb{R}.
\end{equation*}

В качестве \textit{порогового значения остановки алгоритма} было выбрано $\varepsilon = 10 ^ {-10}$. Критерий остановки был дополнен условием на максимальное количество совершённых итераций, равное $10000$.

\clearpage

\begin{figure}[ht]
	\begin{center}
		
		\includegraphics[width = 0.8 \linewidth]{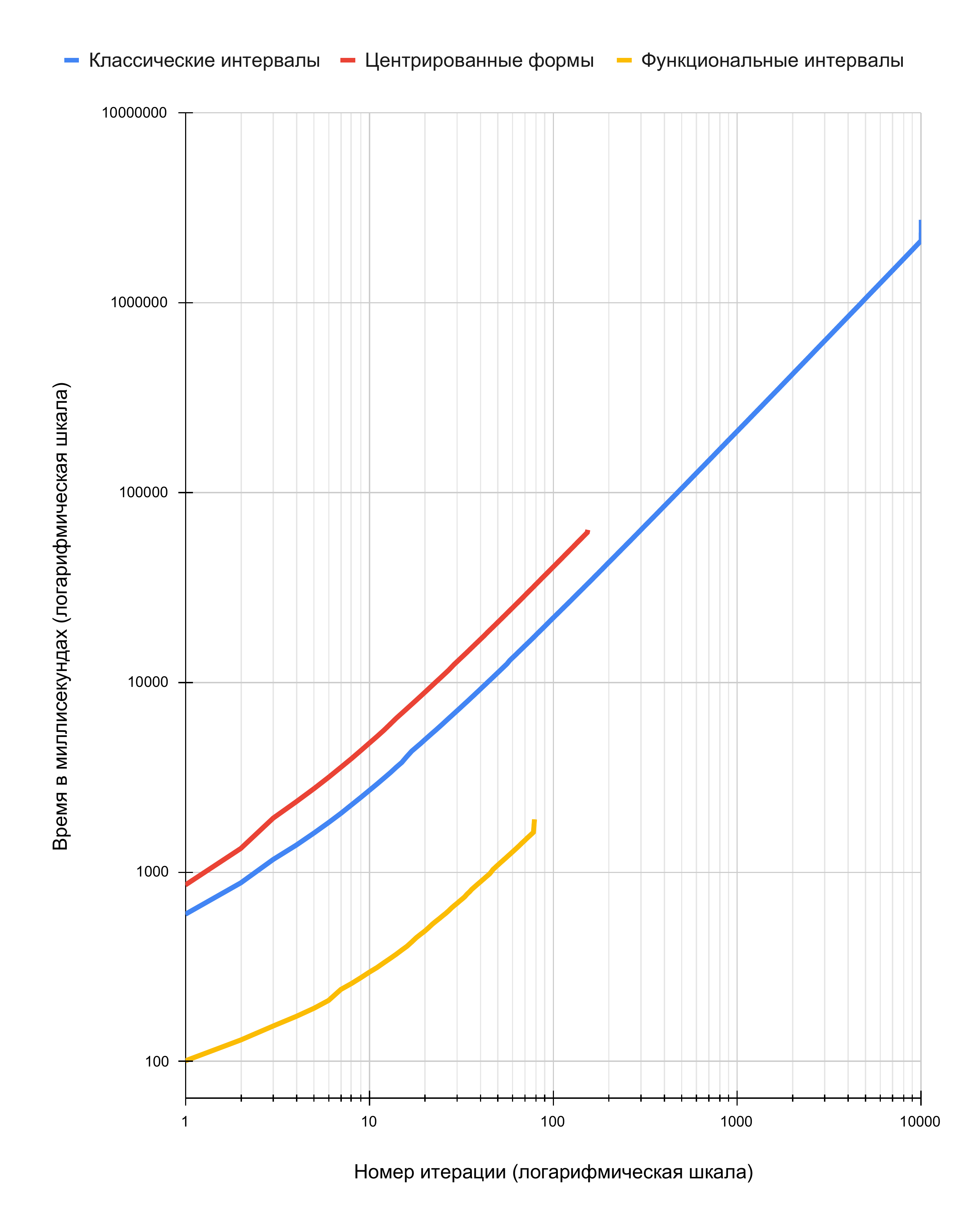}
		\caption{Суммарное время работы в миллисекундах, достигнутое на каждой итерации алгоритма.}
		
	\end{center}
\end{figure}

\begin{figure}[ht]
	\begin{center}
		
		\includegraphics[width = 0.8 \linewidth]{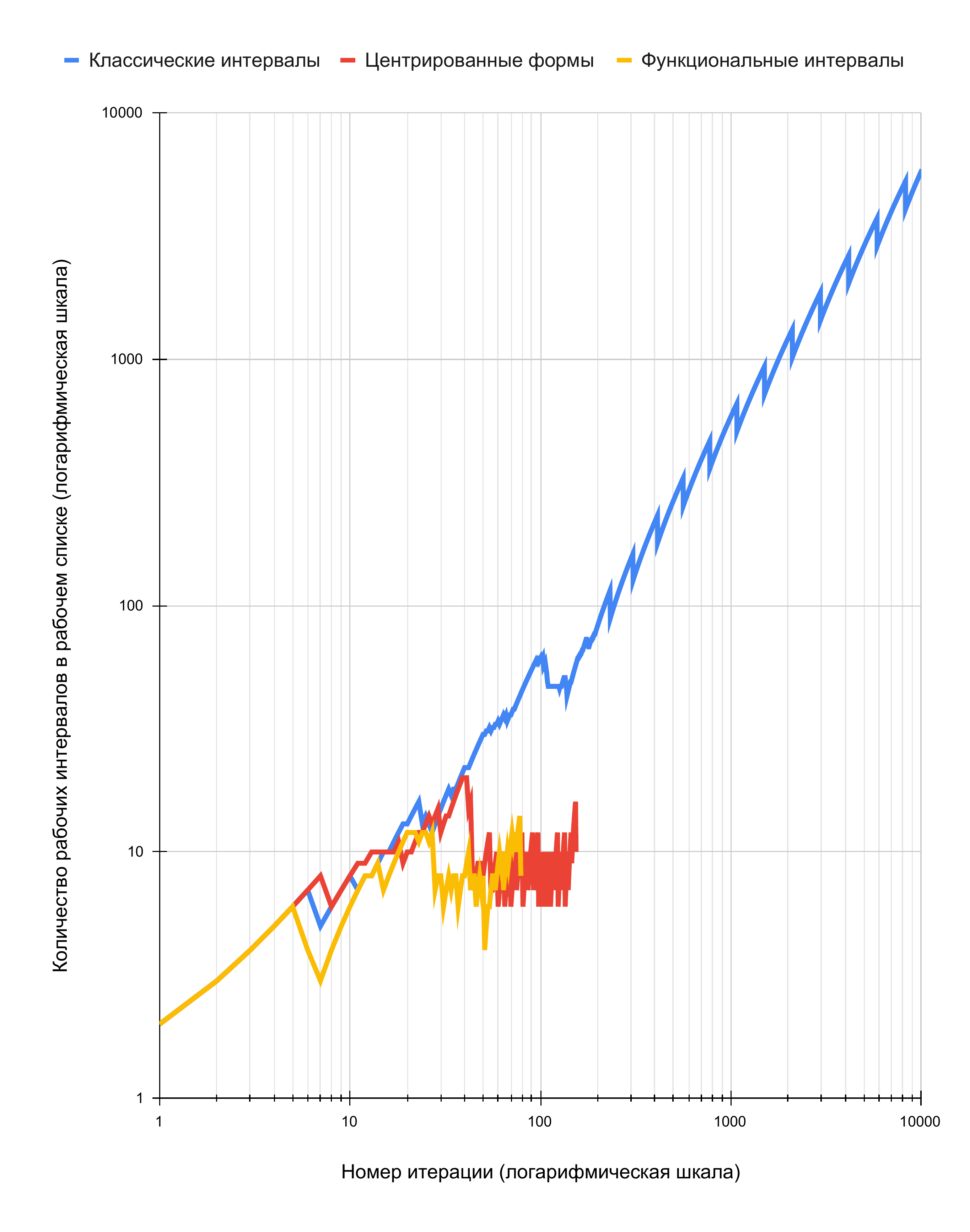}
		\caption{Количество рабочих интервалов в рабочем списке на каждой итерации алгоритма.}
		
	\end{center}
\end{figure}

\begin{figure}[ht]
	\begin{center}
		
		\includegraphics[width = 0.8 \linewidth]{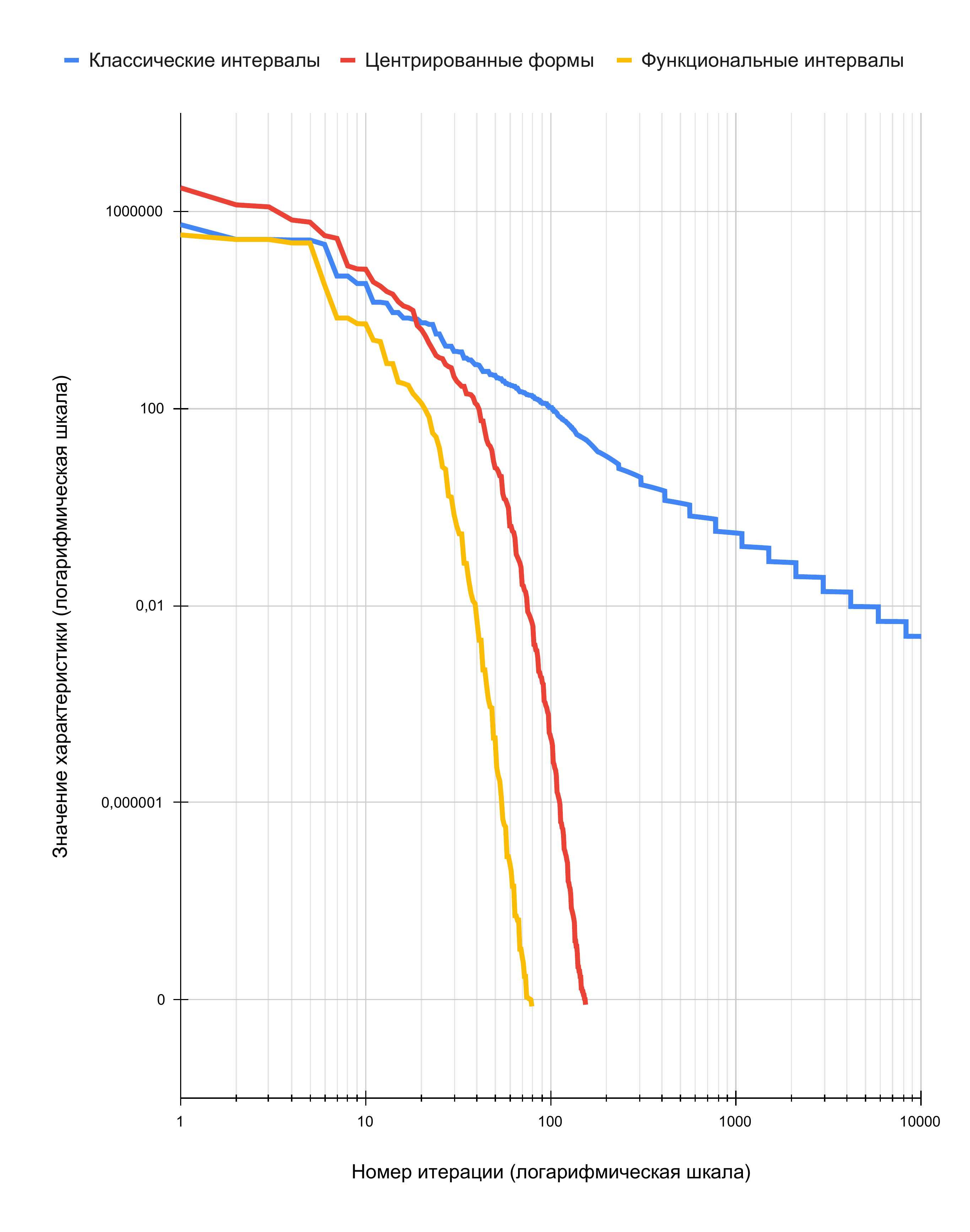}
		\caption{График сравнения максимальной характеристики рабочего интервала в рабочем списке для разных видов интервалов на каждой итерации.}
		
	\end{center}
\end{figure}

\clearpage

\subsection{Выводы}

В данной главе была построена функциональная интервальная арифметика.
Принципы, заложенные в основу построения этой арифметики, позволяют применить к этим интервалам новые различные методы для сжатия функциональных интервалов к корню или к минимуму, а также для доказательства гарантированного существования корней на этом интервале.

Была построена и подробно рассмотрена однопараметрическая линейная функциональная арифметика. Она показала свою эффективность при решении задач нахождения корней и минимума функций по сравнению с использованием классической интервальной арифметики или с использованием центрированных дифференциальных интервальных форм. Так, она вычисляет интервальную оценку функции, которая квадратично зависит от ширины области определения функции, однако при этом арифметика не использует дополнительное вычисление интервальной оценки производной функции. При этом на первых итерациях качество интервальной оценки однопараметрической линейной функциональной арифметики не хуже, чем у классической, а у квадратичных дифференциальных форм хуже.

Дальнейшее развитие данной арифметики может состоять в увеличении размерности базиса, чтобы функции $\eus{L}$ и $\eus{U}$ были линейными, но зависели от нескольких переменных. Другой путь развития состоит в использовании полиномов высоких степеней в этих функциях. 

\clearpage

\section[Многомерный вариант линейной\\ функциональной арифметики]{Многомерный вариант линейной\\ функциональной арифметики}

В предыдущих частях работы подробно рассматривалась арифметика линейных функциональных интервалов одной переменной. Однако чаще на практике используются функции многих переменных. Увеличение мощности базиса граничных функций и аналогичное рассмотрение арифметических операций между ними представляется технически сложным или невозможным, поэтому от этой идеи было решено отказаться.

Автором работы предложена идея модификации существующего метода многомерной оптимизации функции на основе алгоритма <<ветвей-и-границ>>, в котором используется квадратичная интервальная оценивающая форма. Эта модификация основана на геометрической интерпретации функциональных интервалов.

Рассмотрим эту идею подробнее.

\subsection{Теоретическая идея}

Пусть в рассмотрении имеется неравенство вида
\begin{equation}
\label{eq:unequal_1}
    \mbf{a} x \leq b,\text{ где }x \in \mbf{x} = [ \, \underline{x}, \overline{x} \, ]\text{ --- переменная, } 0 \in \mbf{a}, \, b < 0.
\end{equation}
Для определённости положим $\mbf{a} = [ \, -1, \, 2 \, ]$, $b = -1$.

Заметим, что $\mbf{a}x$ можно рассматривать как функциональный интервал, у которого 
\begin{equation*}
    \eus{L}(x) = |x| \cdot (-\text{wid} \, \mbf{a}) + x \cdot \text{mid} \, \mbf{a} \qquad \text{и} \qquad \eus{U}(x) = |x| \cdot \text{wid} \, \mbf{a} + x \cdot \text{mid} \, \mbf{a}.
\end{equation*}

\begin{figure}[ht]
	\begin{center}
		
		\includegraphics[width = 0.20 \linewidth]{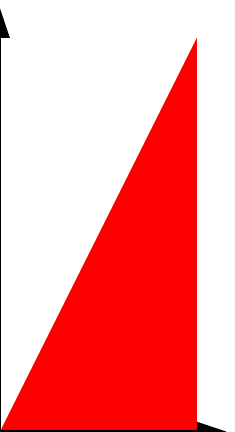}
		\caption{Множество точек, задаваемое выражением $[ \, -1, \, 2 \,  ]x \leq -1$.}
		\label{fig:label_ax_2}
	
	\end{center}
\end{figure}

Разрешим неравенство (\ref{eq:unequal_1}). На рисунке (рис. \ref{fig:label_ax_1}) показано множество точек, задаваемое выражением $\mbf{a} x$, а на рисунке (рис. \ref{fig:label_ax_2}) показано множество точек, задаваемое неравенством $\mbf{a} x \leq b$:

\begin{figure}[ht]
	\begin{center}
		
		\includegraphics[width = 0.20 \linewidth]{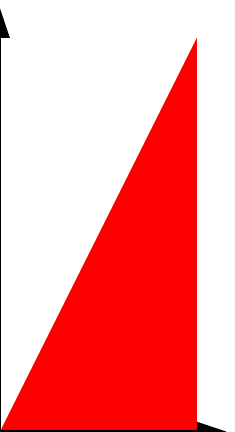}
		\caption{Множество точек, задаваемое выражением $[ \, -1, \, 2 \,  ]x$.}
		\label{fig:label_ax_1}
	
	\end{center}
\end{figure}

Видно, что можно достаточно просто определить те $x \in \mbf{x}$, где неравенство $\mbf{a} x \leq b$ не выполняется, используя представление функции $\eus{L}(x)$ для функционального интервала $\mbf{a}x$. В рассматриваемом случае неравенство не выполняется при $x \in [ \, -0.5, \, 1 \, ]$.

В общем случае, формула исключенного интервала такова:
\begin{equation*}
    \tilde{\mbf{x}} = \Big[ \, \frac{b}{\overline{a}}, \, \frac{-b}{\underline{a}} \, \Big]
\end{equation*}

Таким образом, кандидатами на решения нашего исходного неравенства (\ref{eq:unequal_1}) будут два интервала:

\begin{equation*}
    \mbf{x} \cap [ \, -\infty, \, \underline{\tilde{\mbf{x}}} \, ] \qquad \text{и} \qquad \mbf{x} \cap [ \, \overline{\tilde{\mbf{x}}}, \, +\infty \, ].
\end{equation*}

Как это можно распространить на случай неравенства многих переменных? Рассмотрим неравенство вида:
\begin{equation}
\label{eq:unequal_2}
    \sum\limits_{i = 1}^{n} \mbf{a}_{i} \, x_{i} \leq b.
\end{equation}

Заметим, что это неравенство можно рассматривать, как результат сложения неравенств вида

\begin{equation}
\label{eq:split_uneq}
    \mbf{a}_{i} \, x_{i} \leq \frac{b}{n} \quad \forall i = 1, \dots, n \quad \Rightarrow \quad \sum\limits_{i = 1}^{n} \mbf{a}_{i} \, x_{i} \leq b.
\end{equation}

Затем для каждого из этих слагаемых неравенств можно получить интервал $\tilde{\mbf{x}}_{i}$, в котором гарантированно не содержится решение неравенства для координаты $x_{i}$. Прямое произведение всех этих интервалов будет представлять собой брус $\tilde{\mbf{x}}_{1} \times \dots \times \tilde{\mbf{x}}_{n}$, в котором гарантированно не содержится решение исходного неравенства.

Итак, мы получили метод для решения неравенств (\ref{eq:unequal_2}). Рассмотрим, как получить данное неравенство из интервальной центрированной квадратичной формы для задач оптимизации функции?

Пусть у нас имеется дифференцируемая функция $f(x_{1}, \dots, x_{n})$, где $x_{1} \in \mbf{x}_{1}, \dots, x_{n} \in \mbf{x}_{n}$.
Тогда центрированная форма для данной функции будет выглядеть как
\begin{equation}
\label{eq:center_form}
    f_{c}(x_{1}, \dots, x_{n}) \in f(\text{mid} \, \mbf{x}_{1}, \dots, \text{mid} \, \mbf{x}_{n}) + \sum\limits_{k = 1}^{n} (x_{k} - \text{mid} \, \mbf{x}_{k}) \cdot f^{'}_{x_{k}}(\mbf{x}_{1}, \dots, \mbf{x}_{n}).
\end{equation}

Так как мы решаем задачу оптимизации, то случаи, когда интервальная оценка какой-либо частной производной $f^{'}_{x_{k}}(\mbf{x}_{1}, \dots, \mbf{x}_{n})$ не содержит $0$, можно сразу произвести сжатие исходного бруса 
$\tilde{\mbf{x}}_{1} \times \dots \times \tilde{\mbf{x}}_{n}$ в зависимости от знака этой производной.

Поэтому мы рассматриваем случай, когда все интервальные оценки частных производных содержат в себе $0$.

Представим, что путём бросания случайных точек в брус $\mbf{x}_{1} \times \mbf{x}_{2} \times \dots \times \mbf{x}_{n}$ мы нашли такую точку $(x_{1}^{*}, \dots, x_{n}^{*})$, для которой верно
\begin{equation*}
    f(x_{1}^{*}, \dots, x_{n}^{*}) < f(\text{mid} \, \mbf{x}_{1}, \dots, \text{mid} \, \mbf{x}_{n}).
\end{equation*}
Так как мы достигли рекордное значение минимума $f(x_{1}^{*}, \dots, x_{n}^{*})$, то можно отсеять те точки исходного бруса, в которых это значение не может быть достигнуто. То есть необходимо решить неравенство
\begin{equation}
\label{eq:min_uneq}
    \sum\limits_{k = 1}^{n} (x_{k} - \text{mid} \, \mbf{x}_{k}) \cdot f^{'}_{x_{k}}(\mbf{x}_{1}, \dots, \mbf{x}_{n}) \leq f(x_{1}^{*}, \dots, x_{n}^{*}) - f(\text{mid} \, \mbf{x}_{1}, \dots, \text{mid} \, \mbf{x}_{n}).
\end{equation}

В обозначениях неравенства (\ref{eq:unequal_2}) получаем, что
\begin{equation*}
    \begin{array}{lll}
        \mbf{a}_{i} & \leftrightarrow &  f^{'}_{x_{k}}(\mbf{x}_{1}, \dots, \mbf{x}_{n}) \ni 0, \\
        x_{i} & \leftrightarrow & x_{k} - \text{mid} \, \mbf{x}_{k}, \\
        b & \leftrightarrow & f(x_{1}^{*}, \dots, x_{n}^{*}) - f(\text{mid} \, \mbf{x}_{1}, \dots, \text{mid} \, \mbf{x}_{n}) < 0.
    \end{array}
\end{equation*}

Далее это неравенство решаем с помощью разбиения на простые неравенства (\ref{eq:split_uneq}), тем самым мы найдем брус, в котором минимум функции не содержится.

\subsection{Практическая часть}

Применим разработанную выше теорию по исключению брусов на практике. Рассмотрим задачу безусловной глобальной оптимизации дифференцируемой функции $f(x_{1}, \dots, x_{n})$, где $x_{1} \in \mbf{x}_{1}, \dots, x_{n} \in \mbf{x}_{n}$.

Будем действовать следующим образом. Будем использовать алгоритм для нахождения глобального минимума из пункта работы (\ref{min_alg}). В качестве критерия остановки зададим условие на рабочий брус
\begin{equation}
\label{eq:criteria}
    \max\limits_{i = 1, \dots, n} \text{wid} \, \mbf{x}_{i} < \varepsilon.
\end{equation}

Интервальное расширение в алгоритме будем находить для функции на интервале с помощью центрированной формы, а процедуру сжатия будем проводить с помощью решения неравенства (\ref{eq:min_uneq}). 

В таблице ниже перечислены функции, которые будут тестироваться, а также начальные 
области поиска и значение $\varepsilon$ для критерия остановки (\ref{eq:criteria}). 
  
\begin{table}[!h]
\centering
\begin{small}
    \begin{tabular}{c|c|c|c}
        $i$ & $f_{i}(x, \, y)$ & $\mbf{x} = \mbf{y}$ & $\varepsilon$ \\ \hline
        
        1 & $x^{2} \cdot \cos{y} + y ^ {2} \cdot \sin{x}$ & $[ \, -10 ^ {2}, \, 10 ^ {2} \, ]$ & $10^{-12}$ \\ \hline
        
        2 & $\sin(x + y) + (x - y) ^ {2} - 1.5x + 2.5y + 1$ & $[ \, -10 ^ {5}, \, 10 ^ {5} \, ] $ & $10 ^ {-6}$ \\ \hline
        
        3 & $- \cos{x} \cdot \cos{y} \cdot \exp(-((x - \pi) ^ {2} + (y - \pi) ^ {2}))$ & $[ \, -10 ^ {9}, \, 10 ^ {9} \, ] $ & $10^{-6}$
        
    \end{tabular}
\end{small}
\caption{Функции для тестирования}
\label{tab:functions_table}
\end{table}

Далее приведена таблица с результатами поиска глобального минимума с помощью 
классического и модифицированного алгоритмов. В ячейках таблицы приведено время 
в секундах, посчитанное, как среднее время от 20 запусков соответствующего алгоритма. 
  
\begin{table}[]
\centering
\begin{tabular}{c|c|c}
    $i$ & Классический & Модифицированный \\ \hline
    $1$ & $0.59$ с & $0.36$ с \\ \hline
    $2$ & $0.83$ с & $0.74$ с \\ \hline
    $3$ & $0.11$ с & $0.11$ с
\end{tabular}
    \caption{Время работы алгоритмов}
    \label{tab:times_table}
\end{table} 
 
  
\begin{figure}[ht]
	\begin{center}
	
		\unitlength=1mm  
		
		\begin{picture}(120, 80)
		\put(0, 0){\includegraphics[height = 80mm]{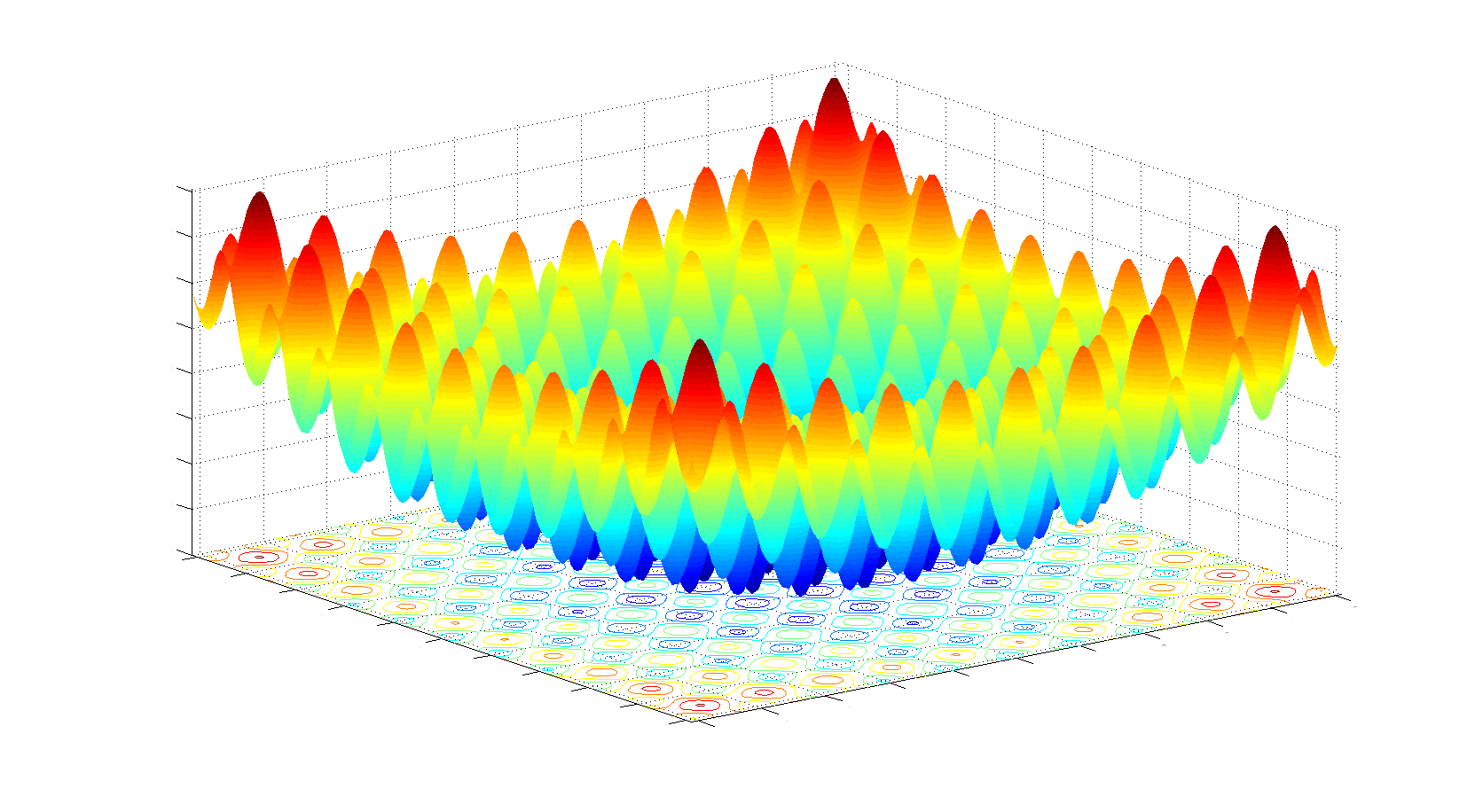}}
		
        \put(35, 10){$x_{2}$}
        \put(105, 10){$x_{1}$}
        \put(-5, 40){$f(x_{1}, x_{2})$}
        
        \end{picture} 		
		\caption{График функции Растригина для случая двух переменных.}
		\label{fig:rastrigin}
	
	\end{center}
\end{figure} 
 
\begin{figure}[ht]
	\begin{center}
		
		\unitlength=1mm  
		
		\begin{picture}(120, 80)
		\put(0, 0){\includegraphics[height = 80mm]{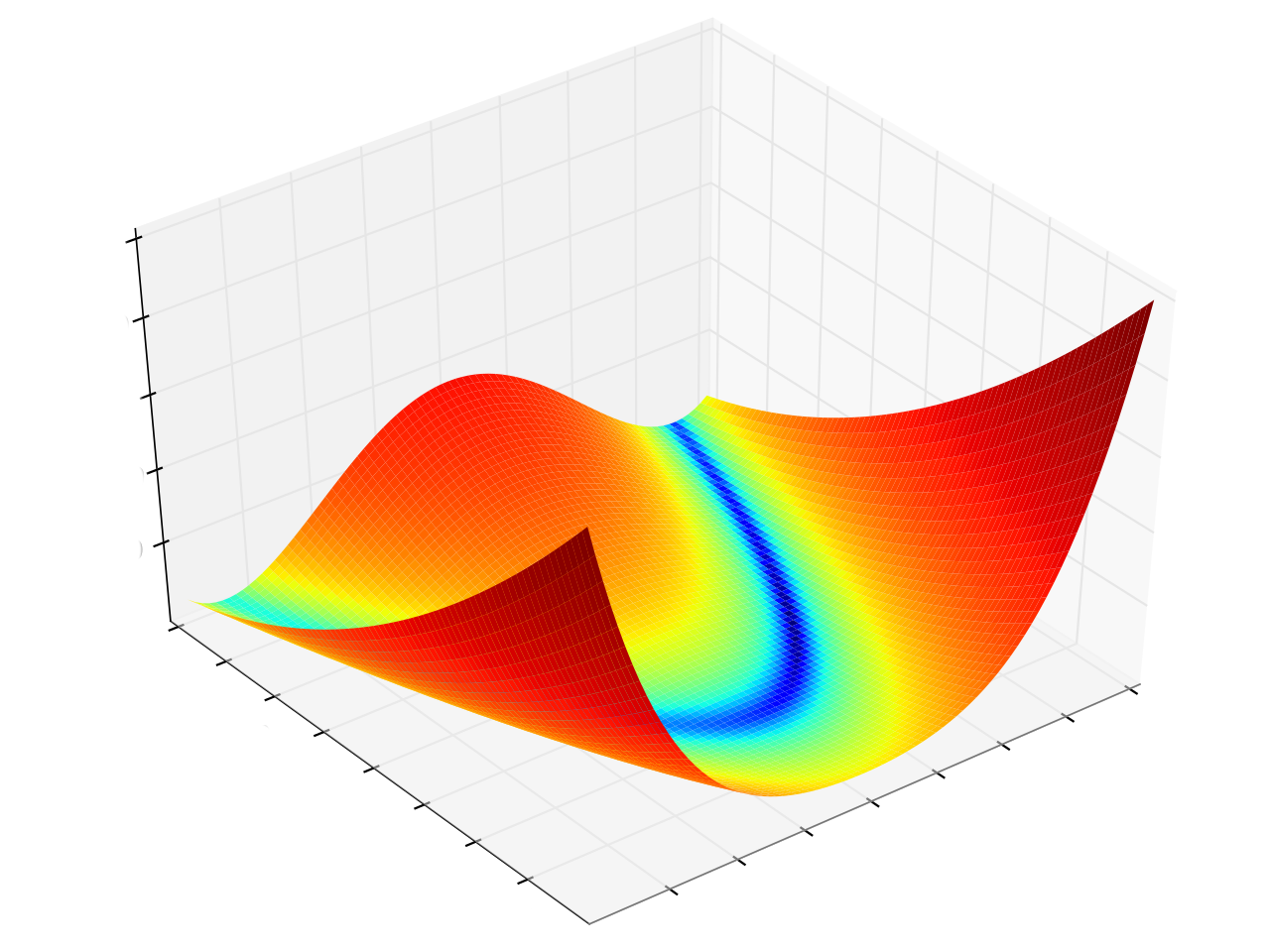}}
		
        \put(25, 10){$x_{2}$}
        \put(80, 10){$x_{1}$}
        \put(-10, 40){$f(x_{1}, x_{2})$}
        
        \end{picture}
		
		\caption{График функции Розенброка для случая двух переменных.}
		\label{fig:rozenbrock}
	
	\end{center}
\end{figure}

  
В качестве тестовых функций нескольких переменных рассмотрим известные функции 
Растригина и Розенброка. Функция Растригина (рис. \ref{fig:rastrigin}): 
\begin{equation}
\label{eq:rastrigin}
    f(x_{1}, \dots, x_{n}) = An + \sum\limits_{k = 1}^{n} (x_{i} ^ {2} - A \cos(2 \pi x_{i})), \qquad A = 10.
\end{equation}
Функция Розенброка (рис. \ref{fig:rozenbrock}):
\begin{equation}
\label{eq:rozenbrock}
    f(x_{1}, \dots, x_{n}) = \sum\limits_{k = 1}^{n - 1} \big( 100(x_{i + 1} - x_{i} ^ {2}) ^ {2} + (x_{i} - 1) ^ {2} \big).
\end{equation} 
  
\begin{table}[ht]
    \centering
    \begin{tabular}{c|c|c}
        $n$ & Классический & Модифицированный \\ \hline
        $1$ & $0.20$ с & $0.18$ с \\ \hline
        $2$ & $0.36$ с & $0.26$ с \\ \hline
        $3$ & $7.20$ с & $6.84$ с \\ \hline
        $4$ & $192.43$ с & $214.67$ с
    \end{tabular}
    \caption{Результаты запуска на функции Растригина (\ref{eq:rastrigin})}
    \label{tab:my_label}
\end{table}

Для функции Розенброка исходный брус поиска минимума состоял из интервалов 
вида $[ -10^{2},  10^{2}]$ для каждой координаты, а $\varepsilon = 10^{-4}$. 

\begin{table}[ht]
    \centering
    \begin{tabular}{c|c|c}
        $n$ & Классический & Модифицированный \\ \hline
        $2$ & $1.57$ с & $1.52$ с \\ \hline
        $3$ & $23.10$ с & $23.39$ с \\ \hline
        $4$ & $909.16$ с & $947.36$ с
    \end{tabular}
    \caption{Результаты запуска на функции Розенброка (\ref{eq:rozenbrock})}
    \label{tab:my_label}
\end{table}

\subsection{Выводы}

На практике при оптимизации функций многих переменных к применению рекомендуется 
модифицированный алгоритм, поскольку он не требует дополнительных вычислений значений 
функции, её производных или наклонов, а использует только уже полученную на текущей 
итерации информацию о поведении функции. 

Модификация алгоритма нуждается в дальнейшей доработке, поскольку при оптимизации 
функций от более, чем 3 переменных, дополнительные трудозатраты по вычислению бруса, 
где глобальный минимум не содержится, начинают замедлять вычисления основного алгоритма. 
  
Возможно, имеет смысл применять данную модификацию после нахождения какого-то 
<<предварительного>> значения минимума каким-нибудь классическим <<точечным>> 
алгоритмом.

\clearpage
  
\section{Выводы}
  
Главным результатом работы является построение однопараметрической линейной функциональной арифметики, которая обобщает классическую интервальную арифметику и позволяет проводить 
более точное оценивание областей значений функции. Подробно рассмотрены конструкции 
линейной функциональной арифметики для случаев одной и нескольких переменных. 
  
Реализованные арифметики были использованы для модификации адаптивных интервальных 
алгоритмов глобальной оптимизации функций и алгоритмов глобального решения 
уравнений.  

Модифицированные алгоритмы продемонстрировали высокую скорость сходимости и более точные 
двусторонние оценки минимума функции и решений уравнений в сравнении с их классическими 
версиями. 


\end{document}